\documentclass[a4paper, twoside, headsepline, 11pt, BCOR10mm, DIV10, footinclude, headinclude]{scrbook}

\usepackage{amsmath}  
\usepackage{amsthm}  
\usepackage{color}
\usepackage{ngerman}  
\usepackage{scrpage2}
\usepackage{url}
\usepackage{amssymb}
\usepackage{mathtools} 
\usepackage{tikz}
\usepackage{ifpdf}

\ifpdf
	\usepackage[pdftex,colorlinks,linkcolor=blue,citecolor=blue,urlcolor=blue,a4paper,hypertexnames=true,plainpages=false]{hyperref} 
	\usepackage[all]{xy}   
\SelectTips{cm}{}    
  \hypersetup{pdfauthor=Philipp K"uhl,pdftitle=Isomorphismusvermutungen und 3-Mannigfaltigkeiten}
\else
  \usepackage[all, dvips]{xy}
\SelectTips{cm}{}    
	\usepackage[colorlinks,linkcolor=blue,citecolor=blue,urlcolor=blue,a4paper,hypertexnames=true,plainpages=false]{hyperref}
\fi

	\automark[section]{chapter}

  \ihead{\headmark}
  \ohead[]{\pagemark}
  \addtokomafont{pagehead}{\scshape}
  \addtokomafont{sectioning}{\normalfont\bf}
  \deftripstyle{myheadings}[0pt][0pt]{}{}{\headmark}{}{}{\pagemark}
				     
  \theoremstyle{theorem}
  \newtheorem{satz}{Satz}[chapter]
  \newtheorem{lemma}[satz]{Lemma}
  \newtheorem{korollar}[satz]{Korollar}
  
  \newtheorem*{zielsatz}{Zielsatz}
  \newtheorem*{schlusssatz}{Schlusssatz}
  \newtheorem{etappe}{Etappe}

  \theoremstyle{definition}
  \newtheorem{definition}[satz]{Definition}
  \newtheorem{bemerkung}[satz]{Bemerkung}
  \newtheorem{behauptung}[satz]{Behauptung}

  \newcounter{fallc}
    \newenvironment{fallumgebung}[1]
  {\setcounter{fallc}{0}
  \newcommand{\fall}[1]{\stepcounter{fallc}\noindent {\bf
  \ifcase #1 Fall \Alph{fallc}\or\arabic{fallc}. Fall\fi :}\;{\it ##1}\\ \indent}
  }
  {}

  \newenvironment{beweis}
  {\begin{proof}[Beweis]}
  {\end{proof}}
	\newenvironment{refbeweis}[1]
	{\begin{proof}[Beweis von #1]}
	{\end{proof}}

	\newcommand{\schneekugel}{	
\begin{scope}
	\begin{scope}
				\clip (-1.3,0.9) rectangle (0.8,-0.9);
				\draw (0,0) ellipse (1.2cm and 0.8cm);
			\end{scope}
			\draw (0.8,0.6) .. controls (1,0.55) .. (1.2,0.6);
			\draw (0.8,-0.6) .. controls (1,-0.55) .. (1.2,-0.6);

			\begin{scope}
				\clip (1.2,0.9) rectangle (3.3,-0.9);
				\draw (2,0) ellipse (1.2cm and 0.8cm);
			\end{scope}
			\begin{scope}
					\clip (1,-0.01) rectangle (3,0.7);
					\draw (2,-0.21) ellipse (0.5cm and 0.35cm);
			\end{scope}
			\begin{scope}
					\clip (1,0.09) rectangle (3,-0.8);
					\draw (2,0.59) ellipse (0.8cm and 0.7cm );
			\end{scope}
			\filldraw[white] (3,0) ellipse (1.3cm and 1.3cm);
			\draw (3,0) ellipse (1.3cm and 1.4cm)
						(3,0) ellipse (0.7cm and 0.8cm);	
			\begin{scope}
				  \clip (3,0) ellipse (0.7cm and 0.7cm);	
					\draw (2,0) ellipse (1.2cm and 0.8cm);
					\begin{scope}
						\clip (1,-0.01) rectangle (3,0.7);
						\draw (2,-0.21) ellipse (0.5cm and 0.35cm);
					\end{scope}
			\end{scope}
			\begin{scope}
					\clip (1.5,-0.08) rectangle (2.7,-0.85);
					\filldraw[white] (1,0.5) rectangle (3,-0.9);
					\draw (2,0) ellipse (1.2cm and 0.8cm);
				\end{scope}
				\begin{scope}
					\clip (1,0.09) rectangle (3,-0.8);
					\draw (2,0.59) ellipse (0.8cm and 0.7cm );
				\end{scope}
				\begin{scope}
					\clip (2,0) rectangle (4.4,-1);
					\draw(4,0) ellipse (0.3cm and 0.2cm);
				\end{scope}
				\begin{scope}
					\clip (2,0) rectangle (4.4,1);
					\draw(4,0)[dashed] ellipse (0.3cm and 0.2cm);
				\end{scope}
				\draw (1.5,0) circle (3.2cm);
		\begin{scope}
				\clip (-1.8,-3) rectangle (4.8,0);
				\draw (1.5,0) ellipse (3.2cm and 2cm);
		\end{scope}
		\begin{scope}
				\clip (-1.8,0) rectangle (4.8,3.3);
				\draw[dashed] (1.5,0) ellipse (3.2cm and 2cm);
		\end{scope}	
	\end{scope}
}
\newcommand{\schneekugeln}{
\[
	\pi_1\left(
	\parbox{5.3cm}{
		\begin{tikzpicture}[scale=0.8]
			\schneekugel
					\begin{scope}
					\clip (-1,-0.01) rectangle (1,0.7);
					\draw (0,-0.21) ellipse (0.5cm and 0.35cm);
			\end{scope}
			\begin{scope}
					\clip (-1,0.09) rectangle (1,-0.8);
					\draw (0,0.59) ellipse (0.8cm and 0.7cm );
			\end{scope}
				\begin{scope}
					\clip (-1,0.06) rectangle (1,0.7);
					\draw (0,-0.21) ellipse (0.6cm and 0.45cm); 
					\draw (0,0.4) node {$_{D}$};
				\end{scope}
				\begin{scope}
					\clip (0,-0.21) ellipse (0.6cm and 0.45cm);
					\clip (0,0.59) ellipse (0.8cm and 0.7cm );
					\draw[step=0.8mm, rotate=30] (-1,-1) grid (1,1);
				\end{scope}
			\end{tikzpicture}}
			\right)
\cong 
\mathbb{Z}*\pi_1\left(
\parbox{5.3cm}{
\begin{tikzpicture}[scale=0.8]
				\schneekugel
				\begin{scope}[rotate=-10]
					\begin{scope}
						\node (mitte) at (0,0) {};
						\clip (mitte)+(0,-1) rectangle (2,1);
						\draw[dashed] (mitte) ellipse (0.3cm and 0.8cm);
					\end{scope}
					\begin{scope}
						\clip (mitte)+(0,-1) rectangle (-1,1);
						\draw (mitte) ellipse (0.3cm and 0.8cm);
					\end{scope}
				\end{scope}
			\end{tikzpicture}}
\right)
\]
}

  \makeatletter
  \let\c@equation\c@satz
  \renewcommand{\tagform@}[1]{\maketag@@@{#1}}
  \let\c@figure\c@satz
  \makeatother

	\newcommand{\p}{\pi} 
  \newcommand{\N}{{\mathbb N}}
  \newcommand{\R}{\mathbb R}
  \newcommand{\Z}{{\mathbb Z}}
  
  \newcommand{\A}{{\mathcal A}}  
	\newcommand{\F}{\mathcal F}
	\newcommand{\T}{{\mathcal T}}
	\newcommand{\G}{{\mathcal G}}
	\newcommand{\KT}{\textbf{K}}
	\newcommand{\LT}{{\textbf{L}^{\langle-\infty\rangle}}}
  \newcommand{\Homo}{{\mathcal H}}
	\newcommand{\VC}{{\mathcal VC}}
  \DeclareMathOperator*{\colim}{colim}
  
  \DeclareMathOperator{\id}{id} 
	\DeclareMathOperator{\Aut}{Aut}
	\DeclareMathOperator{\ind}{ind}
	\DeclareMathOperator{\CAT}{CAT}
  \DeclareSymbolFont{symbolsC}{U}{txsyc}{m}{n}
  \DeclareMathSymbol{\sdp}{\mathrel}{symbolsC}{89}
	\newcommand*{\s}[2]{\mathord{^{#1}\mathord{#2}}}
	\newcommand*{\so}[2]{\mathbin{^{#1}\mathord{#2}}}
	\newcommand*{\sd}[2]{\,\mathbin{^{#1}\!\mathord{#2}}}

	\DeclareFontFamily{U}{matha}{\hyphenchar\font45}
	\DeclareFontShape{U}{matha}{m}{n}{
      <5> <6> <7> <8> <9> <10> gen * matha
      <10.95> matha10 <12> <14.4> <17.28> <20.74> <24.88> matha12
      }{}
	\DeclareSymbolFont{matha}{U}{matha}{m}{n}
	\DeclareFontSubstitution{U}{matha}{m}{n}

	\DeclareMathSymbol{\nsimeq} {3}{matha}{"1C}

  \newcommand{\invers}{^{-1}}
  \newcommand{\gexseq}[6]{\[1\rightarrow #1\xrightarrow{#4} #2\xrightarrow{#5} #3\rightarrow 1#6\]}
  \newcommand{\exseq}[4]{\[1\rightarrow #1\rightarrow #2\rightarrow #3\rightarrow 1#4\]}
   
	\newcommand{\E}{E_{\F} G}
  \newcommand{\pt}{\mbox{pt}}

 	\newcommand{\rem}[1]{}
	\newcommand{\pushout}[8]{\xymatrix{{#1}\ar@{->}[r]^{#5}\ar@{->}[d]_{#6}&{#2}\ar@{->}[d]^{#7}\\{#3}\ar@{->}[r]^{#8}&{#4}}}

  \author{Philipp K"uhl}
  \title{Isomorphismusvermutungen und 3"=Mannigfaltigkeiten}
  \date{\today}

\begin{document}
\frontmatter
\pagestyle{myheadings}


\begin{titlepage}
	\setlength{\topmargin}{0cm}
\thispagestyle{empty}
\enlargethispage{40mm}
\begin{center}

	{
		\sc  Westf"alische Wilhelms-Universit"at M"unster\\
      \vspace*{1mm} Fachbereich Mathematik und Informatik}
			\vspace*{10em}

		{\Huge Isomorphismusvermutungen\\ \Huge und\\ \Huge 3"=Mannigfaltigkeiten\\}
    \vspace*{6em}
    {\Large  Diplomarbeit}\\
    \vspace{3mm}
    \large vorgelegt von\\
    \vspace{4mm}
    {\Large Philipp K"uhl}\\
    \vspace{10mm}
    {\large September 2008}\\
    \vspace{65mm}
\end{center}
	\begin{flushright}
		\begin{tabular}{rr}
      \vspace{1mm}
			Betreuer: & Prof. Dr. Wolfgang L"uck\\
      Zweitgutachter: & Prof. Dr. Arthur Bartels
	\end{tabular}\\
	\end{flushright}
\end{titlepage}
\clearpage{\pagestyle{empty}\cleardoublepage}


  \thispagestyle{empty}
	\enlargethispage{20mm}
  \hskip 0mm
  \vfill
  \rem{\begin{center}\sffamily\bfseries\Large
    Erkl"arung zur Selbstst"andigkeit
	\end{center}}
  \bigskip\noindent Ich erkl"are hiermit, dass ich die vorliegende Arbeit
  selbstst"andig verfasst und keine anderen als die angegebenen Quellen und Hilfsmittel verwendet habe.\par
  \bigskip\noindent M"unster, \today\par
  \vskip 10mm
  \hfill\hrulefill


\chapter*{Einleitung}
Roushon hat in \cite{roushon} und \cite{roushon-b-gruppen} eine umfassende Sammlung von Resultaten zusammengestellt, die sich mit
Fundamentalgruppen von dreidimensionalen Mannigfaltigkeiten und der Farrell-Jones-Isomorphismusvermutung besch"aftigen. Sein Hauptaugenmerk lag dabei auf der Farrell-Jones-Isomorphismusvermutung bez"uglich des stabilen topologischen Pseudo-Iso\-topie-Funktors. Wir wollen uns dieser Thematik abweichend von der urspr"unglichen Definition von Farrell-Jones \cite{farrell-jones} mit einer allgemeineren Formulierung mit Hilfe von "aquivarianten Homologietheorien n"ahern, die beispielsweise auch die Baum-Connes-Vermutung mit einschlie"st. Eine ausf"uhrliche Einf"uhrung in diese Konstruktion findet sich bei L"uck und Reich \cite{lueck-reich}. Wir werden mit einer erweiterten Form der dort definierten gefaserten Isomor\-phismus\-vermutung arbeiten, die bez"uglich endlichen Erweiterungen gute Vererbungseigenschaften hat. 

Die vorliegende Arbeit kristallisiert aus den Ergebnissen von Roushon die Voraussetzungen heraus, die solch eine beliebige  "`endlich erweiterbare, gefaserte Isomorph\-ismus\-vermutung"' erf"ullen muss, damit sie f"ur Fundamentalgruppen von beliebigen dreidimensionalen Mannigfaltigkeiten gilt.  
Ein besonderes Augenmerk wird dabei auf der Farrell-Jones-Isomorphismusvermutung in L-Theorie und algebraischer K-Theorie liegen, da sie bereits einige Ergebnisse mitbringt (siehe beispielsweise \cite{lueckmain}), die die n"otigen Voraussetzungen daf"ur erheblich einschr"anken.

Das erste Kapitel, "`Das Fundament"' auf dem diese Arbeit ruht, besteht aus den ben"otigten Grundlagen "uber dreidimensionale Mannigfaltigkeiten und Kranzprodukten von Gruppen.
F"ur die Klassifikation von dreidimensionalen Mannigfaltigkeiten spielen zwei Zerlegungen eine wichtige Rolle: Die Zerlegung entlang von Sph"aren, die als Primzerlegung bekannt ist und auf Kneser \cite{kneser} zur"uckgeht, und die Zerlegung entlang von Tori, die Jaco-Shalen \cite{jaco-shalen} und Johannson \cite{johannson} entwickelt haben. Wir geben im ersten Teil des Kapitels einen kurzen "Uberblick "uber diese beiden Zerlegungen und gehen im zweiten Teil auf die grundlegenden Eigenschaften von Kranzprodukten von Gruppen ein, die uns bei dem Umgang mit endlichen Erweiterungen bzw. endlichen "Uberlagerungen n"utzliche Dienste leisten werden.

Das zweite Kapitel f"uhrt kurz und knapp die Definitionen der Isomorphismusvermutungen ein, zu denen im dritten Kapitel, dem "`FIC-Werkzeugkasten"' dann alle Aussagen zusammengestellt werden, die wir "uber Isomophismusvermutungen brauchen werden. Die "`Schraubenzieher"' sind dabei relativ grundlegende, allgemeine Ergebnisse, w"ahrend wir mit den "`Bohrmaschinen"' f"unf Aussagen aufstellen, deren G"ultigkeit wir von einer Isomorphismusvermutung fordern m"ussen, damit sie f"ur Fundamentalgruppen von dreidimensionalen Mannigfaltigkeiten gilt. Sie sind mit W1 bis W5 durchnummeriert und k"onnen abh"angig von der betrachteten Isomorphismusvermutung noch vereinfacht oder teilweise sogar bewiesen werden. F"ur die Farrell-Jones-Isomorphismusvermutung in L-Theorie und algebraischer K-Theorie bleibt dabei vor allem W1 - die G"ultigkeit f"ur Gruppen der Gestalt $\Z^2\sdp \Z$ - als entscheidende Voraussetzung bzw. Hindernis bestehen. 

Wenn wir W1 bis W5 als gegeben voraussetzen, k"onnen wir die Werkzeuge W3 bis W5 ohne Einschr"ankung zu allgemeineren Aussagen versch"arfen, mit *W3 bis *W5 bezeichnet, und au"serdem vier weitere wichtige "`Bohrmaschinen"' - *W6 bis *W9 - beweisen. Im Anhang sind zur "Ubersicht noch einmal alle "`Bohrmaschinen"' aufgelistet.

Als letztes Kapitel schlie"slich f"uhrt "`Der rote Baum"' als verzweigter roter Faden in mehreren Etappen durch den Beweis unseres Zielsatzes:
\begin{zielsatz}
	Eine endlich erweiterbare, gefaserte Isomorphismusvermutung, die die Anforderungen W1 bis W5 erf"ullt, gilt f"ur Fundamentalgruppen von dreidimensionalen Mannigfaltigkeiten.
\end{zielsatz}
Die grundlegende Voraussetzung daf"ur, die mit W2 gefordert wird, ist, dass die Isomorphismusvermutung f"ur nicht-positiv gekr"ummte Mannigfaltigkeiten gilt. Das verhilft uns zu den zwei wichtigsten Resultaten, auf denen alles weitere aufbaut: Die G"ultigkeit f"ur Fundamentalgruppen von geschlossenen Seifert-Mannigfaltigkeiten und von geschlossenen Haken-Mannigfaltigkeiten mit atorischer Komponente. Im "`Mit-Rand-Ast"'-Kapitel werden wir Mannigfaltigkeiten mit Rand auf diese beiden Resultate zur"uck\-f"uhr\-en. Vor allem bei Mannigfaltigkeiten mit Randkomponenten vom Geschlecht $\geq 2$ stellen sich dabei etwas gr"o"sere Schwierigkeiten in den Weg, "uber die auch Roushon ins Stolpern gekommen ist. Wir werden im Abschnitt "`B-Gruppen"' neben Roushons zweifelhaftem noch einen alternativen Beweis f"ur diese Klasse von Mannigfaltigkeiten ausf"uhren.
Im letzten Abschnitt bleiben schlie"slich nur noch geschlossene Graph-Mannigfaltigkeiten "ubrig, f"ur deren Fundamentalgruppen mit Hilfe des Hauptsatzes aus \cite{roushon} die Isomorphismusvermutung gezeigt und damit der Beweis vollendet wird.

\rem{
In den ersten beiden Kapitel Grundlagen "uber Zerlegung von Mannigfaltigkeiten, Kranzprodukte und Isomorphismusvermutungen.\\
FIC-Werkzeugkasten: \\
Schraubenzieher - grundlegende Eigenschaften\\
Bohrmaschinen - n"otige Voraussetzungen\\
Dann in "`Etappen"' durch den roten Baum "`klettern"'
}
\tableofcontents
\deftripstyle{myheadings}[0pt][0pt]{}{}{}{}{}{\pagemark}
\clearpage{\pagestyle{myheadings}\cleardoublepage}

\mainmatter
	
\chapter{Das Fundament}

Wir beginnen mit der Kl"arung von Bezeichnungen und der Einf"uhrung von s"amtlichen grundlegenden Eigenschaften von dreidimensionalen Mannigfaltigkeiten und Kranzprodukten von Gruppen, die wir im Folgenden brauchen werden.
\section{Zerlegung von 3"=Mannigfaltigkeiten}

Sei  $D^d\coloneqq\{x\in\R^d\; |\; \sum_{i=1}^{d} x_i^2\leq 1\}$ die \emph{d-dimensionale Einheitskugel} bzw. \emph{d-Disk}, $S^d\coloneqq\{x\in\R^{d+1}\; | \; \sum_{i=1}^{d+1}x_i^2=1\}$ die \emph{d-dimensionale Sph"are} und $T^d\coloneqq\Pi_{i=1}^dS^1$ der \emph{d-dimensionale Torus}. $S^1\times [0,1]$ bezeichnen wir als \emph{Annulus}.

\emph{3"=Mannigfaltigkeit} wird als Abk"urzung f"ur dreidimensionale Mannigfaltigkeit verwendet. Wenn im Folgenden von Mannigfaltigkeitsgruppen die Rede ist, sind immer die Fundamentalgruppen der Mannigfaltigkeiten gemeint. Des Weiteren verstehen wir in der gesamten Arbeit unter der Bezeichnung
\emph{Mannigfaltigkeit mit Rand} oder \emph{berandete  Mannigfaltigkeit}  den Spezialfall von Mannigfaltigkeiten mit leerem Rand als ausgeschlossen. Eine \emph{geschlossene Mannigfaltigkeit} ist eine kompakte Mannigfaltigkeit ohne Rand.

Wir werden h"aufiger eine berandete Mannigfaltigkeit mit einer Kopie am Rand identifizieren, um eine geschlossene Mannigfaltigkeit zu erhalten. Wir nennen diesen Vorgang \emph{Verdoppeln} und schreiben f"ur die Verdopplung $\widetilde{M}$ einer Mannigfaltigkeit $M$
\[
\widetilde{M} = M\cup_\partial M
\]
bzw. wenn wir auf einer Randfl"ache nicht mit der Identit"at sondern einem Hom"oo\-mor\-phis\-mus $f$ verkleben
\[
\widetilde{M} = M\cup_f M.
\]

F"ur die Klassifikation von 3"=Mannigfaltigkeiten hat man in den vergangenen 80 Jahren Zerlegungen entlang von Sph"aren und Tori entwickelt und studiert. Wir werden beide Zerlegungen ben"otigen und deshalb im Folgenden die Grundz"uge und notwendigen Begriffe kurz zusammenfassen.

Seien $M, M_1, M_2$ zusammenh"angende 3"=Mannigfaltigkeiten mit 3-Disks $D_i\subset \mathring M_i$ und Einbettungen $h_i:M_i\smallsetminus\mathring D_i\rightarrow M$ f"ur $i=1,2$, so dass 
\[
h_1(M_1\smallsetminus\mathring D_1)\cap h_2(M_2\smallsetminus\mathring D_2) = h_1(\partial D_1) = h_2(\partial D_2)
\] und 
\[M=h_1(M_1\smallsetminus\mathring D_1)\cup h_2(M_2\smallsetminus\mathring D_2)\]
gilt.
Wir nennen $M$ dann die zusammenh"angende Summe von $M_1$ und $M_2$ und schreiben daf"ur $M=M_1\#M_2$.
Sie ist eindeutig bis auf Hom"oomorphie, wenn $M$, $M_1$ und $M_2$ orientiert sind und wir fordern, dass $h_1$ und $h_2$ orientierungserhaltend sind.

\begin{definition}[prim]
	Eine kompakte 3"=Mannigfaltigkeit $M$ hei"st \emph{prim}, wenn f"ur jede zusammenh"angende Summe $M\cong M_1 \# M_2$ gilt $M_1\cong S^3$ oder $M_2\cong S^3$.
\end{definition}

Die Zerlegung von 3"=Mannigfaltigkeiten entlang von Sph"aren in Prim-Mannig\-faltig\-keiten geht auf Kneser zur"uck \cite[S. 253]{kneser}. Eine "ubersichtliche Darstellung findet sich bei Hempel \cite[Kapitel 3]{hempel}. Wir werden die Primzerlegung nur f"ur kompakte, orientierbare Mannigfaltigkeiten ben"otigen, aber f"ur die Eindeutigkeit der Zerlegung braucht die Orientierbarkeit nicht gefordert zu werden. 

\begin{satz}[Primzerlegungssatz]\label{thm:primzerlegung}
	Sei $M$ eine kompakte 3"=Mannigfaltigkeit. Dann gibt es endlich viele Prim"=Mannigfaltigkeiten $M_1,\ldots,M_n$ mit
	$M=M_1\#\ldots\#M_n$. Die $M_i$ sind bis auf Hom"oomorphie und Reihenfolge der Summanden eindeutig.
\end{satz}

F"ur einen Beweis siehe beispielsweise \cite[Theorem 3.15 und Theorem 3.21]{hempel}.
Mit folgendem Begriff k"onnen Prim"=Mannigfaltigkeiten noch genauer klassifiziert werden.

\begin{definition}[irreduzibel]
  Eine 3"=Mannigfaltigkeit $M$ hei"st \emph{irreduzibel}, wenn jede $S^2$ in $M$ eine $D^3$ berandet, das hei"st zu jeder Einbettung $f:S^2\rightarrow M$ gibt es eine Einbettung $\overline{f}:D^3\rightarrow M$ mit $\partial\overline{f}(D^3)=f(S^2)$. 
\end{definition}

\begin{lemma}\label{lm:prim-irr}
	Sei $M$ prim. Dann ist $M$ irreduzibel oder ein $S^2$-B"undel "uber $S^1$.
\end{lemma}
Auch dazu findet sich ein Beweis bei Hempel \cite[Lemma 3.13]{hempel}.
\\

F"ur die Zerlegung entlang von Tori m"ussen wir etwas weiter ausholen und zun"achst den Begriff der Seifert-Mannigfaltigkeit einf"uhren. Er geht zur"uck auf eine Arbeit von Seifert \cite{seifertraum}. Wir werden uns bei Einf"uhrung der Definition vor allem an das Buch von Jaco \cite{jaco} halten, das sich auf den f"ur uns ausreichenden orientierbaren Fall beschr"ankt, und dabei die deutschen Bezeichnungen von Orlik, Vogt und Zieschang \cite{orlik} verwenden. 

Sei die 2-Disk $D^2\subset\R^2$ durch Polarkoordinaten $\{(r,\theta)\;|\;0\leq r\leq 1,\, 0\leq\theta<2\pi\}$ gegeben und $(\mu,\nu)$ ein Paar ganzer, teilerfremder Zahlen. Ein \emph{gefaserter Volltorus vom Typ $(\mu,\nu)$} ist ein Vollzylinder $D^2\times[0,1]$, bei dem die Boden- und Deckelfl"ache um den Winkel $\mu/\nu$ verdreht identifiziert sind, also \[((r,\theta), 1) = ((r, \theta + \frac{2 \pi\nu}{\mu}), 0)\] gesetzt wird. Dabei setzen sich jeweils  $\mu$ achsenparallele Strecken zu einfach-geschlos\-senen Kurven, den Fasern 
\rem{\[F\coloneqq\left.\left\{\left(r,\theta+k\frac{2\pi\nu}{\mu}\right)\;\right|\quad k=0,\ldots,\mu -1\right\}\times\left[0,1\right]$\]}
des Volltorus zusammen, die das Loch des Volltorus $\mu$ Mal umlaufen und sich dabei $\nu$ Mal um den Torus winden. Einzige Ausnahme ist die Achse des Vollzylinders, die Seele, die den gefaserten Volltorus nur einmal uml"auft. Bis auf fasertreue Hom"oomorphismen kann $\mu> 0$ und $0\leq\nu< \mu$ angenommen werden. Ist $\mu>1$ hei"st die Seele Ausnahmefaser, ist $\mu = 1$, so hei"st der gefaserte Volltorus regul"ar gefasert und jede Faser ist eine regul"are Faser.
\begin{definition}
  Eine zusammenh"angende orientierbare kompakte 3"=Mannigfaltigkeit $M$ hei"st \emph{Seifert-Mannigfaltigkeit}, wenn es eine Zerlegung von $M$ in disjunkte einfach-geschlossene Kurven, den sogenannten Fasern, gibt, so dass jede Faser eine geschlossene Umgebung von Fasern besitzt, die fasertreu hom"oomorph zu einem gefaserten Volltorus ist. R"ander einer Seifert-Mannigfaltigkeit sind immer gefaserte Tori.

  Der Quotientenraum, den man aus einer Seifert-Mannigfaltigkeit $M$ erh"alt, wenn man jede Faser mit einem Punkt identifiziert, hei"st \emph{Zerlegungsfl"ache} von $M$ und ist eine zusammenh"angende zweidimensionale Mannigfaltigkeit. Auch wenn die Seifert-Mannigfaltigkeit orientierbar war, muss die Zerlegungsfl"ache nicht orientierbar sein.
\end{definition}

Wir werden f"ur die Fundamentalgruppen von Seifert-Mannigfaltigkeiten die Pr"a\-sen\-ta\-tionen aus dem folgenden Satz verwenden, der
beispielsweise von Hempel \cite[Satz 12.1]{hempel} und ausf"uhrlicher von Seifert \cite[\S10]{seifertraum} bewiesen wird.
\begin{satz}\label{thm:seifert-presentation}
  Sei $M$ eine kompakte Seifert-Mannig\-faltig\-keit mit $q$ Ausnahmefasern und $B$ ihre Zerlegungsfl"ache vom Geschlecht $g$ und mit $p$ Randkomponenten. 
  Ist $B$ orientierbar, gilt 
	\begin{eqnarray*}
  	\pi_1(M) &= \quad\left< \right. &a_1,b_1,\ldots,a_g,b_g,c_1,\ldots,c_q, d_1,\ldots,d_p, t; \\
		& & a_ita_i\invers =t^{\varepsilon_i}, b_itb_i\invers = t^{\delta_i}, c_jtc_j\invers = t^{\eta_j}, d_ktd_k\invers = t^{\theta_k}, c_j^{\alpha_j} = t^{\beta_j},\\
		& & \left. c_q=[a_1,b_1]\ldots [a_g,b_g]c_1\ldots c_{q-1} d_1\ldots d_p\quad\right>
	\end{eqnarray*}
\mbox{  bzw. ist $B$ nicht orientierbar, gilt}
\begin{eqnarray*}
 \pi_1(M) &= \quad\left< \right. &a_1,\ldots,a_g,c_1,\ldots,c_q, d_1,\ldots,d_p, t; \\
 & &a_ita_i\invers =t^{\varepsilon_i}, c_jtc_j\invers = t^{\delta_j}, d_ktd_k\invers = t^{\eta_k}, c_j^{n_j} = t^{s_j},\hspace{5em}\\
 & & \left. c_q=a_1^2\ldots a_g^2c_1\ldots c_{q-1} d_1\ldots d_p\quad\right>,
  \end{eqnarray*}
	wobei alle $\varepsilon_i,\delta_i,\eta_i,\theta_i\in\{-1,+1\}$ und f"ur jede Ausnahmefaser $j$ gilt $0<n_j<s_j$.
\end{satz}

\rem{ TODO raus? TODO
\begin{definition}[asph"arisch]
	Eine 3"=Mannigfaltigkeit $M$ hei"st \emph{asph"arisch}, wenn $M$ zusammenh"angend ist und die universelle "Uberlagerung kontraktibel ist.
\end{definition}}
Eine \emph{Fl"ache} ist im Folgenden immer eine zusammenh"angende, kompakte, orientierbare, zweidimensionale Mannigfaltigkeit. Die Fundamentalgruppen von Fl"achen (vom Geschlecht $\geq x$) werden abk"urzend auch als \emph{Fl"achengruppen} (vom Geschlecht $\geq x$) bezeichnet.

\begin{definition}[unkomprimierbar]
	Sei $M$ eine 3"=Mannigfaltigkeit. Eine in $M$ eingebette Fl"ache $F$ mit $\partial M\cap F = \partial F$ oder $F\subset \partial M$ hei"st \emph{unkomprimierbar}, wenn $M$ keine 2-Sph"are ist und die Einbettung $F\rightarrow M$ eine injektive Abbildung $\pi_1(F)\rightarrow\pi_1(M)$ auf den Fundamentalgruppen induziert. Andernfalls hei"st sie \emph{komprimierbar}.
\end{definition}
\begin{definition}[zweiseitig] Sei $h:F\rightarrow M$ eine Einbettung einer Fl"ache $F$ in eine 3"=Mannigfaltigkeit $M$. $F$ hei"st \emph{zweiseitig} in $M$ eingebettet, wenn es eine Einbettung $\widetilde{h}:F\times [0,1]\rightarrow M$ gibt, so dass  $\widetilde{h}(x,\frac{1}{2})=h(x)$ f"ur alle $x\in F$, $\widetilde{h}(\partial F \times [0,1]) \subset\partial M$ und $\widetilde{h}(F\times [0,1])$ eine Umgebung von $F$ in $M$ ist. Dies ist genau dann der Fall, wenn das Normalenb"undel von $F$ trivial ist.
\end{definition}
\begin{definition}[Haken-Mannigfaltigkeit]
  Eine irreduzible 3"=Mannigfaltigkeit M hei"st \emph{Haken-Mannigfaltigkeit}\footnote{In der englischen Literatur wird h"aufig auch der Begriff \emph{sufficiently large manifold} benutzt.}, wenn sie eine unkomprimierbare, zweiseitige Fl"ache enth"alt.
\end{definition}
\begin{definition}[parallel]
  Seien $F_0$ und $F_1$ zwei in einer 3-Mannigfaltigkeit $M$ zweiseitig eingebettete Fl"achen. $F_0$ und $F_1$ hei"sen \emph{parallel}, wenn es eine Einbettung $h: F_0 \times [0,1] \rightarrow M$ gibt, so dass  $h(F_0 \times \{i\}) = F_i$ f"ur $i=0,1$. Eine zweiseitig ein\-gebette Fl"ache $F\subset M$ hei"st \emph{randparallel}, wenn es eine Einbettung $h: F\times [0,1] \rightarrow M$ gibt mit $h(F\times {0}) = F$ und $\partial (h(F\times [0,1]))=F$. 
\end{definition}

\begin{definition}[atorisch]
  	Eine 3"=Mannigfaltigkeit $M$ hei"st \emph{atorisch}, wenn jeder unkomprimierbar eingebetteter Torus bereits randparallel ist.
		
\end{definition}

Sei $M$ eine kompakte orientierbare Haken-Mannigfaltigkeit. Wir wollen $M$ entlang von disjunkt eingebetteten, unkomprimierbaren Tori auftrennen.
Diese Zerlegung ist im Allgemeinen erst einmal nicht eindeutig. Johannson \cite{johannson} und Jaco-Shalen \cite{jaco-shalen} haben aber gezeigt, dass lediglich Seifert-Mannigfaltigkeiten f"ur die Nichteindeutigkeit verantwortlich sind. Man kann deshalb $M$ eindeutig bis auf Isotopie in eine minimale Menge von Seifert- und atorischen Mannigfaltigkeiten zerlegen. Diese Zerlegung wollen wir \emph{Toruszerlegung} nennen. Da es unter den $|SL_2(\Z)|$-vielen M"oglichkeiten, einen Volltorus einzukleben, keine kanonische Wahl gibt, behalten die Teilst"ucke von der Toruszerlegung unkomprimierbare Torusr"ander zur"uck. Wir werden diese Teilst"ucke im Folgenden als \emph{Seifert- und atorische Komponenten} (der Toruszerlegung von $M$) bezeichnen. Es gibt ein paar Mannigfaltigkeiten, die beides sind, die wir aber den Seifertkomponenten zuordnen wollen. Wenn wir also im weiteren Verlauf der Arbeit von einer atorischen Komponente sprechen, soll der Fall einer atorischen Seifert-Mannigfaltigkeit immer ausgeschlossen sein.
Alle Seifert-Komponenten zusammen hei"sen \emph{charakteristische Untermannigfaltigkeit} von $M$. Eine Mannigfaltigkeit, deren Toruszerlegung nur aus Seifert-Komponenten besteht, hei"st \emph{Graph-Mannig\-faltigkeit} und wir nennen sie \emph{nicht-triviale Graph-Mannigfaltigkeit}, wenn die Toruszerlegung nicht trivial ist, also aus mindestens zwei Seifert-Komponente besteht.
Besitzt $M$ einen nicht-leeren Rand, kann man zus"atzlich noch entlang unkomprimierbarer Annuli zerlegen. Wir werden aber haupts"achlich nur die Zerlegung entlang von Tori verwenden und fassen das abschlie"send in einem Satz zusammen. Neumann und Swarup haben in \cite{neumann-swarup} daf"ur einen neuen "ubersichtlicheren Beweis gef"uhrt.

\begin{satz}[Toruszerlegung]\label{thm:toruszerlegung}
	Sei $M$ eine kompakte orientierbare Haken"=Mannigfaltigkeit. Dann gibt es in $M$ eine Menge $T$ von unkomprimierbaren, disjunkten Tori, so dass jede Komponente  von $M\smallsetminus T$ entweder atorisch oder seifertsch ist. So eine minimale Menge $T$ ist eindeutig bis auf Isotopie.
\end{satz}

\section{Kranzprodukt}
Wir werden viel mit Kranzprodukten und auch semidirekten Produkten zu tun haben und neben allgemein bekannten auch ein paar spezielle Eigenschaften ben"otigen, die wir zusammen mit der Definition im Folgenden einf"uhren.

    F"ur eine Gruppe $G$ und eine Menge $M$ bezeichnet $G^M$ die Gruppe der Abbildungen von $M$ nach $G$. Die Verkn"upfung ist durch punktweise Verkn"upfung in $G$ gegeben:
    \[
        (fg)(m) = f(m)g(m) \qquad\text{f"ur $f,g\in G^M$ und $m\in M$.}
    \]
		Operiert $G$ auf $M$ schreiben wir f"ur die Operation eines Elements $g\in G$ auf $m\in M$ kurz $\s{g}{m}$.

\begin{definition}[Semidirektes Produkt]
  Seien $A$ und $Q$ beliebige Gruppen und 
  	\[
	\phi: Q\rightarrow \Aut(A)
	\]
  ein Gruppenhomomorphismus. $Q$ operiert also auf $A$ und mit Hilfe dieser Operation konstruieren wir eine Erweiterung $E$ von $A$ durch $Q$.
  Die Menge 
  \[
  	E=A\times Q = \{(a,q)\;|\;a\in A, q\in Q\}
  \]
  wird mit der Verkn"upfung
  \[
	(a,q)(a',q') = (a\sd{q}a',qq') = (a\phi(q)(a'),qq')
  \]
  zu einer Gruppe.
  Mit
	\[	\alpha: A  \rightarrow   E,\quad a  \mapsto (a,1)\quad\mbox{ und }\quad	\beta:	E  \rightarrow Q, \quad(a,q)  \mapsto  q
	\]
  sieht man, dass $E$ eine Erweiterung von $A$ durch $Q$ ist. $E$ hei"st \emph{semidirektes Produkt} von $A$ und $Q$ und wir schreiben daf"ur 
  \[
  	E = A \sdp_\phi Q
  \]
  bzw., wenn die verwendete Gruppenoperation klar oder irrelavant ist,
  \[
  	E = A \sdp Q.
  \]
\end{definition}
\begin{lemma}\label{lm:semidirekt-spaltende-sequenz}
	Sei 
	\gexseq{K}{G}{H}{}{p}{}
	eine exakte Sequenz von Gruppen, die spaltet. Das hei"st es gibt einen 
	Gruppenhomomorphismus $s:H\rightarrow G$ mit $p\circ s=\id_H$.
	Dann ist 
	\[
		G\cong K \sdp_\phi H,
	\]
	wobei
	\begin{eqnarray*}
			\phi: H  & \rightarrow & \Aut(K) \\
						h & \mapsto & (k \mapsto s(h)ks(h)\invers).
 	\end{eqnarray*}
\end{lemma}
\begin{beweis}
	Wir fassen $K$ als Untergruppe von $G$ auf und definieren 
	\[
		\Phi: K\sdp_\phi H \rightarrow G, \quad(k,h)\mapsto ks(h).
	\]
	Wie man leicht nachrechnet, ist $\Phi$ ein Isomorphismus: $\Phi$ ist surjektiv, denn f"ur $g\in G$ ist $gs(p(g)\invers)\in\ker(p)\cong K$ und es gilt $\Phi(gs(p(g)\invers),p(g))=g$.\\
	$\Phi$ ist injektiv, denn
	\begin{eqnarray*}
		1 = \Phi(k,h)  = ks(h) & \Rightarrow & s(h)=k\invers \Rightarrow s(h)\in\ker(p)\\
		& \Rightarrow & 1=p(s(h))=h\Rightarrow k=1.
	\end{eqnarray*}
	$\Phi$ ist ein Homomorphismus, denn 
	\begin{eqnarray*}
		\Phi(k,h)\Phi(k',h') & = & ks(h)k's(h')= ks(h)k's(h)\invers s(h) s(h')\\ & = & k\phi(h)(k')s(h)s(h') = \Phi(k\sd{h}{k'},hh').
	\end{eqnarray*}
\end{beweis}
\begin{definition}[Kranzprodukt]
    Seien $G$ und $H$ Gruppen und $H$ operiere auf einer Menge $M$. Die Operation von $H$ auf $M$ induziert eine Operation von $H$ auf $G^{M}$ durch
    \[
		\s{h}{\!f(m)} \coloneqq f(\sd{h\invers}{m}) \qquad \mbox{f"ur } f\in G^M, h\in H, m\in M.
	\]
	Das \emph{Kranzprodukt} $G\wr_M H$ ist als das semidirekte Produkt von $G^M$ und $H$ bez"uglich dieser induzierten Gruppenoperation definiert.
	Also 
 	\[
		G\wr_M H \coloneqq G^M \sdp H,
	\]
	wobei $G\wr_M H$ als Menge durch $G^M\times H$ gegeben ist und mit der Verkn"upfung 
    \[
		(f,h)(f',h') = (f\;\s{h}{\!f'},hh')
	\] 
	zu einer Gruppe wird:

	Ein neutrales Element ist durch $(c_1,1)\in G\wr_M H$ mit $c_1(m) = 1$ f"ur alle $m\in M$ ge\-geben und zu jedem  $(f,h)\in G\wr H$ findet man in $(\so{h^{-1}}{\!\!f\invers},h\invers)$ ein Inverses, wobei $\s{h\invers}{\!\!f\invers(x)}\coloneqq(f(\sd{h\invers}{x}))\invers$ f"ur $x\in M$. Die Assoziativit"at rechnet man wie folgt nach:	
	\begin{eqnarray*}
		(f_1,h_1)((f_2,h_2)(f_3,h_3)) & = & (f_1,h_1)(f_2\,(\s{h_2}{\!f_3}),h_2h_3) \\
		& = & (f_1\,\s{h_1}{\!f_2}\,\s{h_1h_2}{\!f_3},h_1h_2h_3) \\
		& = & (f_1\sd{h_1}{f_2},h_1h_2)(f_3,h_3)\\
									& = & (((f_1,h_1)(f_2,h_2))(f_3,h_3).
	\end{eqnarray*}

	Wir werden bei Kranzprodukten immer die Operation von $H$ auf sich selbst durch Linksmultiplikation verwenden und schreiben f"ur $G\wr_H H$ kurz $G\wr H$.
  \end{definition}

\begin{bemerkung}
	Mit dieser Definition ist $G^{H}\times\{1\}$ offensichtlich ein Normalteiler von $G\wr H$, der endlichen Index hat, wenn $H$ endlich ist. 
\end{bemerkung}

Wir werden sp"ater viel mit endlichen Erweiterungen zu tun haben. Dabei leistet das Kranzprodukt hilfreiche Dienste, 
da man f"ur Erweiterungen
\exseq{N}{G}{K}{}
$G$ als Untergruppe in $N\wr K$ wiederfindet:
\begin{satz}\label{thm:Kranzprodukt-Erweiterungen} 
    Sei $G$ eine beliebige Gruppe mit einer normalen Untergruppe $N$ und sei $K\coloneqq G/N$. Dann ist $G$ eine Untergruppe von $N \wr K$.
\end{satz}
\begin{beweis}
  Sei $p:G\rightarrow K$ ein surjektiver Homomorphismus mit Kern $N$. 
  F"ur jedes $u\in K$ w"ahle aus $p\invers(u)$ einen Repr"asentanten $t_u\in G$.
  Dann liegt f"ur jedes $x\in G$ $t_uxt\invers_{up(x)}$ im Kern von $p$ und  wir k"onnen 
  Abbildungen $f_x:K\rightarrow N$ definieren durch
  \[
  f_x(u) \coloneqq t_uxt\invers_{up(x)} \quad\mbox{ f"ur alle } u\in K
  \]
  und erhalten mit
  \begin{eqnarray*}	
	\Phi: G&\rightarrow& N\wr K\\
		  x&\mapsto& (f_x,p(x))
  \end{eqnarray*}
  einen injektiven Gruppenhomomorphismus. Im Detail nachgerechnet wird das in \cite[2.6A Universal embedding theorem]{dixon}.
\end{beweis}

Wir werden noch einige weitere Eigenschaften des Kranzproduktes ben"otigen, die wir im Folgenden auff"uhren.
\begin{lemma}\label{lm:Kranzprodukt-Untergruppen}
  Sei $G$ eine Gruppe, $U$ eine Untergruppe von $G$ und $H$ eine beliebige endliche Gruppe. Dann ist $U\wr H$ eine Untergruppe von $G\wr H$.
\end{lemma}
\begin{beweis}
  	Sei $f:U\rightarrow G$ ein injektiver Gruppenhomomorphismus. Dann ist auch 
	\begin{eqnarray*}
	  \widetilde{f}:U\wr H & \rightarrow & G\wr H \\
	  (u,h) & \mapsto & (f\circ u,h)
	\end{eqnarray*}
	ein injektiver Gruppenhomomorphismus, wie man leicht nachrechnet.
\end{beweis}

\rem{
TODO raus? wird glaube ich nicht verwendet, ist das "uberhaupt richtig? TODO
Das gleiche gilt auch in der zweiten Komponente des Kranzproduktes:
\begin{lemma}\label{lm:Kranzprodukt-Untergruppen-zweite-Komponente}
	Sei $G$ eine Gruppe, $H$ eine beliebige endliche Gruppe und $U$ eine Untergruppe von $H$. Dann ist $G\wr U$ eine Untergruppe von $G\wr H$.
\end{lemma}
\begin{beweis}
	Sei $f:U\rightarrow H$ ein injektiver Gruppenhomomorphismus. Dann ist durch
	\begin{eqnarray*}
	  \overline{f}:G\wr U & \rightarrow & G\wr H \\
	  (\widetilde{g},u) & \mapsto & (\widetilde{g}\circ f\invers,f(h))
	\end{eqnarray*}
	ein injektiver Gruppenhomomorphismus definiert. 
\end{beweis}}
\begin{lemma}\label{lm:Kranzprodukt-Kreuzprodukt}
	Seien $G_1,G_2$ und $H$ Gruppen. Dann ist $(G_1\times G_2)\wr H$ eine Untergruppe von $(G_1\wr H)\times(G_2\wr H)$.
\end{lemma}
\begin{beweis}
	Sei $(\phi, h)\in (G_1\times G_2)\wr H$, also $\phi: H\rightarrow G_1\times G_2$ eine Abbildung von Mengen und $h\in H$.
	$p_i$ bezeichne die Projektion auf die i-te Komponente. Wie man leicht nachrechnet, ist durch
	\begin{eqnarray*}
		\Phi: (G_1\times G_2)\wr H 	& \rightarrow & (G_1\wr H)\times (G_2\wr H) \\
				(\phi, h)			& \mapsto 	  & ( (p_1\circ\phi, h), (p_2\circ\phi, h))
	 \end{eqnarray*}
	ein injektiver Gruppenhomomorphismus gegeben. 
\end{beweis}

\begin{lemma}\label{lm:Kranzprodukt-Limes}
	Sei $\{G_i\}_{i\in I}$ ein gerichtetes System von Gruppen, dass durch Gruppenhomomorphismen $j_i:G_i \rightarrow G_{i+1}$ gerichtet ist. Sei $H$ eine endliche Gruppe der Ordnung $n$. F"ur das gerichtete System von Gruppen $\{G_i\wr H\}_{i\in I}$, gerichtet durch die Gruppenhomomorphismen
	\[
		J_i:G_i\wr H \rightarrow G_{i+1}\wr H, \quad (g,h) \mapsto (j_i\circ g,h),
	\]
	gilt dann
	\[ 
	\colim_{i\in I} (G_i \wr H) = (\colim_{i\in I} G_i) \wr H.
	\]
\end{lemma}
\begin{beweis}
	Wir fassen im folgenden Beweis $G^H$ immer als $G^n$ auf. 
	Mit den Homomorphismen $\psi_i: G_i \rightarrow \colim_{i\in I} G_i$ erhalten wir auch Homomorphismen
	\begin{eqnarray*}
				\Psi_i:\quad\qquad G_i \wr H & \rightarrow & (\colim_{i\in I}G_i)\wr H\\
							(g_1,\ldots,g_n,h) &\mapsto & (\psi_i(g_1),\ldots,\psi_i(g_n), h).
	\end{eqnarray*}
	Wir zeigen, dass $(\colim_{i\in I} G_i)\wr H$ zusammen mit den $\Psi_i$ die universelle Eigenschaft von $\colim_{i\in I}(G_i \wr H)$ erf"ullt.
	
	Sei $L$ eine Gruppe und $\phi_i:G_i\wr H\rightarrow L$ Homomorphismen, so dass $\phi_i = \phi_{i+1}\circ J_i$ f"ur alle $i\in I$.	
	Zu zeigen ist, dass es genau einen Homomorphismus 
	\[
	\Phi:(\colim_{i\in I}G_i)\wr H\rightarrow L
	\]gibt, so dass
  $\Phi\circ \Psi_{i} = \phi_i$ f"ur alle $i\in I$ gilt, also das folgendes Diagramm kommutiert
	
	
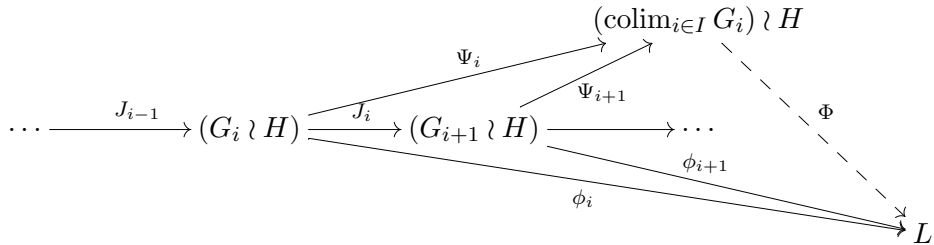
\begin{align*}
	\xymatrix@!C=60pt{
	   & & & 
		(\colim_{i\in I} G_i)\wr H
			\ar@{-->}[rdd]^{\Phi}
			&
		\\
			\ldots
			\ar@{->}[r]^{J_{i-1}}
		&
			(G_i\wr H)
			\ar[rrrd]_{\phi_i}
			\ar[rru]^{\Psi_i}
			\ar@{->}[r]^{J_i}
		&
			(G_{i+1}\wr H)
			\ar@{->}[r]
			\ar[rrd]^{\phi_{i+1}}
			\ar[ru]_{\Psi_{i+1}}
		& 
			\ldots
		&
		\\
		& & &  &
		 L	
	}
	\end{align*}	
	Sei $(g_1,\ldots,g_n,h)\in (\colim_{i\in I} G_i)\wr H$.
	F"ur jedes $g_s$ finden wir ein $g^{i_s}_s\in G_{i_s}$, $i_s\in I$ mit $\psi_{i_s}(g^{i_s}_s) = g_s$. Sei $N\coloneqq \max\{i_s\; |\; 1\leq s\leq n\}$ und $\widetilde{g}_s\coloneqq j_{N-1}\circ\dots\circ j_{i_s}(g_s)\in G_N$.
	Dann gilt $\Psi_N(\widetilde{g}_1,\ldots,\widetilde{g}_n,h) = (g_1,\ldots,g_n,h)$ und
	es bleibt keine andere Wahl als
	\[
	\Phi\left( (g_1,\ldots,g_n,h) \right) \coloneqq \phi_N\left((\widetilde{g}_1,\ldots,\widetilde{g}_n,h)\right)
	\]zu setzen.

	Da die $\phi_i$ Gruppenhomomorphismen sind, ist auch $\Phi$ Gruppenhomomorphismus, denn auch f"ur zwei Tupel $(a_1,\ldots,a_n,h_1),(b_1,\ldots,b_n,h_2)\in (\colim_{i\in I} G_i)\wr H$ finden wir ein $N\in I$ und $(\widetilde{a}_1,\ldots,\widetilde{a}_n,h_1), (\widetilde{b}_1,\ldots,\widetilde{b}_n,h_2)\in G_N\wr H$, so dass gilt
	\[
	\Phi( (a_1,\ldots,a_n,h_1) (b_1,\ldots,b_n,h_2) )=\phi_N( (\widetilde{a}_1,\ldots,\widetilde{a}_n,h_1)(\widetilde{b}_1,\ldots,\widetilde{b}_n,h_2)).
	\]
\end{beweis}
\rem{TODO }Das folgende Lemma wird sich im weiteren Verlauf der Arbeit noch als sehr n"utzlich erweisen. So ersetzt es zum Beispiel \cite[Lemma 5.4]{roushon} im Beweis von Satz~\ref{thm:FICwF-freiesProdukt} und es erm"oglicht au"serdem einen unkomplizierten Beweis von Satz~\ref{thm:FICwF-virtuell-normal-FICwF}, der uns mit Korollar \ref{kr:FICwF-Ueberlagerung} im Beweis des Hauptsatzes~\ref{thm:maintheorem} etwas Arbeit erspart.\rem{vereinfacht, indem es \cite[Lemma 6.6]{roushon} "uberfl"ussig macht.}
\begin{lemma}\label{lm:Kranzprodukt-Einbettung}
	Seien $G,H$ und $K$ Gruppen. $(G\wr H)\wr K$ ist eine Untergruppe von $G\wr(H\wr K)$.
\end{lemma}
\begin{beweis}
	Seien $(\phi, x),(\psi, y)\in (G\wr H)\wr K$, d.h. $x,y\in K$ und 
	\begin{eqnarray*}
	  \phi,\psi:K & \rightarrow & (G\wr H) \\
	  			k & \mapsto		& (f_k:H\rightarrow G, h_k).
	\end{eqnarray*}
	Mit $[\cdot]_i,~i=1,2$ wird im Folgenden die Projektion auf die $i$-te Komponente bezeichnet bzw. mit $\phi_i,\psi_i,~i=1,2$ die Verkn"upfung von $\phi$ und $\psi$ mit der Projektion auf die $i$-te Komponente. $\phi_1,\psi_1: K\rightarrow\left( H\rightarrow G\right)$ und $ \phi_2,\psi_2: K\rightarrow H$ sind Abbildung von Mengen.
	Wir definieren 
	\begin{eqnarray*}
		\Phi:(G\wr H)\wr K 	& \rightarrow	& G\wr(H\wr K) \\
		(\phi,x)		& \mapsto	& \left(\alpha_{\phi} ,\left( \phi_2, x\right)\right)
	\end{eqnarray*}
	mit
	\begin{eqnarray*}
	  	\alpha_{\phi}: \qquad\qquad H\wr K 
	  & \rightarrow & 
	  	G  \\
	  	(\widetilde{h}:K\rightarrow H, k) 
	  & \mapsto & 
	  	[\phi(k)]_1(\widetilde{h}(k))
	\end{eqnarray*}
	Zun"achst zeigen wir, dass $\Phi$ ein Gruppenhomomorphismus ist.
	\begin{eqnarray*}
		\Phi\left( \left( \phi,x\right)\left(\psi, y\right) \right) & = &\Phi\left( \phi\sd{x}{\psi},xy\right)  \\
		& = &\left( \alpha_{\phi\sd{x}{\psi}},\left( \left[\phi\sd{x}{\psi}\right]_2, xy \right) \right) \\
		& = &\left( \alpha_{\phi\sd{x}{\psi}},\left( \phi_2\sd{x}{\psi_2}, xy \right) \right) 
	\end{eqnarray*}
	\begin{eqnarray*}	
		\Phi\left( \phi,x \right) \Phi \left(\psi,y\right)
		& = & \left(\alpha_\phi, (\phi_2,x)\right)\left( \alpha_\psi, (\psi_2,y) \right) \\
		& = & \left(\alpha_\phi \sd{\left(\phi_2,x\right)}{\alpha_\psi}, (\phi_2,x)(\psi_2,y) \right) \\
		& = & \left(\alpha_\phi \sd{\left(\phi_2,x\right)}{\alpha_\psi}, (\phi_2\sd{x}{\psi_2},xy) \right)
		\end{eqnarray*}
	In der zweiten Komponente sieht man bereits, dass die beiden Ausdr"ucke "ubereinstimmen. F"ur die erste Komponente m"ussen wir noch etwas genauer hinsehen:

	Sei $(\widetilde{h}:K\rightarrow H,k)\in H\wr K$. Es gilt
	\begin{eqnarray*}
		\alpha_{\phi\sd{x}{\psi}}(\widetilde{h},k)	
		& = & [\underbrace{\phi(k)}_{\in G\wr H}\underbrace{\psi\left(x\invers k\right)}_{\in G\wr H}]_1 \left( \widetilde{h}(k)\right) \\
		& = & \left( [\phi(k)]_1[\sd{\phi_2(k)}{\psi(x^{-1} k)}]_1\right) \left(\widetilde{h}(k)\right)  \\ 
		& = & \left( \phi_1(k)\left( \widetilde{h}(k)\right) \right) \left( \psi_1(x\invers k) \left( \phi_2\invers(k)\widetilde{h}(k)\right) \right)
	\end{eqnarray*}
und
	\begin{eqnarray*}
		(\alpha_\phi \sd{\left(\phi_2,x\right)}{\alpha_\psi})(\widetilde{h},k) & = & \alpha_\phi(\widetilde{h},k)\s{\left(\phi_2,x\right)}{\!\alpha_\psi(\widetilde{h},k)} \\
		& = & \alpha_\phi(\widetilde{h},k)\alpha_\psi( (\phi_2,x)\invers (\widetilde{h},k)) \\ 
		& = & \alpha_\phi(\widetilde{h},k)\alpha_\psi( (\sd{x\invers}{\phi_2\invers},x\invers)(\widetilde{h},k)) \\ 
		& = & \alpha_\phi(\widetilde{h},k)\alpha_\psi( (\sd{x\invers}{\phi_2^{-1}}\sd{x\invers}{\widetilde{h}},x\invers k)) \\ 
		& = & \left( \phi_1(k)(\widetilde{h}(k))\right) \left(\psi_1(x\invers k)( (\sd{x\invers}{\phi_2}\sd{x\invers}{\widetilde{h}})(x\invers k)) \right) \\
		& = & \left( \phi_1(k)(\widetilde{h}(k))\right) \left( \psi_1(x\invers k) (\phi_2\invers(k) \widetilde{h}(k)) \right).
	\end{eqnarray*}
	
	Die Injektivit"at ist relativ offensichtlich, man muss sich nur einmal durch die Definition k"ampfen:
	
	Sei $(\phi,x)\in \ker(\Phi)$. $\phi_2$ und $x$ tauchen in $\Phi\left( (\phi,x) \right)=(\alpha_{\phi_1},(\phi_2,x))$ in der zweiten Komponente unver"andert auf, m"ussen also trivial sein. Angenommen $\phi_1$ ist nicht trivial. Dann gibt es ein $k\in K$ und ein $h\in H$ mit $\phi_1(k)(h) \neq 1$. 
	Sei $\widetilde{h}:K\rightarrow H$ eine Mengenabbildung mit $\widetilde{h}(k) = h$. Dann gilt $\alpha_{\phi_1}(\widetilde{h},k)=\phi_1(k)(\widetilde{h}(k))=\phi_1(k)(h)\neq 1$,  was im Widerspruch zu $(\phi,x)\in\ker(\Phi)$ steht. 
  \end{beweis}

\chapter{Isomorphismusvermutungen}
Die urspr"ungliche Formulierung der Farrel-Jones-Isomorphismusvermutung findet sich bei Farrell und Jones \cite[1.6]{farrell-jones}. Wir verwenden eine abgewandelte Definition mit Hilfe von "aquivarianten Homologietheorien, die beispielsweise bei L"uck-Reich \cite{lueck-reich} ausf"uhrlich beschrieben ist. \rem{TODO das hier auch noch? (vgl. Kapitel 1, \cite{bartels-lueck1}) TODO}
Wir beschr"anken uns dabei auf diskrete Gruppen. Ringe sind immer assoziativ mit Eins, aber nicht notwendigerweise kommutativ. 
Sei also im Folgenden $G$ eine diskrete Gruppe und $R$ ein assoziativer Ring mit Eins.

Wir ben"otigen zwei Konstruktionen:
\begin{enumerate}
	\item Den klassifizierenden Raum $\E$ zu einer Gruppe $G$ und einer Familie $\F$ von Untergruppen von $G$ und
	\item eine "aquivariante Homologietheorie $\Homo^?_*$.
\end{enumerate}

F"ur eine Gruppe $G$ ist eine \emph{Familie $\F$ von Untergruppen} eine Menge von Untergruppen von $G$, die abgeschlossen ist unter Konjugation und Untergruppenbildung. 
Ist $\phi:K\rightarrow G$ ein Gruppenhomomorphismus und $\F$ eine Familie von Untergruppen von $G$, dann definiere die Familie $\phi^*\F$ von Untergruppen von $K$ durch
\[
\phi^*\F \coloneqq \{H\subseteq K\; |\; \phi(H)\in \F\}.
\]
\begin{definition}
Sei $G$ eine Gruppe und $\F$ eine Familie von Untergruppen von $G$. Dann gibt es eindeutig bis auf $G$-Homotopie"aquivalenz einen $G$-CW-Komplex $\E$ mit der Eigenschaft, dass
\begin{enumerate}
  \item f"ur jedes $x\in \E$ die Standgruppe von $x$ in $\F$ ist und
  \item f"ur jedes $H\in \F$ die Fixpunktmenge $(E_{\F}G)^H$ kontraktibel ist.
\end{enumerate}
$\E$ hei"st \emph{Modell f"ur den klassifizierenden Raum der Familie $\F$}.
\end{definition}
Alternativ kann $\E$ auch mit Hilfe der universellen Eigenschaft charakterisiert werden, dass es f"ur jeden $G$-CW-Komplex $X$, dessen Standgruppen alle zu $\F$ geh"oren, bis auf $G$-Homotopie genau eine $G$-Abbildung von $X$ nach $\E$ existiert.
F"ur ausf"uhrlichere Details zu klassifizierenden R"aumen siehe \cite{lueck-cs}.

\pagebreak
Aufbauend auf der Definition einer $G$-Homologietheorie (siehe \cite[2.1.4]{lueck-reich})
definieren wir eine "aquivariante Homologietheorie.
\begin{definition}
	Eine \emph{"aquivariante Homologietheorie} $\Homo^?_*$ mit Werten in $R$-Moduln ordnet jeder Gruppe $G$ eine $G$-Homologietheorie $\Homo^G_*$ mit Werten in $R$-Moduln zu mit folgender Induktionsstruktur:
	Zu jedem Gruppenhomomorphismus $\alpha:H\rightarrow G$ gibt es f"ur $n\in\Z$ und alle $H$-CW-Paare $(X,A)$ einen nat"urlichen Homomorphismus
	\[
		\ind_\alpha: \Homo^H_n(X,A)\rightarrow \Homo^G_n(\ind_\alpha(X,A))
	\]
	mit $\ind_\alpha(X,A)\coloneqq G\times_\alpha(X,A)$ 
	und zwar derart, dass folgende Axiome erf"ullt sind:
	\begin{enumerate}
		\item
			Vertr"aglichkeit mit dem Randhomomorphismus

			$\partial^G_n\circ\ind_\alpha = \ind_\alpha\circ\;\partial^H_n$
		\item
			Funktorialit"at in $\alpha$

			Sei $\beta: G\rightarrow K$ ein Gruppenhomomorphismus, so dass $\ker(\beta\circ\alpha)$ frei auf $X$ operiert und
			\[
					f_1:\ind_\beta\ind_\alpha X\rightarrow\ind_{\alpha\circ\beta} X,\quad (k,g,x)\mapsto (k\beta(g),x)\]
			der nat"urliche $K$-Hom"oomorphismus. Dann gilt
	\[
	\ind_{\beta\circ\alpha} = \Homo^K_n(f_1)\circ\ind_\beta\circ\ind_\alpha:\Homo^H_n(X,A)\rightarrow \Homo^K_n(\ind_{\beta\circ\alpha}(X,A)).
  \]
	\item
		Vertr"aglichkeit mit Konjugation

		Sei $c(g):G\rightarrow G$ die Konjugation die $g'$ auf $gg'g\invers$ schickt.
		F"ur $n\in\Z, g\in G$ und ein $G$-CW-Paar $(X,A)$ stimmen die Homomorphismen
		\[
			\ind_{c(g)}:\Homo^G_n(X,A)\rightarrow\Homo^G_n(\ind_{c(g)}(X,A))
			\]
			und  $\Homo^G_n(f_2)$ "uberein, wobei $f_2:(X,A)\rightarrow\ind_{c(g)}(X,A)$ der durch\\$x\mapsto (1,g\invers x)$ definierte $G$-Hom"oomorphismus ist.
		\item
			Bijektivit"at

			Falls $\ker\alpha$ frei auf $X$ operiert, so ist 
			\[
				\ind_\alpha:\Homo^H_n(X,A)\rightarrow\Homo^G_n(\ind_\alpha(X,A))
			\]
				f"ur alle $n\in\Z$ bijektiv.
	\end{enumerate}
\end{definition}
\pagebreak

\begin{definition}[Isomorphismusvermutung]
  Sei 
	\begin{enumerate}
		\item $G$ eine diskrete Gruppe,
		\item $\F$ eine Familie von Untergruppen von $G$ und
		\item $\Homo^?_*$ eine "aquivariante Homologietheorie mit Werten in $R$-Moduln f"ur einen assoziativen Ring $R$.
	\end{enumerate}
	Das Tripel $(G,\F,\Homo^?_*)$ erf"ullt die Isomorphismusvermutung (IC), wenn die Projektion $\E\rightarrow\pt$ einen Isomorphismus
			\[
			\Homo^G_n(\E)\rightarrow\Homo^G_n(pt)
			\]
			f"ur alle $n\in\Z$ induziert.
\end{definition}
\begin{definition}[gefaserte Isomorphismusvermutung]
		$(G,\F,\Homo^?_*)$ erf"ullt die gefaserte Isomorphismusvermutung (FIC), wenn f"ur jeden Gruppenhomomorphismus $\phi:K\rightarrow G$ das Tripel $(K, \phi^*\F,\Homo^?_*)$ die Isomorphismusvermutung erf"ullt.
\end{definition}

		Aus der gefaserten Isomorphismusvermutung folgt schon die Isomorphismusvermutung, indem man als Gruppenhomomorphismus die Identit"at w"ahlt.
		
Im Folgenden sei $\F$ nicht mehr eine Familie speziell f"ur eine Gruppe gew"ahlt, sondern eine Klasse von Gruppen, die abgeschlossen unter Isomorphismen und Unter\-grup\-pen- und Quotientenbildung ist. F"ur eine Gruppe $G$ bezeichnen wir dann mit $\F_G$ die Familie, die alle Untergruppen von $G$ enth"alt, die zu $\F$ geh"oren. $\F_G$ ist dann wieder wie im urspr"unglichen Sinn eine Familie von Untergruppen von $G$.

Wir werden haupts"achlich mit einer noch st"arkeren Variante der gefaserten Isomorphismusvermutung arbeiten:	

\begin{definition}[endlich erweiterbare, gefaserte Isomorphismusvermutung]
	Ein Tripel $(G,\F_G,\Homo^?_*)$ erf"ullt $FICwF$, wenn $FIC$ f"ur jede endliche Gruppe $H$ f"ur das Tripel $(G\wr H,\F_{G\wr H},\Homo^?_*)$ erf"ullt ist.
\end{definition}

Indem man $H$ trivial w"ahlt, sieht man sofort, dass aus $FICwF$ schon $FIC$ folgt.
$FICwF$ hat die n"utzliche Eigenschaft, dass sie sich, wie wir mit Satz~\ref{thm:FICwF-virtuell-normal-FICwF} sehen werden, auf endliche Erweiterung vererbt.

F"ur beliebige Tripel sind diese Isomorphismusvermutungen falsch. Das Tripel muss also richtig gew"ahlt werden. 
Ein besonderes Augenmerk werden wir auf die Farrel-Jones-Isomorphismus\-ver\-mutung f"ur $L$-Theorie und algebraische $K$-Theorie richten. Daf"ur fixieren wir f"ur $\F$ die Klasse aller virtuell zyklischen Gruppen $\VC$ und als "aquivariante Homologietheorien w"ahlen wir $\Homo^?_n(-;\KT_R)$ und $\Homo^?_n(-;\LT_R)$ wie beispielsweise von L"uck und Reich in \cite[Proposition 6.7]{lueck-reich} definiert.\rem{ mit der Eigenschaft, dass $H^G_n(\pt;\KT_R)\cong \KT_n(RG)$ bzw. $H^G_n(\pt;\LT_R)\cong\LT_n(RG)$. }\rem{TODO}

		\begin{definition}[$L$-$FIC$ und $K$-$FIC$]
			Eine Gruppe $G$ erf"ullt die \emph{Farrell-Jones-Isomorphismusvermutung f"ur $L$-Theorie bzw. algebraische $K$-Theorie}, wenn
		 die Isomorphismusvermutung von $(G, \VC_G, \Homo^?_*(-,\LT_R))$ bzw. $(G, \VC_G, \Homo^?_*(-,\KT_R))$ er\-f"ullt wird, also die Projektion $E_{\VC_G}G\rightarrow\pt$ f"ur alle $n\in\Z$ Isomorphismen
			\begin{eqnarray*}
  				\Homo^{G}_n(E_{\VC_G}G;\KT_R)& \xrightarrow{\cong} &\Homo^{G}_n(\mbox{pt};\KT_R) \\
					\Homo^{G}_n(E_{\VC_G}G;\LT_R)& \xrightarrow{\cong} &\Homo^{G}_n(\mbox{pt};\LT_R)
			\end{eqnarray*}
			induziert.

			$G$ erf"ullt die (endlich erweiterbare) gefaserte Farrell-Jones-Isomorphismus\-ver\-mu\-tung f"ur $L$-Theorie bzw. algebraische $K$-Theorie, wenn $(G, \VC_G, \Homo^?_*(-,\LT_R))$ bzw.\\ $(G, \VC_G, \Homo^?_*(-,\KT_R))$ die (endlich erweiterbare) gefaserte Isomorphismusvermutung erf"ullen.
			Wir schreiben daf"ur kurz $G$ erf"ullt $L$-$FIC(wF)$ bzw. $K$-$FIC(wF)$.
	\end{definition}
\chapter{FIC-Werkzeugkasten}
\section{Schraubenzieher}
Wir wenden uns zun"achst nicht allzu tiefliegenden Resultaten zu, die uns aber nichtsdestotrotz
im weiteren Verlauf dieser Arbeit noch wichtige Dienste leisten werden. Die meisten der folgenden Resultate finden sich "ubersichtlich zusammengefasst in \cite[Abschnitt 1]{lueckmain}. Wir erweitern einige lediglich von $FIC$ auf $FICwF$.

Es sei bei den folgenden Aussagen immer eine "aquivariante Homologietheorie fest gew"ahlt und teilweise auch eine Klasse $\F$. Wir sagen ein Paar $(G,\F_G)$ erf"ullt $FIC$ ($FICwF$), wenn $(G,\F_G)$ bez"uglich der gew"ahlten Homologietheorie die (endlich erweiterbare) gefaserte Isomorphismusvermutung erf"ullt bzw. wenn zus"atzlich eine Klasse $\F$ fest gew"ahlt auch nur kurz "`$G$ erf"ullt $FIC(wF)$"'.

Schon in die Definition der gefaserten Isomorphismusvermutung eingebaut, ist 
\begin{lemma}\label{lm:FIC-zurueckziehen}
	Sei $\phi: K\rightarrow G$ ein Gruppenhomomorphismus. Wenn $FIC$ f"ur $(G,\F_G)$ erf"ullt ist, dann ist $FIC$ auch f"ur $(K,\phi^*\F_G)$ erf"ullt.
\end{lemma}
\begin{beweis}
	Sei $L$ eine Gruppe und $\psi:L\rightarrow K$ ein Gruppenhomomorphismus und f"ur $(G,\F_G)$ gelte $FIC$. Nach Definition von $FIC$ ist f"ur $(K, (\phi\circ\psi)^*\F_G)$ die Isomorphismusvermutung erf"ullt. Es gilt aber $(K,(\phi\circ\psi)^*\F_G)=(K,\psi^*(\phi^*\F_G))$, womit $FIC$ f"ur $(K,\phi^*\F_G)$ erf"ullt ist.
\end{beweis}

\begin{lemma}\label{lm:FICwF-subgroup}
	Ist $FICwF$ f"ur ein Paar $(G,\F_G)$ erf"ullt, dann wird $FICwF$ auch von $(K,\F_K)$ erf"ullt f"ur jede Untergruppe $K$ von $G$.
\end{lemma}
\begin{beweis}
	F"ur $FIC$ erhalten wir die Aussage direkt aus Lemma~\ref{lm:FIC-zurueckziehen}, denn mit der Inklusion $i:K\rightarrow G$ gilt $FIC$ f"ur 
	\[(K,i^*\F_G)=(K,\F_K).\]

  Sei nun $G$ eine Gruppe, die $FICwF$ erf"ullt, $K$ eine Untergruppe von $G$ und $H$ eine beliebige endliche Gruppe. Aufgrund von Lemma~\ref{lm:Kranzprodukt-Untergruppen} ist $K\wr H$ eine Untergruppe von $G\wr H$ und da $FIC$ von $G\wr H$ erf"ullt wird, ist $FIC$ auch f"ur $K\wr H$ wahr. Es gilt also $FICwF$ f"ur $K$.
\end{beweis}

Betrachtet man $FICwF$ bez"uglich Familien virtuell zyklischer Untergruppen oder anderer Familien, die endliche Gruppen einschlie"sen, dann gilt $FICwF$ schon direkt f"ur endliche Gruppen.

\begin{lemma}\label{lm:FICwF-endliche-Gruppen}
	Sei $G$ eine endliche Gruppe und $\F$ eine Klasse von Gruppen, die alle endlichen Gruppen enth"alt. Dann wird $FICwF$ von $(G,\F_G)$ erf"ullt.
\end{lemma}
\begin{beweis}
	Sei $\phi:L\rightarrow G$ ein Gruppenhomomorphismus.
	Da $G$ endlich ist, gilt $G\in\F_G$ und damit auch $L\in\phi^*\F_G$.
	Somit ist schon ein Punkt ein Modell f"ur den klassifizierenden Raum von $E_{\phi^*\F_G}L$ und insbesondere
		\[
		\Homo^G_n(E_{\phi^*\F_G}G)\rightarrow\Homo^G_n(pt)
	\]
ein Isomorphismus.
	F"ur eine beliebige endliche Gruppe $H$ ist auch $G\wr H$ endlich und mit dem gleichen Argument gilt $FIC$ f"ur $G\wr H$, also $FICwF$ f"ur $G$.	
\end{beweis}

Wir brauchen die Vertr"aglichkeit von $FICwF$ mit Kolimiten bez"uglich gerichteter Systeme von Gruppen. Die erhalten wir, wenn wir eine bestimmte Voraussetzung an die "aquivariante Homologietheorie stellen:

\begin{definition}[streng stetige "aquivariante Homologietheorie]
  Ein "aquivariante Homologietheorie $\Homo^{?}_*$ hei"st \emph{streng stetig}, wenn f"ur jede Gruppe $G$ und ein gerichtetes System von Gruppen $\{G_i \; | \; i\in I\}$ mit $G=\colim_{i\in I} G_i$ f"ur jedes $n\in\Z$ die Abbildung
  \[
  \colim_{i\in I} j_{i}:\colim_{i\in I}\Homo^{G_i}_n(\pt)\rightarrow\Homo^G_n(\pt)
  \]
  ein Isomorphismus ist, wobei $j_i:\Homo^{G_i}_n(\pt)\rightarrow \Homo^G_n(\pt)$ die Komposition des Induktionshomomorphismus $\Homo^{G_i}_n(\pt)\xrightarrow{\cong}\Homo^G_n(G/G_i)$ mit der Abbildung ist, die von der Projektion $G/G_i \rightarrow \pt$ induziert wird.
\end{definition}

\begin{satz}\label{thm:FICwF-limes}
	Sei $\Homo^?_*$ eine streng stetige "aquivariante Homologietheorie. Sei $\{G_i\; |\; i\in I\}$ ein gerichtetes System von Gruppen, das durch Homomorphismen $j_i:G_i\rightarrow G_{i+1}$ gerichtet ist. 
	Falls $FICwF$ f"ur $(G_i,\F_{G_i},\Homo^?_*)$ f"ur alle $i\in I$ erf"ullt ist, dann ist $FICwF$ auch erf"ullt f"ur $(\colim_{i\in I}G_i,\F_{\colim_{i\in I}G_i},\Homo^?_*)$.
\end{satz}
\begin{beweis}
	Sei $H$ eine beliebige endliche Gruppe. Nach der Definition von $FICwF$ ist $FIC$ wahr f"ur $(G_i\wr H,\F_{G_i\wr H},\Homo^?_*)$ f"ur alle $i\in I$. 
	Mit den Homomorphismen 
	\[J_i:G_i\wr H\rightarrow G_{i+1}\wr H, \quad (f,h)\mapsto (j_i\circ f, h)\] bilden die $\{G_i\wr H\; | \; i\in I\}$ ein gerichtetes System von Gruppen. 

	Sei $G\coloneqq\colim_{i\in I}(G_i\wr H)$.
	Da wir $\Homo^?_*$ als streng stetig vorausgesetzt haben, gilt $FIC$ f"ur $(G,\F_G,\Homo^?_*)$.
Ein Beweis dazu findet sich bei Bartels, Echterhoff und L"uck \cite[Proposition 4.6]{echterhoff}.
	Es gilt mit Lemma~\ref{lm:Kranzprodukt-Limes} 
	\[
	G=\colim_{i\in I} (G_i\wr H)\cong (\colim_{i\in I} G_i)\wr H.\]
	Also ist $FIC$ wahr f"ur $(\colim_{i\in I} G_i)\wr H$ und damit $FICwF$ f"ur $\colim_{i\in I}G_i$.
\end{beweis}
\rem{\begin{beweis}
	Sei $H$ eine beliebige endliche Gruppe. Nach der Definition von $FICwF$ ist $FIC$ wahr f"ur $G_i\wr H$ f"ur alle $i\in I$. 
	Mit den Homomorphismen \[J_i:G_i\wr H\rightarrow G_{i+1}\wr H, \quad (f,h)\mapsto (j_i\circ f, h)\] bilden die $\{G_i\wr H\; | \; i\in I\}$ ein gerichtetes System von Gruppen. F"ur streng stetige "aquivariante Homologietheorien wissen wir mit \cite[Proposition 4.6]{echterhoff}, dass dann $FIC$ f"ur $\colim_{i\in I} (G_i\wr H)$ wahr ist. Es gilt mit Lemma~\ref{lm:Kranzprodukt-Limes} $\colim_{i\in I} (G_i\wr H)\cong (\colim_{i\in I} G_i)\wr H$, also ist $FICwF$ wahr f"ur $G$.
\end{beweis}}

F"ur das folgende Lemma nutzen wir das Transitivit"atsprinzip, das beispielsweise in \cite[Theorem 1.5]{bartels-lueck2} \rem{FTODO} bewiesen wird.
 \begin{satz}[Transitivit"atsprinzip]\label{thm:FIC-transitiv}
	Seien $\F\subseteq\G$ zwei Familien von Untergruppen von $G$ und f"ur jede Gruppe $H\in\G$ sei $FIC$ erf"ullt f"ur $(H,\F\cap H)$.  Dann ist $FIC$ f"ur $(G,\G)$ genau dann wahr, wenn $FIC$ f"ur $(G,\F)$ wahr ist.
\end{satz}

\begin{lemma}\label{lm:2.2}
	Sei $p:K \rightarrow G$ ein Gruppenhomomorphismus. $FIC$ sei erf"ullt f"ur $(G,\F_G)$ und f"ur $(p^{-1}(C),\F_{p\invers(C)})$ f"ur alle $C\in\F_G$. Dann wird $FIC$ auch von $(K,\F_K)$ erf"ullt.
\end{lemma}
\begin{beweis}
	Sei $q:L\rightarrow K$ ein Gruppenhomomorphismus. Wir haben zu zeigen, dass $IC$ von $(L,q^*\F_K)$ erf"ullt wird.
	Es gilt $\F_K \subseteq p^*\F_G$, denn da $\F$ abgeschlossen unter Quotientenbildung ist, gilt $p(A)\cong A/\ker(p|_A)\in\F$  f"ur $A\in\F_K$ ist, also $p(A)\in\F_G$. Damit ist auch $q^*\F_K \subseteq q^*p^*\F_G$.
	
	Wir wollen nun f"ur die Familien $q^*\F_K$ und $q^*p^*\F_G$ das Transitivit"atsprinzip anwenden.
	Sei $H\in q^*p^*\F_G$. Dann gibt es ein $V\in\F_G$ mit $q(H)\subseteq p\invers(V)$.
	F"ur das Transitivit"atsprinzip brauchen wir $FIC$ f"ur $(H,q^*\F_K\cap H)$.
	Es gilt aber 
	\[q^*\F_K\cap H = (q|_H)^*\F_K = (q|_H)^*\F_{p\invers(V)}.\]
	F"ur $(p\invers(V),\F_{p\invers(V)})$ ist $FIC$ nach Voraussetzung erf"ullt.
	Wir wenden Lemma~\ref{lm:FIC-zurueckziehen} auf $q|_H$ an und erhalten $FIC$ f"ur $(H,(q|_H)^*\F_{p\invers(V)})$.
\end{beweis}

Wir werden sp"ater mit Lemma~\ref{lm:FICwF-zurueckziehen} noch zeigen, dass sich auch dieses Resultat auf $FICwF$ erweitern l"asst, wenn wir einige zus"atzliche Voraussetzungen annehmen.

Der n"achste Satz rechtfertigt schlie"slich, warum wir uns "uberhaupt mit Kranzprodukten und $FICwF$ statt $FIC$ besch"aftigen.

\begin{satz}\label{thm:FICwF-virtuell-normal-FICwF}
	Sei $G$ eine Gruppe und $N$ ein Normalteiler von $G$ mit endlichem Index. Ist $FICwF$ f"ur $N$ erf"ullt, dann ist $FICwF$ auch schon f"ur $G$ wahr.
\end{satz}
\begin{beweis}	
	Wir haben folgende exakte Sequenz von Gruppen \exseq{N}{G}{G/N}{.}
	Nach Satz~\ref{thm:Kranzprodukt-Erweiterungen} ist $G$ Untergruppe von $N\wr (G/N)$.
	F"ur eine beliebige endliche Gruppe $H$ gilt also mit Lemma~\ref{lm:Kranzprodukt-Untergruppen} und Lemma~\ref{lm:Kranzprodukt-Einbettung}
	\begin{eqnarray*}
		G\wr H  \stackrel{}{\subset}  (N\wr (G/N))\wr H 
		\stackrel{}{\subset}  N\wr ( (G/N)\wr H) 
	\end{eqnarray*}
  $(G/N)\wr H$ ist eine endliche Gruppe, da sowohl $G/N$ als auch $H$ endlich sind.
  Nach Definition von $FICwF$ ist $FIC$ f"ur $N\wr ( (G/N)\wr H)$ wahr, 
  was sich mit Lemma~\ref{lm:FICwF-subgroup} auf die Untergruppe $G\wr H$ "ubertr"agt. 
  Damit wird $FICwF$ von $G$ erf"ullt.
\end{beweis}

Die Voraussetzung, dass die Untergruppe normal sein muss, kann fallengelassen werden, wenn wir uns auf endlich erzeugte Gruppen einschr"anken, wie wir im Folgenden sehen. 

\begin{korollar}\label{kr:FICwF-endl-erz-virtuell-FICwF}
	Sei $G$ eine endlich erzeugte Gruppe.	
	Erf"ullt $G$ virtuell $FICwF$, d.h. es gibt eine Untergruppe mit endlichem Index in $G$, die $FICwF$ erf"ullt, dann gilt $FICwF$ auch f"ur $G$.
\end{korollar}
\begin{beweis}
  Sei $H$ eine Untergruppe von $G$ mit endlichem Index $n$, die $FICwF$ erf"ullt.
  Definiere $N\coloneqq \bigcap_{g\in G}gHg\invers$. 
  $N$ ist Normalteiler von $G$ und eine Untergruppe von $H$.
	Da $G$ endlich erzeugt ist, gibt es zu einem festen endlichen Index $n$ nur endliche viele Untergruppen (siehe \cite{baumslag}, Kapitel III, Theorem 3), das hei"st wir finden eine endliche Teilmenge $M\subset G$ mit $N\cong \bigcap_{g\in M}gHg\invers$.
	Als endlicher Schnitt von Gruppen mit endlichem Index hat $N$ auch endlichen Index in $G$.
	
	$N$ erf"ullt als Untergruppe von $H$ mit Lemma~\ref{lm:FICwF-subgroup} auch $FICwF$ und mit Satz~\ref{thm:FICwF-virtuell-normal-FICwF} gilt $FICwF$ dann auch f"ur $G$.
 \end{beweis}

	Zum Abschluss dieses Abschnittes liefert uns ein zweites Korollar eine erste Aussage "uber Fundamentalgruppen von Mannigfaltigkeiten, die wir sp"ater verwenden werden.
\begin{korollar}\label{kr:FICwF-Ueberlagerung}
  Sei $M$ eine kompakte 3"=Mannigfaltigkeit mit endlicher "Uberlagerung $N$, deren Fundamentalgruppe $FICwF$ erf"ullt. Dann gilt $FICwF$ auch schon f"ur die Fundamentalgruppe von $M$.
\end{korollar}
\begin{beweis}
  	Als endliche "Uberlagerung einer kompakten Mannigfaltigkeit ist $N$ auch wieder kompakt. Damit ist $\pi_1(N)$ endlich erzeugt. 
 	Die lange exakte Homotopiesequenz liefert uns $\pi_1(N)$ als Untergruppe von $\pi_1(M)$ mit endlichem Index. Mit dem vorherigen Korollar wird damit $FICwF$ von $\pi_1(M)$ erf"ullt. 
\end{beweis}
\vfill
\pagebreak

\section{Bohrmaschinen}
Wir kommen nun von den leichten Schraubenziehern zu schwererem Ger"at. Als Ziel hatten wir uns aufgestellt:
\begin{zielsatz}
	Eine endlich erweiterbare, gefaserte Isomorphismusvermutung, die die Anforderungen W1 bis W5 erf"ullt, gilt f"ur Fundamentalgruppen von dreidimensionalen Mannigfaltigkeiten.
\end{zielsatz}
Wir werden nun erl"autern, was es mit den Werkzeugen W1 bis W5, die wir f"ur unseren Beweis des Zielsatzes ben"otigen, auf sich hat, und im Anschluss mit deren Hilfe noch einige weitere wichtige Werkzeuge beweisen.
\begin{enumerate}
	\item[W1:] $FICwF$ gilt f"ur alle Gruppen der Gestalt $\Z^2\sdp\Z$.
	\item[W2:] $FICwF$ gilt f"ur Fundamentalgruppen von nicht-positiv ge\-kr"umm\-ten, geschlossenen Mannigfaltigkeiten. 
	\item[W3:] $FIC$ gilt f"ur virtuell abelsche Gruppen.
	\item[W4:] Sei $\{G_i\; | \; i\in I\}$ ein gerichtetes System von Gruppen und f"ur jedes $G_i$ gelte $FIC$. Dann gilt $FIC$ auch f"ur $\colim_{i\in I}G_i$. 
	\item[W5:] Sei $p:L\rightarrow G$ ein Gruppenhomomorphismus und $FIC$ gelte f"ur $G$. Ist $FIC$ wahr f"ur $p\invers(V)$ f"ur jede virtuell zyklische Untergruppe $V$ von $G$, dann gilt $FIC$ auch f"ur $L$.
\end{enumerate}

F"ur eine beliebige Isomorphismusvermutung wissen wir nicht, ob tats"achlich all diese Anforderungen erf"ullt werden. Aber wenn wir uns bei der Wahl der Isomorphisvermutung etwas einschr"anken, bekommen wir einige dieser Werkzeuge zur Verf"ugung gestellt, was wir im Folgenden genauer ausf"uhren werden. 

W5 ist in den zentralen S"atzen~\ref{thm:FICwF-Seifert} und \ref{thm:maintheorem} ein entscheidendes Beweismittel. Betrachten wir
Isomorphismusvermutungen bez"uglich Familien virtuell zyklischer Untergruppen, ist W5 nichts anderes als die Aussage von Lemma~\ref{lm:2.2}. F"ur jede Unterfamilie, wie beispielsweise die Familie endlicher Untergruppen oder die Familie der trivialen Untergruppe, gilt ebenfalls W5.


Spezialisieren wir uns noch weiter auf die Farrell-Jones-Isomorphismusvermutung f"ur $L$-Theorie oder algebraische $K$-Theorie, bekommen wir mit Satz~\ref{thm:FICwF-limes} auch noch W4, da in beiden F"allen die "aquivarianten Homologietheorien streng stetig sind. Ein Beweis daf"ur findet sich von Bartels, L"uck und Echterhoff in \cite[Lemma 5.2]{echterhoff}, die dort auch noch einige andere streng stetige "aquivariante Homologietheorien auflisten.

W2 erhalten wir, wenn wir wissen, dass $FIC$ f"ur Gruppen gilt, die proper, kokompakt und isometrisch auf einem $\CAT(0)$-Raum operieren. Wir werden das in Satz~\ref{thm:FICwF-CAT0} beweisen. Bei Bartels und L"uck schlummert ein unver"offentlichter Beweis\footnote{inzwischen zu finden in: A. Bartels, W. L"uck, The Borel Conjecture for hyperbolic and $\CAT(0)$-groups, arXiv:math.GT/0901.0442, 2009.}, dass $L$-$FIC$ f"ur $\CAT(0)$-Gruppen gilt und $K$-$FIC$ immerhin f"ur $n<1$ f"ur $\CAT(0)$-Gruppen wahr ist.

\rem{*W8 ist aber f"ur $K$-$FIC$ in \cite[Lemma 2.5]{lueckmain}}
W3 hat Quinn in \cite[Theorem 1.2.2]{quinn} schon f"ur $K$-$FIC$ gezeigt. F"ur die Farrell-Jones-Isomorphismusvermutung in algebraischer $K$-Theorie m"ussen wir also lediglich folgende zwei Werkzeuge als Voraussetzung fordern:

\begin{enumerate}
	\item[K-W1:] Alle Gruppen der Gestalt $\Z^2\sdp\Z$ erf"ullen $K$-$FICwF$.
	\item[K-W2:] F"ur $n\geq 1$ gilt $K$-$FIC$ f"ur Gruppen, die proper, kokompakt und isometrisch auf einem endlich dimensionalen $\CAT(0)$-Raum operieren.
\end{enumerate}

F"ur $L$-$FIC$ fehlt uns leider ein "ahnliches Resultate zu W3, aber daf"ur k"onnen wir davon ausgehen, dass $L$-$FIC$ f"ur $\CAT(0)$-Gruppen gilt. Es bleiben deshalb f"ur den $L$-Theorie-Fall folgende zwei Voraussetzungen "ubrig:

\begin{enumerate}
	\item[L-W1:] Alle Gruppen der Gestalt $\Z^2\sdp\Z$ erf"ullen $L$-$FICwF$.
	\item[L-W2:] $L$-$FIC$ gilt f"ur virtuell abelsche Gruppen.
\end{enumerate}

Gruppen der Gestalt $\Z^2\sdp\Z$ sind ein Spezialfall von virtuellen poly-$\Z$-Gruppen, das hei"st Gruppen, die mit endlichem Index eine Untergruppe enthalten, die eine Normalreihe besitzt mit Faktoren isomorph zu $\Z$. Zur Farrell-Jones-Isomorphis\-musvermutung und virtuell poly-$\Z$-Gruppen gibt es zwar von Farrell und Jones mit \cite[Proposition 2.2]{farrell-jones} eine Aussage, auf die auch Roushon in seiner Arbeit "uber 3"=Mannigfaltigkeiten verweist (\cite[Proposition 2.1]{roushon}), allerdings ist das im Beweis von Farrell und Jones verwendete Theorem 4.8 (\cite[S. 283]{farrell-jones}) nachwievor unbewiesen, so dass wir nicht darum herum kommen K-W1 und L-W1 als Voraussetzung zu fordern.

F"ur die (teilweise) Elimination von W2 als Voraussetzung f"ur $L$-$FIC$ und $K$-$FIC$ m"ussen wir noch folgenden Satz nachreichen:
\begin{satz}\label{thm:FICwF-CAT0}
	$FIC$ gelte f"ur Gruppen die proper, kokompakt und isometrisch auf einem $\CAT(0)$-Raum operieren. Dann gilt $FICwF$ f"ur Fundamentalgruppen von geschlossenen, nicht-positiv gekr"ummten Mannigfaltigkeiten.
\end{satz}
\begin{beweis}
	Sei $M$ eine geschlossene nicht-positiv gekr"ummte Mannigfaltigkeit, $\widetilde{M}$ ihre, ebenfalls nicht-positiv gekr"ummte, universelle "Uberlagerung, die insbesondere ein $\CAT(0)$-Raum ist. 
	Sei $G\coloneqq \pi_1(M)$ und $H$ eine endliche Gruppe. 
	Wir m"ussen $FIC$ f"ur $G\wr H$ zeigen und definieren daf"ur eine propere, kokompakte und isometrische Gruppenoperation von $G\wr H$ auf $\widetilde{M}^H$\!.
	$G$ operiert auf $\widetilde{M}$ bereits proper, kokompakt und isometrisch und zudem frei und eigentlich diskontinuierlich\rem{ TODO isometrisch? TODO.} und ebenso $G^H$ auf $\widetilde{M}^H$. Wir m"ussen uns also nur noch n"aher anschauen, was das $H$ aus $G\wr H$ bewirkt.
	Wir definieren die Operation von $G\wr H$ auf $\widetilde{M}^H$ durch 
	\begin{eqnarray*}
		\s{(f,0)}{\left(m_h\right)_{h\in H}}& = &\left(^{f(h)}m_{h}\right)_{h\in H} \\
		\s{(0,a)}{\left(m_h\right)_{h\in H}}& = &\left(m_{a\invers h}\right)_{h\in H}	
	\end{eqnarray*}	
	f"ur $(f,0),(0,a)\in G\wr H$ und  $(m_h)_{h\in H}\in \widetilde{M}^H$.
Es l"asst sich leicht nachrechnen, dass dies tats"achlich eine Gruppenoperation definiert.

$(0,a)$ permutiert lediglich die Eintr"age von $(m_h)_{h\in H}$.
	W"ahlen wir auf $\widetilde{M}^n$ die Produktmetrik $d^n$ (mit $d^n((m_h)_{h\in H}, (m'_h)_{h\in H})^2 = \sum_{h\in H}d^n(m_h,m'_h)^2$), so ist diese Permutation isometrisch. Au"serdem ist $\widetilde{M}^H/(G\wr H) \cong M^H/H$ als Projektion einer kompakten Mannigfaltigkeit kompakt. Bleibt zu zeigen, dass $G\wr H$ auch proper operiert.

	Sei $(m,n)\coloneqq ( (m_h)_{h\in H}, (n_h)_{h\in H} ) \in \widetilde{M}^H \times \widetilde{M}^H$.
	Da $G\wr H$ diskret ist, operiert $G\wr H$  nach \cite[Korollar 3.1.6]{proper-group-actions} proper, wenn es f"ur $(m,n)$ eine Umgebung $W$ gibt, so dass
	\[
	Z\coloneqq \left\{ (f,a)\in G\wr H \; \left| \; ^{(f,a)}m' = n' \mbox{ f"ur } (m',n')\in W\right. \right\}
	\]
	endlich ist.
\rem{
	Daf"ur ist zun"achst zu zeigen, dass die Standgruppen endlich sind. Sei $p:(G\wr H)\rightarrow H$ die Projektion und $H'\coloneqq p( (G\wr H)_x)$ f"ur das Bild der Standgruppe von $x\in\widetilde{M}^H$ unter $p$..
	Dann ist
	\gexseq{\ker(p')}{(G\wr H)_x}{H'}{}{p'}{}
	exakt, wobei $p'= p|_{(G\wr H)_x}$.
	Es gilt 
	\[\ker(p')=\left\{(f,a)\in \left(G\wr H\right)_x \;|\; a = 0\right\} = \left(G^H\right)_x.\]
	Da $G^H$ frei operiert, ist $\ker(p')$ trivial und damit $(G\wr H)_x \cong H' \subseteq H$ endlich.}

$G$ operiert eigentlich diskontinuierlich auf $\widetilde{M}$, das hei"st wir finden f"ur alle $m_h$ aus $m$ Umgebungen  $\widetilde{U}_h$, so dass $\widetilde{U}_h\cap g\widetilde{U}_h=\emptyset$ f"ur alle $g\in G$ und analog auch Umgebungen $\widetilde{V}_h$ f"ur alle $n_h$ aus $n$. 
	W"ahle $\frac{\varepsilon}{2}$-Umgebungen $U_h$ und $V_h$ f"ur alle $h\in H$ (also $\varepsilon>0$ so, dass $d(m_h,m'_h)< \frac{\varepsilon}{2} \Rightarrow m'_h\in U_h$ und $d(m_h,m'_h)<\varepsilon \Rightarrow 	m'_h\in\widetilde{U_h}$ und analog f"ur $V_h$). 
	
	Wir wollen zeigen, dass es f"ur zwei solche $\frac{\varepsilon}{2}$-Umgebungen $U_i$ und $V_j, i,j\in H$, h"ochstens ein $g\in G$ existiert mit $gU_i\cap V_j\neq\emptyset$.
	Wir nehmen also an zu $a_1,a_2\in U_i$ und $b_1,b_2\in V_j$ gebe es $p,q\in G, p\neq q$ mit $^pa_1=b_1$ und $^qa_2=b_2$.
	Es gilt 
	\[
	d(m_i,a_1)<\frac{\varepsilon}{2} \mbox{ und } d(m_i,a_2)<\frac{\varepsilon}{2}.
	\]
	Da $G$ isometrisch operiert, gilt dann auch 
	\[d(^pm_i,b_1) = d(^pm_i,\s{p}{a_1})<\frac{\varepsilon}{2}\mbox{ und } d(^qm_i,b_2) = d(^qm_i,\s{q}{a_2})<\frac{\varepsilon}{2}. 
	\]
	Also $d(^pm_i,n_j),d(^qm_i,n_j)<\varepsilon$ und damit $^pm_i,\s{q}{m_i}\in \widetilde{V}_j$, was zum Widerspruch $\widetilde{V}_j\cap\, {pq\invers}{\widetilde{V}_j} \neq\emptyset$ f"uhrt.
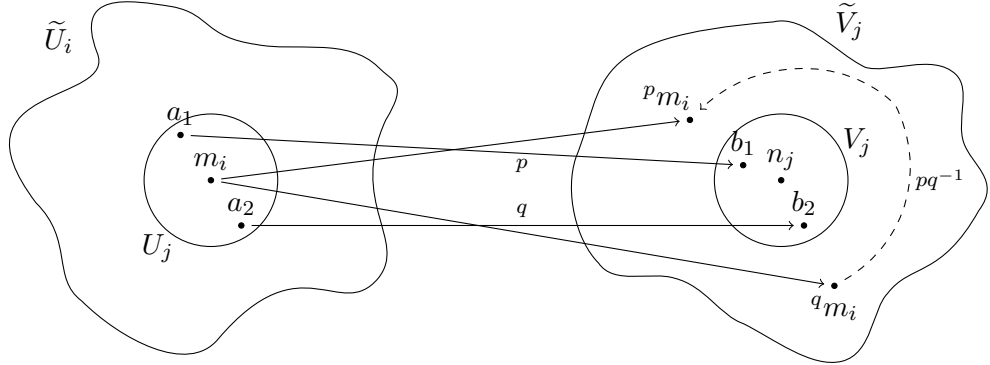
\begin{center}
\begin{tikzpicture}
	\node (links) at (-4.5,0) {};
	\node (rechts) at (3,0) {};
	\node (a1) at (-4.9,0.6) {};
	\node (a2) at (-4.1,-0.6) {};
	\node (b1) at (2.5,0.2) {};
	\node (b2) at (3.3,-0.6) {};
	\node (pm) at (1.8,0.8) {};
	\node (qm) at (3.7,-1.4) {};
	\node (p) at (-0.4,0.2) {${_p}$};
	\node (q) at (-0.4,-0.4) {${_q}$};
	\draw (-6.5,1.9) node {$\widetilde{U}_i$}
				(3.9,2.1) node {$\widetilde{V}_j$}
				(-5.2,-0.9) node {$U_j$}
				(4,0.5) node {$V_j$};
	\draw  (-6.2,-1.6) .. controls (-6.6,-1.3) and (-6.3,-0.8) ..	(-6.8,-0.3) 
								.. controls (-8,1.3) and (-5.7,1) .. (-6,1.7) 
								.. controls (-6.3,2.5) and (-5.6,2.4) .. (-5,2.3) 
								.. controls (-4.5,2.2) and (-4.6,1.9) .. (-3.9,1.8)
								.. controls (-3.4,1.7) and (-2.7,1.9) .. (-2.3,1.5) 
								.. controls (-2.1,1.3) and (-1.9,1.3) .. (-2.2,0.8)
								.. controls (-2.7,-0.3) and (-1.9,-0.8) .. (-2.3,-1.3)
								.. controls (-2.8,-1.9) and (-3.3,-0.8) .. (-4.4,-2.1)
								.. controls (-4.8,-2.6) and (-5.7,-2.1) .. (-6.2,-1.6)	-- cycle;
	\draw  (0.3,0.3) .. controls (0.7,1.3) and (0.4,1.9) ..	(2.9,2.1) 
								.. controls (3.3,2.2) and (3.7,1.7) .. (4.1,1.5) 
								.. controls (4.4,1.4) and (5,1.7) .. (5.3,1.4) 
								.. controls (5.5,1.2) and (5.3,0.8) .. (5.4,0.5)
								.. controls (5.5,0.2) and (5.9,0) .. (5.6,-0.4) 
								.. controls (5.1,-1.3) and (4.7,-1.1) .. (4.4,-1.5)
								.. controls (4.2,-1.7) and (4.2,-2.3) .. (3.8,-2.4)
								.. controls (3.4,-2.5) and (2.9,-2.1) .. (2.5,-1.9)
								.. controls (2.2,-1.8) and (2,-1.4) .. (1.6,-1.3)
								.. controls (1.3,-1.2) and (0.8,-1.4) .. (0.6,-1.1)
								.. controls (0.3,-0.8) and (0.15,-0.15) .. (0.3,0.3) -- cycle;
			\draw (rechts) circle (25pt)
						(links) circle (25pt);
						\filldraw	(links) circle (1pt) node[above] {$m_i$}
						(rechts) circle (1pt) node[above] {$n_j$}
								(a1)  circle (1pt) node[above] {$a_1$}
								(a2) circle (1pt) node[above] {$a_2$}
								(b1) circle (1pt) node[above] {$b_1$}
								(b2) circle (1pt) node[above] {$b_2$}
								(pm) circle (1pt) 
								(pm)+(-0.3,0) node[above] {$^pm_i$}
								(qm) circle (1pt) node[below] {$^qm_i$};
	\draw[dashed, ->] (qm) .. controls (4.7,-0.8) and (4.9,0.2) .. (4.5,1)
												.. controls (3.8,1.7) and (2.6,1.6) .. (pm);
												\draw (5.1,0) node {$_{pq\invers}$};
	\draw[->] (a1) -- (b1);
	\draw[->]	(a2) -- (b2);
	\draw[->]	(links) -- (pm);
	\draw[->]	(links) -- (qm);
\end{tikzpicture}
\end{center}

	F"ur $i,j\in H$ gibt es also h"ochstens ein $g_{i,j}\in G$ mit ${g_{i,j}}U_i\cap V_j\neq\emptyset$.
	W"ahle  $W=(U_h,V_h)_{h\in H}$  als Umgebung von $(m,n)$. 
	Sei $(f,a)\in Z$, also
	\begin{eqnarray*}
   &  ^{(f,a)}(m'_h)_{h\in H} = (n'_h)_{h\in H}  & \mbox{ f"ur ein }(m',n')\in W \\
		\Leftrightarrow &	f(h)m'_h = n'_{a\invers h} & \mbox{ f"ur alle }h\in H \\
		\Leftrightarrow & f(h) = g_{h,a\invers h} &  \mbox{ f"ur alle } h\in H.
\end{eqnarray*}
	Damit ist $f$ durch $a$ schon eindeutig bestimmt. Da $H$	endlich ist, muss also auch $Z$ endlich sein.
\end{beweis}

\subsection{Die *-Werkzeuge}
Zus"atzlich zu den Werkzeugen W1 bis W5 werden wir noch die folgenden Aussagen brauchen, die wir aber alle beweisen k"onnen, wenn wir die Werkzeuge W2 bis W5 voraussetzen.
\begin{enumerate}
	\item[*W6:] $FICwF$ gilt f"ur endliche Gruppen.
  \item[*W7:] $FICwF$ gilt f"ur abz"ahlbare, freie Gruppen.	
	\item[*W8:] Seien $G_1$ und $G_2$ Gruppen die $FICwF$ erf"ullen. Dann gilt $FICwF$ auch f"ur $G_1\times G_2$.
	\item[*W9:] Seien $G_1$ und $G_2$ abz"ahlbare Gruppen die $FICwF$ erf"ullen. Dann gilt $FICwF$ auch f"ur $G_1*G_2$.
\end{enumerate}

*W6 gilt sowieso in den meisten F"allen aufgrund von Lemma~\ref{lm:FICwF-endliche-Gruppen}, aber auch sonst ist *W6 einfach ein Spezialfall von W3. Die restlichen *-Werkzeugen werden wir in diesem Kapitel beweisen.
Au"serdem werden wir zeigen, dass wir die Aussagen W3, W4 und W5 ohne Einschr"ankung von $FIC$ auf $FICwF$ erweitern k"onnen:
\begin{enumerate}
	\item[*W3:] $FICwF$ gilt f"ur virtuell abelsche Gruppen.
	\item[*W4:] Sei $\{G_i\; | \; i\in I\}$ ein gerichtetes System von Gruppen und f"ur jedes $G_i$ gelte $FICwF$. Dann gilt $FICwF$ auch f"ur $\colim_{i\in I}G_i$. 
	\item[*W5:] Sei $p:L\rightarrow G$ ein Gruppenhomomorphismus und $FICwF$ gelte f"ur $G$. Ist $FICwF$ wahr f"ur $p\invers(V)$ f"ur jede virtuell zyklische Untergruppe $V$ von $G$, dann gilt $FICwF$ auch f"ur $L$.
\end{enumerate}

Wenn W3 gilt, erhalten wir *W3 folgenderma"sen.
\begin{lemma}[*W3]\label{FICwF-virtuell-abelsch}
	$FIC$ gelte f"ur virtuell abelsche Gruppen. Dann gilt auch $FICwF$ f"ur virtuell abelsche Gruppen.
	\end{lemma}
\begin{beweis}
	Sei $G$ eine virtuell abelsche Gruppe und $A$ eine abelsche Untergruppe von $G$ mit endlichem Index. Sei $H$ eine beliebige endliche Gruppe. $A^H$ ist auch wieder abelsch und hat endlichen Index in $G^H\!.$ $G^H$ wiederum hat endlichen Index in $G\wr H$.
	Es ist also auch $G\wr H$ virtuell abelsch und erf"ullt nach Voraussetzung $FIC$.
\end{beweis}

Ist die G"ultigkeit von W4 gegeben, k"onnen wir f"ur *W4 den Beweis von Lemma~\ref{thm:FICwF-limes} nutzen: Da wir mit W4 die Vertr"aglichkeit von $FIC$ mit Kolimiten voraussetzen, brauchen wir nicht mehr wie in Lemma~\ref{thm:FICwF-limes} eine streng stetige Homologietheorie fordern und der Rest des Beweises f"ur die Vertr"aglichkeit von $FICwF$ mit Kolimiten geht analog durch.
\rem{
\begin{lemma}[*W4]\label{thm:W4-*W4}
	Sei $FIC$ eine gefaserte Isomorphismusvermutung, f"ur die W4 gilt, dann gilt auch *W4, d.h. wenn $FIC$ mit Kolimiten vertr"aglich ist, dann ist auch $FICwF$ mit Kolimiten vertr"aglich.
\end{lemma}
\begin{beweis}
\end{beweis}}

F"ur die Erweiterung von W5 auf *W5, m"ussen wir zuerst *W7 und *W8 beweisen.
Wir ben"otigen daf"ur folgendes Hilfslemma:
\begin{lemma}\label{lm:FICwF-kompakte-flaechen}
  Sei $FICwF$ eine endlich erweiterbare, gefaserte Isomorphismusvermutung, f"ur die W2 und *W6 gilt.
	Sei $M$ eine kompakte, geschlossene Fl"ache. Dann wird $FICwF$ von $\pi_1(M)$ erf"ullt.
\end{lemma}
\begin{beweis}
	$S^2$ und $\R P^2$ erf"ullen mit trivialer bzw. mit zyklischer Fundamentalgruppe der Ordnung 2 aufgrund von *W6 $FICwF$ und alle anderen kompakten, geschlossene Fl"achen k"onnen mit einer nicht-positiven Metrik versehen werden, erf"ullen also $FICwF$ aufgrund von W2.
\end{beweis}

\begin{lemma}[*W7]\label{lm:FICwF-freieGruppen}
Sei $FICwF$ eine endlich erweiterbare, gefaserte Isomorphismusvermutung, f"ur die W2, W4 und *W6 gilt.
Dann ist $FICwF$ wahr f"ur freie Gruppen mit abz"ahlbar vielen Erzeugern, d.h. es gilt *W7.
\end{lemma}
\begin{beweis}
Sei $F=\left< a_1,\ldots,a_g \right>$ eine endlich erzeugte freie Gruppe mit $g$ Erzeugern, $g\geq 1$, und $M$ eine kompakte Fl"ache vom Geschlecht $g$. Nach dem letzten Lemma~\ref{lm:FICwF-kompakte-flaechen} ist $FICwF$ wahr f"ur 
\[
\pi_1(M)=\left< a_1,b_1,\ldots,a_g,b_g;\; [a_1,b_1]\cdots[a_g,b_g]=1 \right>.
\]
Aufgrund des Freiheitssatzes von Magnus (siehe beispielsweise \cite{magnus}, Theorem 4.10) ist $F$ eine frei erzeugte Untergruppe von $\pi_1(M)$ und erf"ullt somit nach Lemma~\ref{lm:FICwF-subgroup} $FICwF$.

Sei $F$ nun eine freie Gruppe mit abz"ahlbar vielen Erzeugern $x_i$, $i\in\N$. $F_n$ sei die von $\{x_1,\cdots,x_n\}$ erzeugte freie Gruppe. $F_i$ ist Untergruppe von $F_{i+1}$ und bez"uglich dieser Inklusionen gilt
\[
	F\cong\colim_{i\rightarrow\infty}{F_i}.
\]
Wie haben bereits gezeigt, dass $FICwF$ wahr ist f"ur alle $F_i$, also aufgrund von *W4 auch f"ur $\colim_{i\in\N}F_i=F$.
\end{beweis}

Im Folgenden sei jetzt $FICwF$ immer eine endlich erweiterbare, gefaserte Isomorphismusvermutung, f"ur die W2, W3, W4 und W5 gilt.
W1 wird erst wesentlich sp"ater eine Rolle spielen.

\begin{lemma}[*W8]\label{lm:FICwF-Kreuzprodukt}
  	Seien $G_1$ und $G_2$ Gruppen, f"ur die $FICwF$ gilt, dann wird $FICwF$ auch von $G_1 \times G_2$ erf"ullt, d.h.
		es gilt *W8.
\end{lemma}
\begin{beweis}
		Wir zeigen zun"achst, dass $FIC$ f"ur $G_1\times G_2$ erf"ullt ist und betrachten dazu die Projektion 
		\[
			p_1:G_1 \times G_2 \rightarrow G_1.
		\]
		Mit W5 gen"ugt es zu zeigen, dass f"ur jede virtuell zyklische Untergruppe $V_1$ von $G_1$ $FIC$ von $p_1^{-1}(V_1)=V_1 \times G_2$ erf"ullt wird. Dazu betrachten wir die Projektion
		\[
			p_2:V_1 \times G_2 \rightarrow G_2.
		\]
		Wiederum mit W5 gen"ugt es nun zu zeigen, dass $FIC$ von $p_2^{-1}(V_2)$ erf"ullt wird f"ur alle virtuell zyklischen Untergruppen $V_2$ von $G_2$. Es gilt $p_2^{-1}(V_2)=V_1 \times V_2$. Da $V_1$ und $V_2$ virtuell zyklisch sind, existieren zyklische Untergruppen $C_1\subset V_1$ und $C_2\subset V_2$ mit endlichem Index. Dann ist $C_1\times C_2$ eine Untergruppe von $V_1\times V_2$ mit endlichem Index, die insbesondere abelsch ist. Mit W3 gilt also $FIC$ f"ur $V_1\times V_2$ und damit auch f"ur $G_1\times G_2$.

  Sei $H$ eine beliebige endliche Gruppe. Nach Definition von $FICwF$ ist $FIC$ wahr f"ur $G_1\wr H$ und $G_2\wr H$. Es gilt also $FIC$ f"ur $(G_1\wr H) \times (G_2\wr H)$. Aufgrund von Lemma~\ref{lm:Kranzprodukt-Kreuzprodukt} ist $(G_1\times G_2)\wr H$  eine Untergruppe von $(G_1\wr H)\times (G_2\wr H)$. Mit Lemma~\ref{lm:FICwF-subgroup} "ubertr"agt sich $FIC$ auf Untergruppen und damit ist $FICwF$ wahr f"ur $G_1\times G_2$.
\end{beweis}

Mit *W8 k"onnen wir die Aussage von Lemma~\ref{lm:2.2} auf $FICwF$ erweitern und damit dann auch W5 auf *W5.

\begin{lemma}\label{lm:FICwF-zurueckziehen}
	Sei $FICwF$ eine endlich erweiterbare, gefaserte Isomorphismusvermutung, f"ur die *W8 gilt.
		Sei $p:K\rightarrow G$ ein Gruppenhomomorphismus. $FICwF$ sei erf"ullt f"ur $(G,\F_G)$ und f"ur $(p\invers(L),\F_{f\invers(L)})$ f"ur alle  $L\in\F_G$. Dann wird $FICwF$ auch von $(K,\F_K)$ erf"ullt.
\end{lemma}
\begin{beweis}
	Sei $H$ ein beliebige endliche Gruppe der Ordnung $n$. Dann induziert $p$ einen Gruppenhomomorphismus
	\[
		\widetilde{p}:K\wr H\rightarrow G\wr H, \qquad (k,h)\mapsto (f\circ k,h).
	\]

	Sei $V\in\F_{G\wr H}$. F"ur $G\wr H$ gilt nach Voraussetzung $FIC$. Um $FIC$ f"ur $(K\wr H)$ zu zeigen, reicht es aufgrund von Lemma~\ref{lm:2.2} $FIC$ f"ur $\widetilde{p}\invers(V)$ zu zeigen. 
	
	$G^n$ ist Normalteiler von $G\wr H$ mit endlichem Index.
	$V\cap G^n\in\F_{G\wr H}$. Sei $P_i$ die Projektion von $V\cap G^n$ auf die $i$-te Komponente. 
	Da $\F$ abgeschlossen unter Quotientenbildung ist, gilt $P_i\in\F_G$. Nach Voraussetzung gilt $FICwF$ f"ur $\widetilde{p}\invers(P_i)$, mit  *W8 gilt $FICwF$ dann f"ur $\prod_{1\leq i\leq n}\widetilde{p}\invers(P_i)$ und somit auch f"ur $V\cap G^n$, da sich mit Lemma~\ref{lm:FICwF-subgroup} $FICwF$ auf Untergruppen "ubertr"agt.

	Da $G^n$ als Normalteiler mit endlichem Index in $G\wr H$ liegt, ist $\widetilde{p}\invers(V\cap G^n)$ normal in $\widetilde{p}\invers(V)$ (f"ur $x\in\widetilde{p}\invers(v)$ und $y\in \widetilde{p}\invers(w)$ mit $v\in V$ und $w\in V\cap G^n$ gilt $\widetilde{p}(xyx\invers) = vwv\invers\in V\cap G^n$) und hat auch endlichen Index, denn es gilt
	\[
			|V/(V\cap G^n)|=|VG^n/G^n|\leq|(G\wr H)/G^n|<\infty.
			\]
	Wir erhalten also die exakte Sequenz
	\gexseq{\widetilde{p}\invers(V\cap G^n)}{\widetilde{p}\invers(V)}{V/(V\cap G^n)}{i}{pr\,\circ\,\widetilde{p}}{}
	und somit gilt aufgrund von Satz~\ref{thm:FICwF-virtuell-normal-FICwF} $FICwF$ f"ur $\widetilde{p}\invers(V)$ und insbesondere auch $FIC$, was zu zeigen war.
\end{beweis}

\begin{korollar}[*W5]\label{kr:FICwF-zurueckziehen}
	Sei $FICwF$ eine gefaserte Isomorphismusvermutung, f"ur die W5 und *W8 gilt. Dann gilt auch *W5, d.h. ist $p:K\rightarrow G$ ein Gruppenhomomorphismus und $FICwF$ gilt sowohl f"ur $G$ als auch f"ur $p\invers(V)$ f"ur alle virtuell zyklischen Untergruppe $V$ von $G$, dann gilt $FICwF$ auch f"ur $K$.
\end{korollar}
\begin{beweis}
	Die Voraussetzung von Lemma~\ref{lm:FICwF-zurueckziehen} sind erf"ullt. W"ahlen wir $V$ im Beweis von Lemma~\ref{lm:FICwF-zurueckziehen} als virtuell zyklische Untergruppe von $G\wr H$, dann k"onnen wir statt Lemma~\ref{lm:2.2} einfach W5 anwenden und der Beweis geht v"ollig analog durch.
\end{beweis}
Die Aussage von *W9 zeigen wir zun"achst f"ur $FIC$, bevor wir sie mit Satz~\ref{thm:FICwF-freiesProdukt} dann auch f"ur $FICwF$ beweisen. Der Aufbau des Beweises ist leicht abge"andert, aber sonst von Roushon "ubernommen \cite[Reduction Theorem und Theorem 3.1]{roushon}.

\begin{satz}\label{thm:reduktion}
  Seien $G_1$ und $G_2$ abz"ahlbare Gruppen, die $FIC$ erf"ullen. Dann wird $FIC$ auch von dem freien Produkt $G_1 * G_2$ erf"ullt.
\end{satz}

F"ur den Beweis ben"otigen wir Gruppenoperationen auf B"aumen. F"ur Details sei auf das Buch von Serre \cite{serre} verwiesen.
Wir greifen zun"achst auf folgenden Satz zur"uck.

\begin{satz}[{\cite[Theorem 7]{serre}}]\label{th:amalgam-tree}
	Sei $G=G_P*_A G_Q$ ein amalgamiertes Produkt von zwei Gruppen. Dann gibt es bis auf Isomorphie genau einen Baum $X$ auf dem $G$ operiert mit einem Segment
	$T=\xymatrix@M=0pt@R=0pt@L=2pt{ 
	{\bullet} \ar@{-}[r]_{\txt{\tiny a}} & {\bullet} \\ \txt{\tiny P} & \txt{\tiny Q} }$  
	als Fundamentalbereich, wobei $P$ von $G_P$, $Q$ von $G_Q$ und $a$ von $A$ in $X$ festgehalten wird.
\end{satz}
F"ur den Beweis von Satz~\ref{thm:reduktion} basteln wir daraus folgendes Lemma.
\begin{lemma}\label{lm:free-product-on-trees}
	Zu einem freien Produkt $G_P*G_Q$ gibt es einen Baum $X$ auf dem $G_P*G_Q$ operiert mit trivialen Kantenstandgruppen und die Eckenstandgruppen sind alle von der Gestalt $gG_Pg\invers$ oder $gG_Qg\invers$ f"ur gewisse $g\in G_P*G_Q$.
\end{lemma}
\begin{beweis}
	Nach dem letzten Satz gibt es einen Baum $X$ auf dem $G_P*G_Q$ operiert mit einem Segment $T=\xymatrix@M=0pt@R=0pt@L=2pt{ 
	{\bullet} \ar@{-}[r]_{\txt{\tiny a}} & {\bullet} \\ \txt{\tiny P} & \txt{\tiny Q} }$ als Fundamentalbereich. Zu jeder Ecke $E$ von $X$ gibt es daher ein $g\in G_P*G_Q$ mit $gP = E$ oder $gQ = E$ und somit ist jede Eckenstandgruppe zu $G_P$ oder $G_Q$ konjugiert, was man wie folgt einsieht: 
	
	Sei $G_E$ die Standgruppe einer Ecke $E$. Wir nehmen ohne Einschr"ankung an, dass $E$ zur Bahn von $P$ geh"ort. Der andere Fall verl"auft v"ollig analog.
	Es gibt also ein $g\in G_P*G_Q$ mit $gP = E$. F"ur alle $h\in G_E$ wird $P$ von $g\invers hg$ festgehalten und damit $h\in gG_Pg\invers$ und umgekehrt wird $E$ von jedem Element aus $gG_Pg\invers$ festgehalten, denn $gg_1g\invers E=gg_1P=gP=E$ f"ur $g_1\in G_P$. 
	
	$a$ wird nur von der trivialen Gruppe festgehalten. Mit einer analogen Argumentation, bei der lediglich $G_P$ durch die triviale Gruppe ersetzt wird, folgt, dass alle Kantenstandgruppen trivial sein m"ussen.
\end{beweis}

\begin{refbeweis}{Satz~\ref{thm:reduktion}}
  Wir betrachten die kanonische Projektion 
  \[
  		\p: G_1*G_2\rightarrow G_1\times G_2.
  \]
  Nach Voraussetzung ist $FIC$ wahr f"ur $G_1$ und $G_2$ also mit *W8 auch f"ur $G_1\times G_2$. Mit W5 reicht es nun zu zeigen, dass f"ur jede virtuell zyklische Untergruppe $V$ von $G_1 \times G_2$ $FIC$ f"ur das Urbild $\p\invers (V)\subset G_1 * G_2$ wahr ist.
  Wir zeigen mit dem n"achsten Lemma, dass es in dieser Situation eine abz"ahlbare freie Gruppe $F$ und eine endliche Gruppe $H$ gibt, so dass
  \[
   \p\invers(V) \subset F\wr H
  \]
  gilt.
	Mit Lemma~\ref{lm:FICwF-freieGruppen} erhalten wir $FICwF$ f"ur $F$.
	$FIC$ vererbt sich dann mit Lemma~\ref{lm:FICwF-subgroup} von $F\wr H$ auf die Untergruppe $\p\invers(V)$.
\end{refbeweis}  

\begin{lemma}\label{lm:fh}
Seien $G_1$ und $G_2$ abz"ahlbare Gruppen, $V$ eine virtuell zyklische Untergruppe von $G_1\times G_2$ und \[\p:G_1*G_2\rightarrow G_1\times G_2\] der kanonisch gegebene surjektive Gruppenhomomorphismus. Dann gibt es eine ab\-z"ahl\-bare freie Gruppe $F$ und eine endliche Gruppe $H$, so dass $\p\invers(V)$ Untergruppe von $F\wr H$ ist.
\end{lemma}
\begin{beweis}{}
	\begin{fallumgebung}{1}
		\fall{$V$ hat endliche Ordnung.} 
Wir betrachten folgende exakte Sequenz\exseq{\ker(\p)}{\p\invers(V)}{V}{.}
Nach Satz~\ref{thm:Kranzprodukt-Erweiterungen} ist $\p\invers(V)$ Untergruppe von $\ker(\p)\wr V$.
Nach Voraussetzung ist $V$ endlich. Bleibt zu zeigen, dass $\ker(\p)$ frei ist.
Wir benutzen Lemma~\ref{lm:free-product-on-trees}. Demnach operiert das freie Produkt $G_1*G_2$ auf einem Baum $\mathcal{T}$ mit trivialen Kantenstandgruppen und die Eckenstandgruppen haben die Gestalt $gG_1g\invers$ oder $gG_2g\invers$ f"ur bestimmte $g\in G_1*G_2$.
$\ker(\p)$ operiert als Untergruppe von $G_1*G_2$ ebenfalls auf $\mathcal{T}$ und zwar mit Eckenstandgruppen $\ker(\p)\cap gG_ig\invers$, $i\in\{1,2\}$.

Sei $ghg\invers\in \ker(\p)\cap gG_ig\invers$ ein Element, dass eine Ecke von $\mathcal{T}$ festh"alt und sei $(g_1,g_2)\coloneqq\pi(g)$. Wir nehmen $h\in G_1$ an. Der andere Fall geht analog. Es gilt  
\[
1=\p(ghg\invers)=\p(g)\p(h)\p(g\invers)=(g_1hg_1\invers,g_2g_{2}\invers).
\]
Aus $g_1hg_1\invers = 1$ folgt $h=1$ und damit sind alle Standgruppen $\ker(\p)\cap gG_ig\invers$ trivial. $\ker(\p)$ operiert also frei auf $\mathcal{T}$ und ist damit auch als Gruppe frei (\cite[Theorem 4]{serre}).

\fall{$V$ hat unendliche Ordnung.}
Da $V$ virtuell zyklisch ist, gibt es eine unendliche zyklische Untergruppe $C$ von $V$ mit endlichem Index. $C$ kann ohne Einschr"ankung als normale Untergruppe gew"ahlt werden. Wir zeigen in Lemma~\ref{lm:zyklische-normale-Untergruppe}, dass dies tats"achlich keine  Einschr"ankung ist.
Wir betrachten nun folgende exakte Sequenz von Gruppen
 
\gexseq{\p\invers(C)}{\p\invers(V)}{V/C}{i}{pr~\circ ~\p}{,} 
wobei $pr:V\rightarrow V/C$ die kanonische Projektion ist. Mit Satz~\ref{thm:Kranzprodukt-Erweiterungen} ist $\p\invers(V)$ Untergruppe von $\p\invers (C)\wr V/C$. $V/C$  ist endlich. Es bleibt lediglich zu zeigen, dass 
$\p\invers(C)$ frei ist.

Wir nehmen dazu wieder die Operation von $G_1*G_2$ auf einem Baum $\mathcal{T}$ zur Hilfe. Als Untergruppe von $G_1*G_2$ operiert auch $\p\invers(C)$ auf $\mathcal{T}$.  Wir haben bereits im ersten Fall gezeigt, dass $\ker(\p)\cap gG_ig\invers$ trivial ist. $\p$ eingeschr"ankt auf die Eckenstandgruppen $gG_ig\invers$ f"ur $i=1,2$ ist also injektiv und damit insbesondere auf $\p\invers(C)\cap gG_ig\invers$ injektiv. $\p\invers(C)\cap gG_ig\invers$ ist also isomorph zu einer Untergruppe von $C$ und damit entweder trivial oder wie $C$ ebenfalls unendlich zyklisch. $\p\invers(C)$ operiert also auf $\mathcal{T}$ mit trivialen Kantenstandgruppen und die Eckenstandgruppen sind alle entweder trivial oder frei mit einem Erzeuger. $\p\invers(C)$ ist damit ein freies Produkt von freien Gruppen (siehe beispielsweise \cite[VII, Theorem 2]{baumslag} oder \cite[5.4, Theorem 12]{serre}), also selbst eine freie Gruppe.  
\end{fallumgebung}
\end{beweis}

\begin{lemma}\label{lm:zyklische-normale-Untergruppe}
	Sei $V$ eine virtuell zyklische Gruppe. Dann besitzt $V$ einen zyklischen Normalteiler $N$ mit endlichem Index.
\end{lemma}
\begin{beweis}
Nach Voraussetzung existiert eine zyklische Untergruppe $C$ von $V$ mit endlichem Index. 
	$N:=\bigcap_{g\in V}gCg\invers$ ist ein Normalteiler von $V$. Da $V$ endlich erzeugt ist, gibt es nur endlich viele Untergruppen
	von $V$ mit Index $[V:C]$ (siehe beispielsweise \cite[Kapitel III, Theorem 3]{baumslag}). Es gibt also eine endliche Teilmenge $H\subset V$ mit
	$N=\bigcap_{g\in H}gCg\invers$. $N$ ist als endlicher Schnitt von zyklischen Untergruppen mit endlichem Index selbst wieder zyklisch mit endlichem Index. 
\end{beweis}

Wir haben nun die Vertr"aglichkeit von $FIC$ mit dem freien Produkt bewiesen. Jetzt fehlt sie noch f"ur $FICwF$.
\begin{satz}[*W9]\label{thm:FICwF-freiesProdukt}
	Sei $FICwF$ eine endlich erweiterbare, gefaserte Isomorphismusvermutung, f"ur die W2, W3, W4 und W5 gilt. 
	Seien $G_1$ und $G_2$ abz"ahlbare Gruppen, die $FICwF$ erf"ullen. Dann wird $FICwF$ auch von dem freien Produkt $G_1 * G_2$ erf"ullt.
\end{satz}
\begin{beweis}
	Sei $H$ eine endliche Gruppe. Die surjektive Abbildung 
	\[
		\p: G_1*G_2\rightarrow G_1\times G_2
	\]
	l"asst sich erweitern zu der surjektiven Abbildung
	\begin{eqnarray*}
		p: \qquad\quad(G_1*G_2)\wr H&\rightarrow& (G_1\times G_2)\wr H \\
			(f:H\rightarrow G_1*G_2, h)&\mapsto& (\p\circ f, h).
    \end{eqnarray*}

	Mit Lemma~\ref{lm:FICwF-Kreuzprodukt} wissen wir, dass $FICwF$ f"ur $G_1\times G_2$, also $FIC$ f"ur $(G_1\times G_2)\wr H$ gilt. Wir betrachten eine beliebige virtuell zyklische Untergruppe $V$ von $(G_1\times G_2)\wr H$ und zeigen, dass das $FIC$ von dem  Urbild $p\invers(V)$ erf"ullt wird. Da wir W5 gefordert hatten, ist dann auch $FIC$ f"ur $(G_1* G_2)\wr H$ wahr.

		In $V$ finden wir einen zyklischen Normalteiler $C$ mit endlichem Index (Lemma~\ref{lm:zyklische-normale-Untergruppe}). Sei $\gamma=(f,h)$ der Erzeuger von $C$. Da $H$ endlich ist, gibt es ein $k\in\N$ und ein $g\in (G_1\times G_2)^H$ mit $\gamma^k=(g,1)$. 
		Setze 
		\[
			C'\coloneqq \left< \gamma^k\right>=C\cap((G_1\times G_2)^H\times \{1\}).
		\]
		$C'$ hat endlichen Index in $C$ und damit auch in $V$. $C$ ist normal in $V$ und\\ $(G_1\times G_2)^H\times\{1\}$ normal in $(G_1\times G_2)\wr H$ und damit ist $C'$ normal in $V$. Das liefert uns folgende exakte Sequenz von Gruppen:
		\begin{eqnarray}
  			1\rightarrow p\invers(C')\rightarrow p\invers(V)\rightarrow V/C'\rightarrow 1.\label{exseq}
			\end{eqnarray}
			Wir betrachten $C'$ nun als Untergruppe von $(G_1\times G_2)^{|H|}\cong(G_1\times G_2)^H\times\{1\}$. Sei $(\gamma_1,\ldots,\gamma_n), n\coloneqq |H|$ der Erzeuger von $C'$. Es gilt
	
		\[
		p\invers(C') \subset \p\invers(\left< \gamma_1 \right>)\times\p\invers(\left< \gamma_2 \right>)\times \dots \times\p\invers(\left< \gamma_n \right>).
		\]
	Nach Lemma~\ref{lm:fh} ist f"ur jedes $i\in\{1,\dots,n\}$ die Gruppe $\p\invers(\left< \gamma_i \right>)$ eine Untergruppe eines Kranzproduktes $F_i\wr H_i$ mit einer abz"ahlbaren freien Gruppe $F_i$ und einer endlichen Gruppe $H_i$.
	Mit der exakten Sequenz \ref{exseq} erhalten wir
	\begin{eqnarray*}
	  p\invers(V) &\stackrel{\ref{thm:Kranzprodukt-Erweiterungen}}{\subset} &	p\invers(C')\wr(V/C') \\
		& \stackrel{\ref{lm:Kranzprodukt-Untergruppen}}{\subset} &  \left[(F_1\wr H_1)\times(F_2\wr H_2)\times \dots\times(F_n\wr H_n)\right] \wr (V/C') \\
		& \stackrel{\ref{lm:Kranzprodukt-Kreuzprodukt}}{\subset} & \left[ (F_1\wr H_1)\wr (V/C) \right]\times \dots \times\left[ (F_n\wr H_n)\wr (V/C')\right] \\
		& \stackrel{\ref{lm:Kranzprodukt-Einbettung}}{\subset} &
		F_1\wr(H_1\wr(V/C))\times \dots \times F_n\wr(H_n\wr(V/C')).
	\end{eqnarray*}
	Jede Gruppe $F_i\wr(H_i\wr(V/C'))$ erf"ullt nach Lemma~\ref{lm:FICwF-freieGruppen} $FIC$, da die $F_i$ frei und die $H_i\wr (V/C')$ endlich sind. F"ur 
	$F_1\wr(H_1\wr(V/C))\times \dots \times F_n\wr(H_n\wr(V/C'))$ gilt also $FIC$ aufgrund von Lemma~\ref{lm:FICwF-Kreuzprodukt} und mit Lemma~\ref{lm:FICwF-subgroup} "ubertr"agt sich $FIC$ auf die Untergruppe $p\invers(V)$.
\end{beweis}
Wir haben somit in diesem Kapitel gezeigt, dass wenn f"ur eine endlich erweiterbare, gefaserte Isomorphismusvermutung W1 bis W5 gilt, dann gilt auch *W3 bis *W9.
\chapter{Der rote Baum}
\section{Die Wurzel}
$FIC$ sei von nun an eine "`gefaserte Isomorphismusvermutung mit Werkzeugkasten"', d.h. eine gefaserte Isomorphismusvermutung f"ur die W1 bis W5 gilt und damit auch *W3 bis *W9. Entsprechend sei $FICwF$ eine "`endlich erweiterbare gefaserte Isomorphismusvermutung mit Werkzeugkasten"'.

Mit unserem Werkzeugkasten ausger"ustet, wagen wir uns nun an den Beweis f"ur unseren
\begin{zielsatz}
	Sei $M$ eine 3"=Mannigfaltigkeit und $FICwF$ eine endlich erweiterbare, gefaserte Isomorphismusvermutung mit Werkzeugkasten. Dann gilt $FICwF$ f"ur $\pi_1(M)$.
\end{zielsatz}
Die Suche nach einem "`roten Faden"', der die vor allem in \cite{roushon} und \cite{roushon-b-gruppen} verstreuten Baukl"otzchen zu einem nachvollziehbaren Beweis aneinanderreiht, hat einen recht verzweigten "`roten Baum"' hervorgebracht. Zur besseren "Ubersicht wird dieser rote Baum als Diagramm parallel zur Beweisf"uhrung mitentwickelt. Dabei zeigen einfache Verzweigungen immer Fallunterscheidungen an, die aus der Struktur der 3"=Mannigfaltigkeiten gewonnen werden. Durchgezogene Umrandungen kennzeichnen Typen von 3"=Mannigfaltigkeiten, f"ur die wir bereits bewiesen haben, dass $FICwF$ gilt, f"ur gestrichelt umrandete F"alle m"ussen wir den Beweis noch f"uhren. Doppelpfeile machen deutlich, welche F"alle in Beweisen f"ur andere F"alle ben"otigt werden und wo welche Lemmata und S"atze eingehen. Wir werden am Ende sehen, dass wir nur knapp an einem Ringschluss vorbeischrammen. 

Der erste Schritt auf unserem Weg zum Zielsatz besteht darin zu zeigen, dass wir nicht an beliebigen 3"=Mannigfaltigkeiten herumklempnern m"ussen, sondern wir uns auf orientierbare, kompakte, irreduzible 3-Mannigfaltigkeiten beschr"anken k"onnen. Das folgende erste Baumdiagramm verdeutlicht, wie wir dabei vorgehen werden. Anschlie"send haben wir dann zwei gr"o"sere Baukl"otze vor uns, die wir getrennt behandeln werden: Orientierbare, kompakte, irreduzible 3"=Mannigfaltigkeiten mit Rand und ohne Rand.

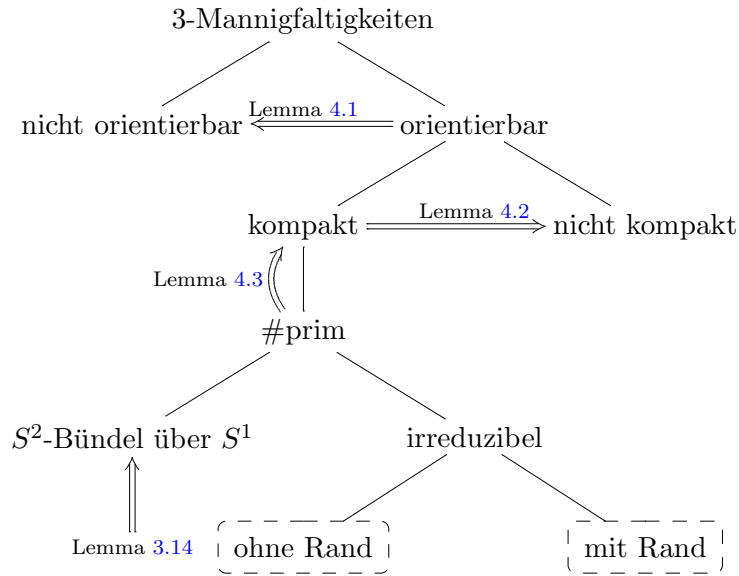
\begin{align*}
	\xymatrix@!C=40pt{
	& {\txt{3"=Mannigfaltigkeiten}}
		\ar@{-}[dl]
		\ar@{-}[dr]  
	& &  
	\\ 
		{\text{nicht orientierbar}} 
	& &  
		{\text{orientierbar}}
		\ar@{-}[dl] 
		\ar@{-}[dr]\ar@{=>}[ll]_{\text{Lemma~\ref{lm:orientierbar=>nichtorientierbar}}} 
	&   
	\\
	&  
		{\text{kompakt}} 
		\ar@{-}[d]
		\ar@{=>}[rr]^{\text{Lemma~\ref{lm:mit-Rand=>nicht-kompakt}}}
	& & 
		{\text{nicht kompakt}}  
	\\
	&  
		{\text{\#prim}}  
		\ar@{=>}@/^1pc/[u]^{\text{Lemma~\ref{lm:Zerlegung-Mfkgruppe-in-freies-Produkt}}} 
		\ar@{-}[dl]
		\ar@{-}[dr]  
		& & 
	\\ 
		{S^2\text{-B"undel "uber }S^1}
	& & 
		{\text{irreduzibel}} 
		\ar@{-}[dl]
		\ar@{-}[dr] 
	&   
	\\ 
	  *+{\text{\scriptsize Lemma~\ref{lm:FICwF-freieGruppen}}}
		\ar@{=>}[u]
	& 
		*++[F--:<3pt>]{\text{ohne Rand}} 
	& & 
		*++[F--:<3pt>]{\text{mit Rand}}
	}
\end{align*}

Mit der Sprechweise, dass "`$FIC(wF)$ f"ur eine Mannigfaltigkeit"' gilt, meinen wir im Folgenden immer, dass $FIC(wF)$ f"ur die Fundamentalgruppe der Mannigfaltigkeit wahr ist.

\begin{lemma}\label{lm:orientierbar=>nichtorientierbar}
  $FICwF$ sei f"ur orientierbare 3"=Mannigfaltigkeiten wahr. Dann gilt $FICwF$ schon f"ur beliebige 3"=Mannigfaltigkeiten.
\end{lemma}

\begin{beweis}
  	Jede nicht-orientierbare 3"=Mannigfaltigkeit $M$ besitzt eine 2-bl"attrige orientierbare "Uberlagerung $\widetilde{M}$, die man "uber die von allen orientierungserhaltenden Wegen erzeugte Untergruppe von $\pi_1(M)$ erh"alt. F"ur Details sei auf das 3. Beispiel in \cite[\S75]{seifert} verwiesen.
		$\pi_1(\widetilde{M})$ hat Index 2 in $\pi_1(M)$, ist also Normalteiler. Gilt $FICwF$ f"ur $\widetilde{M}$, dann "ubertr"agt sich mit Satz~\ref{thm:FICwF-virtuell-normal-FICwF} $FICwF$ auch auf die "uberlagerte Mannigfaltigkeit $M$.
  \end{beweis}

Im Folgenden sind deshalb, auch wenn wir es nicht mehr dazu schreiben, ohne Ausnahme alle Mannigfaltigkeiten, die wir betrachten, orientierbar und au"serdem, aufgrund der beiden folgenden Lemmata, in aller Regel kompakt und irreduzibel.

\begin{lemma}\label{lm:mit-Rand=>nicht-kompakt}
	$FICwF$ sei f"ur kompakte Mannigfaltigkeiten mit Rand wahr. Dann gilt $FICwF$ schon f"ur nicht-kompakte 3"=Mannigfaltigkeiten.
\end{lemma}
\begin{beweis}
	Es gibt ein $n\in\N$ und eine Einbettung $i:M\rightarrow \R^n$, so dass $0\notin i(M)$. Sei $d$ die Standardmetrik auf $\R^n$. Die Verkn"upfung 
	\[
		M\xrightarrow{i}\R^n\xrightarrow{d(-,0)}\R
	\]
	ist differenzierbar	und die Menge der regul"aren Werte von $i\circ d(-,0)$  liegt somit dicht in $\R$ (Satz von Sard, siehe beispielsweise \cite[Satz 6.1]{broecker}).
	Wir finden somit eine bez"uglich Inklusion aufsteigende Folge von kompakten Mannigfaltigkeiten mit Rand 
	\[(M_r)_{r\in I}\coloneqq	i\invers(B_r(0)\cap i(M))_{r\in I},\] wobei $B_r(0)$ der abgeschlossene Ball um $0$ mit Radius $r$ ist, die uns eine Folge von Fundamentalgruppen $(\pi_1(M_r))_{r\in I}$ liefert mit $\colim_{r\in I} \pi_1(M_r) = \pi_1(M)$. 
	Gilt also $FICwF$ f"ur kompakte 3"=Mannigfaltigkeiten mit Rand, dann ist mit *W4 auch $FICwF$ f"ur nicht-kompakte 3-Mannig\-faltig\-keiten wahr.
\end{beweis}

\begin{lemma}\label{lm:Zerlegung-Mfkgruppe-in-freies-Produkt}
	Sei $M$ eine kompakte 3"=Mannigfaltigkeit.
	Dann gilt
	\[
	\pi_1(M) = \pi_1(M_1)*\ldots*\pi_1(M_n)*F,
	\]
	wobei die $M_i$ irreduzible, kompakte 3"=Mannigfaltigkeiten sind und $F$ eine freie Gruppen von endlichem Rang ist.
\end{lemma}
\begin{beweis}
  $M$ l"asst sich als endliche zusammenh"angende Summe $M_1 \# \dots \# M_k = M$ von Prim"=Mannigfaltigkeiten schreiben (Satz~\ref{thm:primzerlegung}).	Weiterhin ist jede Prim"=Mannigfaltigkeit entweder irreduzibel oder ein $S^2$-B"undel "uber $S^1$ (Lemma~\ref{lm:prim-irr}). 
Ist $M_i$ ein $S^2$-B"undel "uber $S^1$, erhalten wir mit der langen exakten Homotopiesequenz $\pi_1(M_i)\cong\pi_1(S^1)\cong\Z$, da die Fasern zusammenh"angend sind und triviale Fundamentalgruppe haben. Mit dem Satz von Seifert und van Kampen (siehe beispielsweise \cite[IV.\S3, Theorem 3.1]{seifert-vankampen}) erhalten wir $\pi_1(M)$ aus der zusammenh"angenden Summe der Mannigfaltigkeiten $M_i$ als ein freies Produkt der Fundamentalgruppen $\pi_1(M_i)$, wobei sich die Fundamentalgruppen der nicht irreduziblen $M_i$ zu einer freien Gruppe $F$ von endlichem Rang zusammenfassen lassen.
\end{beweis}
Wir wissen mit *W9 (Satz~\ref{thm:FICwF-freiesProdukt}) bereits, dass $FICwF$ mit dem freien Produkt vertr"aglich ist und mit *W7 (Lemma~\ref{lm:FICwF-freieGruppen}) hatten wir gezeigt, dass $FICwF$ f"ur abz"ahlbare, freie Gruppen wahr ist. Mit Lemma~\ref{lm:Zerlegung-Mfkgruppe-in-freies-Produkt} brauchen wir $FICwF$ also lediglich noch f"ur kompakte, irreduzible 3"=Mannigfaltigkeits\-gruppen zu zeigen, um $FICwF$ f"ur beliebige 3"=Mannigfaltigkeits\-gruppen zu bekommen.
Halten wir also unser erstes Etappenziel fest:
\begin{etappe}\label{et:kompakt-irreduzibel}
	Gilt $FICwF$ f"ur orientierbare, kompakte, irreduzible 3"=Mannigfaltigkeiten, dann gilt $FICwF$ schon f"ur beliebige 3"=Mannigfaltigkeiten.
\end{etappe}

\section{Ohne-Rand-Ast}
Wir schauen uns zun"achst kompakte, irreduzible 3"=Mannigfaltigkeiten ohne Rand genauer an.

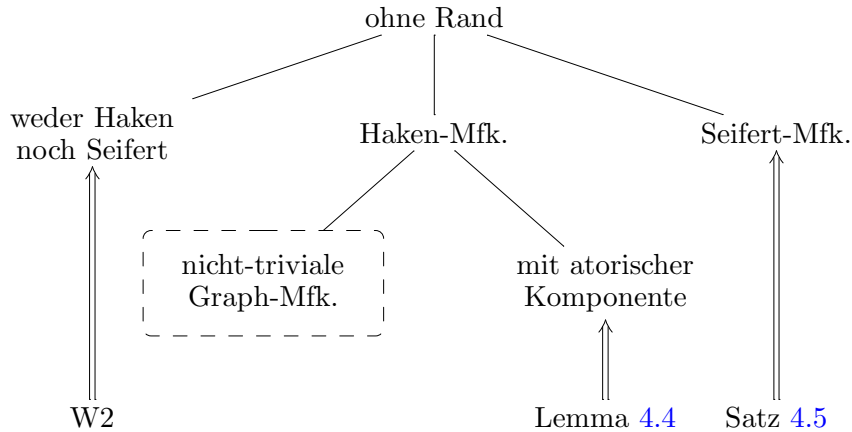
\begin{align*}
	\xymatrix@!C=40pt{
	& & 
		{\text{ohne Rand}}
		\ar@{-}[dll]
		\ar@{-}[d]
		\ar@{-}[drr]
	&  &
	\\
		{\txt<6pc>{weder Haken\\noch Seifert}} 
	& &
		{\txt{Haken-Mfk.}}
		\ar@{-}[dl]
		\ar@{-}[dr] 
	& &
		{\txt<5pc>{Seifert-Mfk.}} 
	\\
	&
		*+++[F--:<3pt>]{\txt<6pc>{nicht-triviale Graph-Mfk.}}
	& &
		{\txt<6pc>{mit atorischer Komponente}}
	&
	\\
		{\text{W2}}
		\ar@{=>}[uu]
	& & &
		{\text{Lemma~\ref{lm:FICwF-atorische-Komponente}}}
		\ar@{=>}[u]
	&
		{\text{Satz~\ref{thm:FICwF-Seifert}}}
		\ar@{=>}[uu]
			\\
	}
\end{align*}

F"ur 3"=Mannigfaltigkeiten mit endlicher Fundamentalgruppe hatten wir mit *W6 gefordert, dass $FICwF$ gilt. 
Im Folgenden ist also ohne Einschr"ankung $\pi_1(M)$ unendlich. F"ur kompakte, orientierbare, irreduzible 3-Mannig\-faltig\-keiten $M$ mit unendlicher Fundamentalgruppe, die weder Haken- noch Seifert-Mannig\-faltig\-keiten sind, besagt die Geometrisierungsvermutung von Thurston (siehe \cite[\S6]{scott}), dass man eine hyperbolische Metrik findet. Dank des Ergebnisses von Perelman \cite{perelman} k"onnen wir also auf 3-Mannig\-faltig\-keiten, die weder Seifert- noch Haken-Mannigfaltigkeit sind, W2 anwenden.

F"ur Haken-Mannigfaltigkeiten machen wir die weitere Unterscheidung, ob die Toruszerlegung eine atorische Komponente besitzt oder nur aus Seifert-Komponenten besteht. 
F"ur den ersten Fall beweisen wir folgendes Lemma:
\begin{lemma}[{\cite[Korollar 4.1]{roushon}}]\label{lm:FICwF-atorische-Komponente}
	Sei $M$ eine geschlossene Haken-3"=Mannigfaltigkeit, deren Toruszerlegung eine atorische Komponente besitzt. Dann wird $FICwF$ von $\pi_1(M)$ erf"ullt.
\end{lemma}

\begin{beweis}
	Leeb hat mit Satz 3.3 in \cite{leeb} gezeigt, dass sich f"ur $M$ eine nicht-positiv gekr"ummte Metrik w"ahlen. Mit W2 gilt also $FICwF$ f"ur $\pi_1(M)$.
\end{beweis}

F"ur den Fall der geschlossenen, nicht-trivialen Graph-Mannigfaltigkeit ben"otigen wir $FICwF$ f"ur 3"=Mannigfaltigkeiten mit Rand, wof"ur wir wiederum $FICwF$ f"ur geschlossene Seifert-Mannigfaltigkeiten brauchen. Denen widmen wir uns deshalb mit dem n"achsten Satz.

\rem{\subsection{Seifert-Mannigfaltigkeiten}}
\begin{satz}\label{thm:FICwF-Seifert}
	Sei $S$ eine geschlossene Seifert-Mannigfaltigkeit. Dann ist $FICwF$ erf"ullt f"ur $\pi_1(S)$.
\end{satz}

\begin{beweis}
		Wir lehnen uns an den Beweis an, den Roushon in \cite[Theorem 4.6]{roushon} f"ur beliebige Seifert-Mannigfaltigkeiten gibt, wandeln ihn aber so ab, dass wir als Voraussetzung $FICwF$ lediglich f"ur virtuell abelsche statt f"ur virtuelle poly-$\Z$-Gruppen ben"otigen. 

	Wir erinnern an die Pr"asentation, die wir in Satz~\ref{thm:seifert-presentation} f"ur Seifert-Mannigfaltigkeiten gesehen haben.
	Da $S$ keine Randkomponenten besitzt, gilt
  \begin{eqnarray*}
  	\pi_1(S) &= \quad\left< \right. &a_1,b_1,\ldots,a_g,b_g,c_1,\ldots,c_q, t; \\
		& & a_ita_i\invers =t^{\varepsilon_i}, b_itb_i\invers = t^{\delta_i}, c_jtc_j\invers = t^{\eta_j}, c_j^{\alpha_j} = t^{\beta_j},\\
		& & \left. c_q=[a_1,b_1]\ldots [a_g,b_g]c_1\ldots c_{q-1}\quad\right>
  \end{eqnarray*}
  oder
  \begin{eqnarray*}
 \pi_1(S) &= \quad\langle &a_1,\ldots,a_g,c_1,\ldots,c_q, t; \\
 & &a_ita_i\invers =t^{\varepsilon_i}, c_jtc_j\invers = t^{\delta_j}, c_j^{n_j} = t^{s_j},\hspace{60pt}\\
 & &c_q=a_1^2\ldots a_g^2c_1\ldots c_{q-1} \quad\rangle.
  \end{eqnarray*}
  
	Offensichtlich ist die von $t$ erzeugte zyklische Untergruppe $\langle t\rangle$ normal in $\pi_1(S)$ und liefert uns somit folgende exakte Sequenz:
	\gexseq{\left< t\right>}{\pi_1(S)}{\pi_1(S)/\left< t \right>}{}{p}{,}
	wobei 
	 \begin{eqnarray*}
  	\pi_1(S)/\left< t\right> & = \quad \left< \right.&a_1,b_1,\ldots,a_g,b_g,c_1,\ldots,c_q;\\
		& & \left. c_1^{n_1}=\dots =c_q^{n_q} = 1, c_q=[a_1,b_1]\ldots [a_g,b_g]c_1\ldots c_{q-1}\quad\right>
  \end{eqnarray*}
  oder
  \begin{eqnarray*}
 \pi_1(S)/\left< t \right>&= \quad\left< \right. &a_1,\ldots,a_g,c_1,\ldots,c_q; \\
 & &  \left. c_1^{n_1} = \dots = c_q^{n_q} = 1, c_q=a_1^2\ldots a_g^2c_1\ldots c_{q-1} \quad\right>\hspace{40pt}
	\end{eqnarray*}
	die Gestalt einer Orbifold-Fundamentalgruppe hat. An der Pr"asentation sieht man schon, dass $\pi_1(S)/\langle t\rangle$ die Fundamentalgruppe einer geschlossenen Fl"ache enth"alt. Tats"achlich findet man auch eine Fl"achengruppe $G$, die endlichen Index in $\pi_1(S)/\langle t\rangle$ hat. Hempel beweist das beispielsweise in \cite[Theorem 12.2(i)]{hempel} oder auch in Bemerkungen verstreut zu finden bei Scott \cite{scott}, die Martino in \cite[Theorem 2.2]{martino} zusammengefasst hat.
	Da mit *W6 $FICwF$ f"ur endliche Gruppen gilt, k"onnen wir $\pi_1(S)$ als unendlich und damit $\left< t \right>$ als unendlich zyklisch annehmen\footnote{Dies ist f"ur den Beweis nicht zwingend erforderlich, erspart aber ein paar Fallunterscheidungen}.\rem{(\cite[VI.11]{jaco})}
	\rem{
	Au"serdem, wenn $\pi_1(S)$ eine unendliche Gruppe ist, enth"alt $\pi_1(S)/\left< t\right>$ eine Untergruppe $G$ mit endlichem Index, die isomorph zur Fundamentalgruppe einer geschlossenen Fl"ache. (siehe \cite[12.2]{hempel}).} Mit Lemma~\ref{lm:FICwF-kompakte-flaechen} ist $FICwF$ f"ur $G$ erf"ullt und da $\pi_1(S)/\left< t\right>$ endlich erzeugt ist, gilt aufgrund von Korollar~\ref{kr:FICwF-endl-erz-virtuell-FICwF} $FICwF$ auch f"ur $\pi_1(S)/\left< t\right>$.

	Sei $V$ eine virtuell zyklische Untergruppe von $\pi_1(S)/\left< t\right>$. $V$ besitzt eine normale zyklische Untergruppe $C$ mit endlichem Index in $V$ (Lemma~\ref{lm:zyklische-normale-Untergruppe}), die uns folgende exakte Sequenz liefert.
\gexseq{\left< t\right>}{p\invers (C)}{C}{}{p}{.}
Ist $C$ endlich, dann ist $\left< t\right>$ ein Normalteiler mit endlichem Index in $p\invers(C)$ und erf"ullt als freie Gruppe dar"uber hinaus $FICwF$. Mit Satz~\ref{thm:FICwF-virtuell-normal-FICwF} gilt also $FICwF$ auch f"ur $p\invers(C)$. 

Ist $C$ unendlich zyklisch, dann gibt es einen Schnitt $s:C\rightarrow p\invers(C)$ und aufgrund von Lemma~\ref{lm:semidirekt-spaltende-sequenz} gilt $p\invers (C)\cong \left< t\right>\sdp C\cong\Z\sdp\Z$. Somit wirkt $\left< t\right>$ auf $C$ entweder trivial oder mit $t^nc\mapsto (-1)^nc$. Im ersten Fall ist $\left< t\right>\sdp C\cong\Z\times\Z$ und erf"ullt damit als abelsche Gruppe $FICwF$ nach Voraussetzung *W3. Im zweiten Fall wirkt $\left< t^2\right>$ trivial auf $C$ und wir haben mit $\left< t^2\right>\times C$ einen abelschen Normalteiler mit endlichem Index von $\left< t\right>\sdp C$. Es ist damit ebenfalls *W3 anwendbar.

Mit
\exseq{p\invers (C)}{p\invers (V)}{V/C}{}
sehen wir, dass wir mit $p\invers(C)$ eine normale Untergruppe von $p\invers(V)$ mit endlichem Index gefunden haben, die $FICwF$ erf"ullt.  Es kommt erneut Satz~\ref{thm:FICwF-virtuell-normal-FICwF} zum Einsatz, womit sich $FICwF$ zun"achst auf $p\invers(V)$ "ubertr"agt und mit *W5 dann weiter auf $\pi_1(S)$.
\end{beweis}

Mit dem gerade Gezeigten und Etappe~\ref{et:kompakt-irreduzibel} k"onnen wir also folgenden Stand festhalten.
\begin{etappe}\label{et:mit-rand-und-graph}
  $FICwF$ gilt f"ur geschlossene Haken-Mannigfaltigkeiten mit atorischer Komponente und f"ur geschlossene Seifert-Mannigfaltigkeiten. 
	Ist $FICwF$ f"ur kompakte, irreduzible 3"=Mannigfaltigkeiten mit Rand und f"ur nicht-triviale Graph-Mannig\-faltigkeiten ohne Rand erf"ullt, dann gilt $FICwF$ schon f"ur beliebige 3"=Mannigfaltigkeiten.
\end{etappe}

\pagebreak
\section{Mit-Rand-Ast}
Mit Lemma~\ref{lm:Zerlegung-Mfkgruppe-in-freies-Produkt} hatten wir gesehen, dass wir kompakte 3-Mannigfaltigkeiten so zerlegen k"onnen, dass wir, bis auf eine zus"atzliche endlich erzeugte freie Gruppe, die Fundamentalgruppe als freies Produkt von Fundamentalgruppen von irreduziblen 3-Mannigfaltigkeiten erhalten. Wir werden jetzt zeigen, dass wir irreduzible 3-Mannigfaltigkeiten mit komprimierbarem Rand noch weiter zerlegen k"onnen, so dass wir, wieder bis auf eine zus"atzliche endlich erzeugte freie Gruppe, die Fundamentalgruppe als freies Produkt von Fundamentalgruppen von irreduziblen Mannigfaltigkeiten mit unkomprimierbarem Rand bekommen.
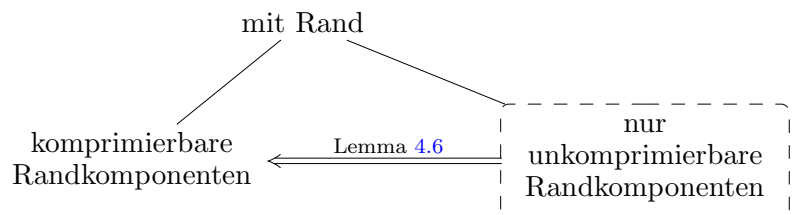
\begin{align*}
  \xymatrix@!C=40pt{
 	&
		{\text{mit Rand}}
		\ar@{-}[dl]
		\ar@{-}[drr]
	& & 
	\\  
		{\txt<8pc>{kompri\-mierbare\\Rand\-komponenten}} 
	& & &
		*++[F--:<3pt>]{\txt<8pc>{nur unkompri\-mierbare\\Rand\-komponenten}}
		\ar@{=>}[lll]_{\text{Lemma~\ref{lm:Zerlegung-Mfkgruppe-unkomprimierbar}}}
	}
\end{align*}
\begin{lemma}\label{lm:Zerlegung-Mfkgruppe-unkomprimierbar}
	Sei $M$ eine kompakte, zusammenh"angende, nicht einfach-zusammen\-h"ang\-ende, irreduzible 3"=Mannigfaltigkeit mit Rand. Dann gilt
	\[
		\pi_1(M)\cong \pi_1(M_1) * \ldots * \pi_1(M_k) * F,
	\]
	wobei die $M_i, i=1,\ldots,k$ kompakte irreduzible 3"=Mannigfaltigkeiten mit unkomprimierbarem Rand sind und $F$ eine freie Gruppe mit endlichem Rang ist.
      \end{lemma}
			Wir benutzen f"ur den Beweis dieses Lemmas den Schleifensatz von Papakyria\-kopoulos~\cite[Theorem 1]{papa} in einer von Stallings erweiterten Fassung \cite[Generalized loop theorem]{stallings}.
	\begin{satz}[Schleifensatz]\label{thm:looptheorem}
	  	Sei $M$ eine kompakte 3"=Mannigfaltigkeit mit Rand und $B$ eine Randkomponente von $M$. Sei $G$ eine normale Untergruppe von $\pi_1(B)$. Gibt es im Kern der Abbildung $\pi_1(B)\rightarrow\pi_1(M)$ Elemente, die nicht in $G$ liegen, dann existiert ein einfacher geschlossener Weg $\gamma$ in $B$, der eine in $M$ eingebettete 2-Disk berandet, so dass $[\gamma]\notin G$.
	\end{satz}
	W"ahlen wir $G=1$, erhalten wir folgendes Korollar:
	\begin{korollar}\label{kr:schleifensatzkorollar}
	  Sei $M$ eine kompakte 3"=Mannigfaltigkeit mit Rand und $B$ eine Randkomponente von $M$. Gibt es einen einfachen geschlossenen, nicht nullhomotopen Weg $\gamma$ in $B$, der in $M$ nullhomotop ist, dann berandet $\gamma$ in $M$ eine eingebettete 2-Disk.
	\end{korollar}
 
	Die Idee f"ur Lemma~\ref{lm:Zerlegung-Mfkgruppe-unkomprimierbar} ist, die Mannigfaltigkeit an komprimierbaren 2-Disks, die uns Korollar~\ref{kr:schleifensatzkorollar} liefert, aufzutrennen, bis alle Randfl"achen unkomprimierbar sind. 

	Beispielsweise f"ur eine $D^3$ aus der ein Volltorus und Vollfl"ache vom Geschlecht 2 entfernt wurden, die ineinander h"angen\footnote{Das ist zwar eigentlich kein zul"assiges Beispiel, weil diese Mannigfaltigkeit nicht irreduzibel ist, aber es veranschaulicht den Vorgang sehr gut.}, findet man eine komprimierbare Disk $D$, die eine Zerlegung in folgendes freies Produkt erm"oglicht:
\schneekugeln

	Damit wir zeigen k"onnen, dass tats"achlich jede Zusammenhangskomponente $M_i$, die dabei entsteht, einen Rand hat, ben"otigen wir f"ur den Beweis von Lemma~\ref{lm:Zerlegung-Mfkgruppe-unkomprimierbar} noch das folgende Lemma.
	In einem sp"ateren Abschnitt werden wir die Aussage von Lemma~\ref{lm:Rand-von-Untermfk} auch f"ur nicht-kompakte Mannigfaltigkeiten brauchen, weshalb wir sie hier gleich mitbeweisen.

\begin{lemma}\label{lm:Rand-von-Untermfk} 
	Sei $M$ eine nicht einfach-zusammenh"angende, irreduzible 3"=Mannigfaltigkeit, die entweder kompakt mit Rand oder nicht-kompakt ist.
	Dann besitzt jede kompakte, zusammenh"angende, aber nicht einfach-zusammenh"angende 3-dimensio\-nale Untermannigfaltigkeit $N$ von $M$ mindestens eine Randkomponente vom Geschlecht $\geq 1$. 
\end{lemma}
\begin{beweis}
Angenommen alle Randkomponenten von $N$ sind vom Geschlecht $0$, also zweidimensionale Sph"aren $S_1,\dots,S_n$. Da $M$ irreduzibel ist, berandet jede der $S_i$ eine eingebette 3-Disk $D_i\subset M$.
	$M\smallsetminus S_i$ ist nicht zusammenh"angend f"ur alle $i\in\{1,\ldots,n\}$, denn andernfalls g"abe es im Widerspruch dazu, dass $M$ irreduzibel ist, ein $S^2$-B"undel "uber $S^1$ in $M$ (\cite[Lemma 3.8]{hempel}). 	
	Da $N$ zusammenh"angend ist, gilt entweder $N\subseteq D_i$ oder $N\subseteq M\smallsetminus \mathring{D_i}$ f"ur alle $i\in\{1,\ldots,n\}$.	

	Nehmen wir $N\subset D_k$ f"ur ein $k\in\{1,\ldots,n\}$ an, dann berandet jede andere Randkomponente $S_i, i\neq k$ eine 3-Disk $D_i\subset \overline{D_k\smallsetminus N}$. Wir k"onnen die $D_i$ einfach entlang der $S_i$ in $N$ einkleben ohne die Fundamentalgruppen von $N$ zu ver"andern und erhalten im Widerspruch dazu, dass $N$ nicht einfach-zusammen\-h"angend ist $\pi_1(N)=\pi_1(N\cup\bigcup_{i\neq k} D_i)=\pi_1(D_k)$.
	Es muss also $N\subset  M\smallsetminus \mathring D_i$ f"ur alle $i\in\{1,\ldots n\}$ gelten. Dann ist aber $N'\coloneqq N\cup\bigcup_{i=1}^nD_i$  eine geschlossene, kompakte 3-dimensionale Untermannigfaltigkeit von $M$. Das ist aber ebenfalls ein Widerspruch, da $M$ entweder als	kompakt mit Rand oder als nicht-kompakt vorausgesetzt war.
\end{beweis}

\begin{refbeweis}{Lemma~\ref{lm:Zerlegung-Mfkgruppe-unkomprimierbar}}
	Da $M$ irreduzibel und nicht einfach-zusammenh"angend ist, sind alle Randkomponenten von $M$ Fl"achen vom Geschlecht $\geq 1$. Sei also $Q$ eine Randkomponente von $M$ mit Geschlecht $\geq 1$, die komprimierbar in $M$ ist, das hei"st der Kern der Abbildung $f:\pi_1(Q)\rightarrow\pi_1(N)$ ist nicht trivial.

Sei $\gamma\in\ker f$.	Mit Korollar~\ref{kr:schleifensatzkorollar} erhalten wir eine eingebette 2-Disk $D\subset M$ mit $\gamma= \partial D \subset Q$. Wir trennen $M$ an dieser Disk auf und erhalten eine Untermannigfaltigkeit $N\subset M$ mit einer oder zwei Zusammenhangskomponenten mit Rand. Einfach-zusammenh"angende Zusammenhangskomponenten von $N$ brauchen wir f"ur die Zerlegung der Fundamentalgruppe von $M$ nicht weiter ber"ucksichtigen. Bei nicht einfach-zusammenh"angenden Zusammenhangskomponenten k"onnen wir $S^2$-L"ocher einfach mit 3-Disks stopfen, ohne dass sich die Fundamentalgruppe "andert. Mit Lemma~\ref{lm:Rand-von-Untermfk} wissen wir, dass dabei mindestens eine Randkomponente vom Geschlecht $\geq 1$ "ubrig bleibt. 
	
	In dem Fall, dass $N$ zwei Zusammenhangskomponenten  $N_1$ und $N_2$ besitzt, gilt mit dem Satz von Seifert und van Kampen $\pi_1(N)\cong\pi_1(N_1)*\pi_1(N_2)$. Im anderen Fall, wenn $N$ nur eine Zusammenhangskomponente hat, geht $M$ aus $N$ hervor, indem man zwei disjunkt in den Rand von $N$ eingebette 2-Disks verklebt. Es kommt in der Fundamentalgruppe lediglich ein freier Erzeuger dazu, also $\pi_1(M)=\pi_1(N)*\Z$.

	$N$ hat wie $M$ nur Randkomponenten vom Geschlecht $\geq 1$. 
	Sind darunter immer noch komprimierbare Randkomponenten, wiederholen wir den Prozess mit $N$. 
	Sollte $N$ nicht irreduzibel sein, zerlegen wir vorher mit Hilfe der Primzerlegung $\pi_1(N)$ in ein freies Produkt von Fundamentalgruppen von irreduziblen 3"=Mannigfaltigkeiten $N_i$ (Lemma~\ref{lm:Zerlegung-Mfkgruppe-in-freies-Produkt}). 
	Seien $D^1_i,\ldots,D^s_i$ die 3-Disks, an denen $N_i$ in der Primzerlegung von $N$ verklebt wird. Dann ist $N_i\smallsetminus(\mathring D^1_i\cup\ldots \mathring D^s_i)$ eine Untermannigfaltigkeit von $M$ und besitzt nach Lemma~\ref{lm:Rand-von-Untermfk} eine Randkomponente vom Geschlecht $\geq 1$. Also besitzen alle $N_i$ einen Rand, der sich dar"uber hinaus durch die Primzerlegung nicht ver"andert hat. Wir suchen uns in einem dieser R"ander einen neuen, in $M$ nullhomotopen Weg $\gamma$, der uns wieder eine komprimierbare 2-Disk liefert, an der wir weiter auftrennen k"onnen.

	Dieser Prozess muss irgendwann stoppen, wie man zum Beispiel mit Hilfe der Eulercharakteristik einsieht:
	Sei $U$ eine offene Umgebung von $\gamma$ in der Randfl"ache $Q$. Setze $X\coloneqq \overline{\partial M\smallsetminus U}$ und sei	$\iota:S^1\amalg S^1\rightarrow X$ die Inklusion, die jeweils eine $S^1$ mit einem durch das Auftrennen an $\gamma$ entstandenen Randkreis von $X$ identifiziert.
	Der "Ubergang von $M$ zu $N$ wird auf dem Rand dann durch die folgenden beiden Pushouts beschrieben.
	\begin{align*}
			\xymatrix{ 
			S^1\amalg S^1\ar@{->}[r]^{\iota}\ar@{->}[d]_{\id\amalg\id}
				&
				X\ar@{^{(}->}[d]
				\\D^1\amalg D^1\ar@{^{(}->}[r]
				&\partial M 
			}
		\qquad
	\xymatrix{ 
				S^1\amalg S^1\ar@{->}[r]^{\iota}\ar@{^{(}->}[d]
				&
				X\ar@{^{(}->}[d]
				\\D^1\amalg D^1\ar@{^{(}->}[r]
				&\partial N 
			}
	\end{align*}
	Im Fall, dass $M$ ein Volltorus ist, erhalten wir zum Beispiel folgendes Bild:
	
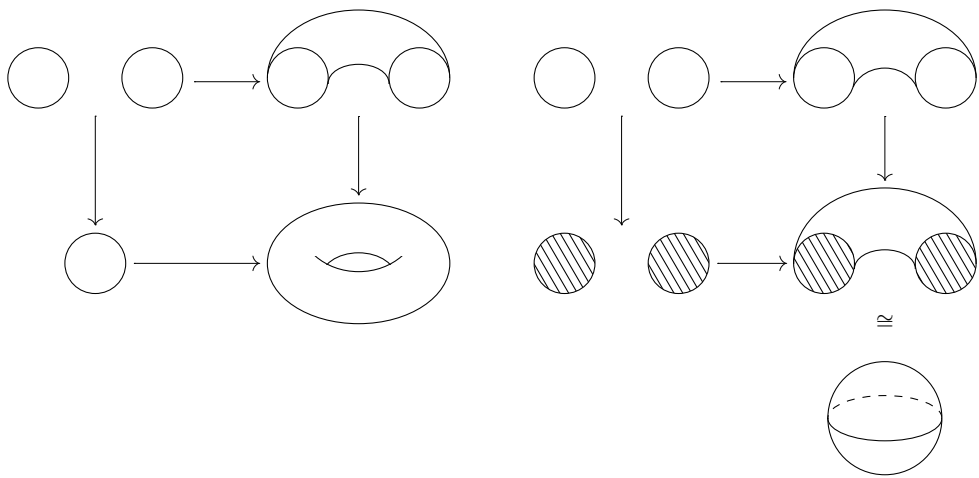
\begin{align*}
	\xymatrix{
			{\begin{tikzpicture}
				\clip (-1.2,-0.4) rectangle (1.2,1.3);
				\draw (-0.75,0) circle (0.4cm)
							(0.75,0) circle (0.4cm);
						\end{tikzpicture}}
			\ar@<-0.5cm>[r]
			\ar[d]
			&
			{
			\begin{tikzpicture}
				\draw (-0.8,0) circle (0.4cm);
				\draw (0.8,0) circle (0.4cm);
				\begin{scope}
					\clip (-1.2,0) rectangle (1.2,1.3);
					\draw (0,0) ellipse (1.2cm and 0.9cm);
				\end{scope}
				\begin{scope}
					\clip (-1.2,-0.07) rectangle (1.2,1.3);
					\draw (0,-0.07) ellipse (0.4cm and 0.25cm);
				\end{scope}
			\end{tikzpicture}}
			\ar[d]
			&
			\begin{tikzpicture}
				\clip (-1.2,-0.4) rectangle (1.2,1.3);
				\draw (-0.75,0) circle (0.4cm)
							(0.75,0) circle (0.4cm);
			\end{tikzpicture}
			\ar@<-0.5cm>[r]
			\ar[d]
			&
			\begin{tikzpicture}
				\draw(0.8,0) circle (0.4cm);
				\draw (-0.8,0) circle (0.4cm);
				\begin{scope}
					\clip (-1.2,0) rectangle (1.2,1.3);
					\draw (0,0) ellipse (1.2cm and 0.9cm);
				\end{scope}
				\begin{scope}
					\clip (-1.2,-0.11) rectangle (1.2,1.3);
					\draw (0,-0.27) ellipse (0.45cm and 0.4cm);
				\end{scope}	
			\end{tikzpicture}
			\ar[d]
			\\
				\begin{tikzpicture}
				\draw (0,0) circle (0.4cm);
			\end{tikzpicture}
			\ar[r]
			&
				\begin{tikzpicture}
				\draw (0,0) ellipse (1.2cm and 0.8cm);
				\begin{scope}
					\clip (-1,-0.01) rectangle (1,0.7);
					\draw (0,-0.21) ellipse (0.5cm and 0.35cm);
				\end{scope}
				\begin{scope}
					\clip (-1,0.09) rectangle (1,-0.8);
					\draw (0,0.59) ellipse (0.8cm and 0.7cm );
				\end{scope}
				\end{tikzpicture}
				&
			\begin{tikzpicture}
				\draw(0.75,0) circle (0.4cm);
				\draw (-0.75,0) circle (0.4cm);
				\begin{scope}
				\clip (-0.75,0) circle (0.4cm)
							(0.75,0) circle (0.4cm);
				\draw[step=1mm,rotate=30] (-2,1) grid (2,-1);
			\end{scope}
			\end{tikzpicture}
			\ar[r]
			&
			{\begin{tikzpicture}
				\clip (-1.2,-1) rectangle (1.2,1);
				\draw(0.8,0) circle (0.4cm);
				\draw (-0.8,0) circle (0.4cm);
				\begin{scope}
					\clip (-1.2,0) rectangle (1.2,1);
					\draw (0,0) ellipse (1.2cm and 1cm);
				\end{scope}
				\begin{scope}
					\clip (-1.2,-0.07) rectangle (1.2,1);
					\draw (0,-0.07) ellipse (0.4cm and 0.25cm);
				\end{scope}
					\begin{scope}
						\clip (-0.8,0) circle (0.4cm)
							(0.8,0) circle (0.4cm);
						\draw[step=1mm,rotate=30] (-2,1) grid (2,-1);
				\end{scope}
			\end{tikzpicture}}
			\\
			& & &
			\save[]+<0cm,0cm>*{
                           	\begin{tikzpicture}
			\draw (0,0) circle (0.75cm);
			\begin{scope}
				\clip (-0.75,-1) rectangle (0.75,0);
				\draw (0,0) ellipse (0.75cm and 0.3cm);
		\end{scope}
		\begin{scope}
				\clip (-0.75,1) rectangle (0.75,0);
				\draw[dashed] (0,0) ellipse (0.75cm and 0.3cm);
		\end{scope}	
	\end{tikzpicture}}\ar@{}[u]+<0cm,0.5cm>|{\cong}
		\restore
			}
		\end{align*}
		
	F"ur die Eulercharakteristik gilt damit
	\begin{eqnarray*}
		\chi(\partial M)& = &\chi(D^1\amalg D^1)+\chi(X)-\chi(S^1\amalg S^1) \\
										 & = &\chi(X) + 2 \\
										 & = &\chi(S^1\amalg S^1)+\chi(\partial N)-2 \chi(S^1) + 2 \\
										 & = &\chi(\partial N) + 2.
	\end{eqnarray*}
	Die Eulercharakteristik w"achst also mit jedem Schritt. Die einzige M"oglichkeit, die einem Ende noch im Wege steht, ist, dass als Rand unendlich viele Kopien von $S^2$ entstehen. Wir betrachten dazu eine einzelne Zusammenhangskomponente $K$ von $\partial M$, die wie oben beschrieben entlang eines einfach geschlossenen, nicht nullhomotopen Weges $\gamma$ aufgeschnitten wird. Zerf"allt $K$ in zwei Zusammenhangskomponenten $A$ und $B$, dann gilt $A\neq S^2\neq B$, denn sonst w"are $\gamma$ schon in $K$ nullhomotop gewesen.
	Es gilt also $\chi(A),\chi(B) \leq 0$	und mit der oben berechneten Beziehung $\chi(A)+\chi(B)=\chi(Q)+2$ folgt 
	\[
	\chi(K) < \chi(A),\chi(B) \leq 0.\]
	Jede Komponente $K$ von $\partial M$ kann also h"ochstens in $2^{|\chi(K)|}$-viele Zusammenhangskomponenten zerfallen, was also insgesamt endlich viele bleiben, da $\partial M$ kompakt ist.
\end{refbeweis} 

Abz"ahlbare freie Gruppen erf"ullen nach Lemma~\ref{lm:FICwF-freieGruppen} $FICwF$ und mit *W9 ist $FICwF$ auch mit dem freien Produkt vertr"aglich. Mit dem gerade bewiesenen Lemma~\ref{lm:Zerlegung-Mfkgruppe-unkomprimierbar} in Kombination mit Lemma~\ref{lm:mit-Rand=>nicht-kompakt} haben wir also f"ur berandete 3"=Mannigfaltigkeiten folgende Zwischenetappe erreicht:
\begin{etappe}\label{et:rand-etappe}
	Gilt $FICwF$ f"ur kompakte, irreduzible 3"=Mannigfaltigkeiten mit unkomprimierbarem Rand, dann gilt $FICwF$ f"ur 3"=Mannigfaltigkeiten mit Rand und f"ur nicht-kompakte 3"=Mannigfaltigkeiten.
\end{etappe}
\rem{
Zusammen mit unserem vorherigen Etappenziel~\ref{et:mit-rand-und-graph} ergibt sich daraus
\begin{etappe}
Gilt $FICwF$ f"ur nicht-triviale Graph-Mannigfaltigkeiten ohne Rand und f"ur kompakte, irreduzible 3"=Mannigfaltigkeiten mit unkomprimierbarem Rand, dann gilt $FICwF$ f"ur beliebige 3"=Mannigfaltigkeiten.
\end{etappe}}

Sei also $M$ eine kompakte, irreduzible 3"=Mannigfaltigkeit mit unkomprimierbarem Rand.
Da $M$ irreduzibel ist, sind alle Randkomponenten Fl"achen vom Geschlecht $\geq 1$. Besitzt $M$ eine Torusrandkomponente, greift folgendes Lemma: 

\begin{lemma}[{\cite[Korollar 4.1]{roushon}}]\label{lm:FICwF-torusrand}
	Sei $M$ eine kompakte, irreduzible 3"=Mannigfaltigkeit mit unkomprimierbarem Rand, der mindestens eine Torusrandkomponente besitzt. Dann wird $FICwF$ von $\pi_1(M)$ erf"ullt.
\end{lemma}
\begin{beweis}
	Sei $M'$ die Mannigfaltigkeit, die entsteht, wenn man $M$ mit einer Kopie von sich selbst an den Randfl"achen vom Geschlecht $\geq 2$ verklebt. Alle Randkomponenten von $M'$ sind also unkomprimierbare Tori.  
	
	Ist $M'$ atorisch k"onnen wir $M'$ entlang der unkomprimierbaren Randtori zu einer geschlossen Haken-Mannigfaltigkeit $M''$ verdoppeln. Die Toruszerlegung von $M''$ enth"alt damit eine atorische Komponente und mit Etappe~\ref{et:mit-rand-und-graph} wissen wir dann bereits, dass $FICwF$ gilt.

	Andernfalls enth"alt $M'$ einen unkomprimierbaren nicht-randparallelen Torus und ist damit eine Haken-Mannigfaltigkeit. Entweder enth"alt nun die Toruszerlegung von $M'$ ein atorisches St"uck und wir k"onnen analog zum vorherigen Fall $M'$ zu $M''$ verdoppeln und sind mit Etappe~\ref{et:mit-rand-und-graph} fertig oder $M'$ ist eine Graph-Mannigfaltigkeit. Da $M'$ einen Rand hat, k"onnen wir $M'$ in diesem Fall nach \cite[Satz 3.2]{leeb} mit einer nicht-positiven Metrik versehen, die auf dem Rand als Produkt gew"ahlt werden kann. Wir k"onnen somit $M'$ zu einer geschlossenen Mannigfaltigkeit mit nicht-positiver Metrik verdoppeln, so dass W2 greift.
   Wegen $\pi_1(M)\subset \pi_1(M')\subset \pi _1(M'')$ bekommen wir in allen F"allen mit Lemma~\ref{lm:FICwF-subgroup} $FICwF$ f"ur $\pi_1(M)$. 
\end{beweis}

Wir brauchen also, wie folgendes Diagramm verdeutlicht, nur noch kompakte, irreduzible 3"=Mannigfaltigkeiten mit unkomprimierbarem Rand, dessen Randkomponenten alle vom  Geschlecht $\geq 2$ sind, zu untersuchen.

\begin{align*}
  \xymatrix@!C=40pt{
 	&
		{\text{mit Rand}}
		\ar@{-}[dl]
		\ar@{-}[drr]
	& & &
	\\  
		{\txt<8pc>{kompri\-mierbare\\Rand\-komponenten}} 
	& & &
		{\txt<8pc>{nur unkompri\-mierbare\\Rand\-komponenten}}
		\ar@{-}[dr]
		\ar@{-}[dl]
		\ar@{=>}[lll]_{\text{Lemma~\ref{lm:Zerlegung-Mfkgruppe-unkomprimierbar}}}
	& &
	\\
	& &
		{\txt<6pc>{alle vom\\Geschlecht $\geq 2$}} 
	& &
		{\txt<6pc>{mindestens eine Toruskomponente}}  
	&  
	\\
	& &
		{\txt{Satz~\ref{thm:FICwF-B-Gruppen}}}
		\ar@{=>}[u]
	& &
		{\txt{Lemma~\ref{lm:FICwF-torusrand}}}
		\ar@{=>}[u]
	\\
	& & &
		{\txt{Etappe \ref{et:mit-rand-und-graph}}}
		\ar@{=>}[ul]
		\ar@{=>}[ur]
	& &
 }
\end{align*}

\subsection{B-Gruppen}
\begin{satz}\label{thm:FICwF-B-Gruppen}
  	Sei $M$ eine kompakte, irreduzible 3"=Mannigfaltigkeit mit unkomprimierbarem Rand, deren Randkomponenten alle Fl"achen vom Geschlecht $\geq 2$ sind.
	Dann ist $FICwF$ wahr f"ur $M$.
\end{satz}
In \cite{roushon} hat Roushon f"ur Fundamentalgruppen solcher Mannigfaltigkeiten die Bezeichnung \emph{B-Gruppe} eingef"uhrt und in \cite{roushon-b-gruppen} einen Beweis ausgef"uhrt, dass f"ur B-Gruppen $FICwF$ gilt. Wir wollen den Beweis im Folgenden nachvollziehen. \rem{TODO bzw.: einen alternative Beweis aufzeigen TODO }\rem{TODO Der Beweis besitzt allerdings eine Schwachstelle, an der er aber nicht zwingend bricht - zumindest konnte ich kein konkretes Gegenbeispiel konstruieren. Wir wollen den Beweis im Folgenden nachvollziehen und anschlie"send eine M"oglichkeit aufzeigen. 
}
Die Grundidee besteht darin $M$ mit einer Kopie zu einer geschlossenen Mannigfaltigkeit zu verkleben und mit Hilfe der Toruszerlegung f"ur die verdoppelte Mannigfaltigkeit $FICwF$ zu beweisen. Wir m"ussen dabei nur so geschickt verkleben, dass keine unkomprimierbaren Tori entstehen, die bei der Toruszerlegung die verklebten Randfl"achen von $M$ zerschneiden, damit wir folgendes Lemma anwenden k"onnen:
\begin{lemma}\label{lm:FICwF-Toruszerlegung-die-unkomprimierbare-Flaeche-nicht-zerschneidet}
  Sei $N$ eine geschlossene Haken-3"=Mannigfaltigkeit. Besitzt die Toruszerlegung von $N$ eine Komponente $K$, die eine unkomprimierbare, geschlossene Fl"ache $F$ vom Geschlecht $\geq 2$ enth"alt, dann wird $FICwF$ von $N$ erf"ullt.
\end{lemma}

\begin{beweis}
	Enth"alt die Toruszerlegung von $N$ eine atorische Komponente oder besteht sie nur aus einer Komponente, die eine Seifert-Mannigfaltigkeit ist, wissen wir bereits mit Etappe~\ref{et:mit-rand-und-graph}, dass $FICwF$ f"ur $N$ gilt. Wir m"ussen deshalb lediglich den Fall ausschlie"sen, dass $N$ eine nicht-triviale Toruszerlegung besitzt, die nur aus Seifert-Komponenten besteht, also dass $N$ eine nicht-triviale Graph-Mannigfaltigkeit ist.
   
	 Wir nehmen also an, $N$ ist eine nicht-triviale Graph-Mannigfaltigkeit. Dann ist $K$ eine Seifert-Komponente und besitzt, da die Toruszerlegung nicht trivial ist, einen nicht-leeren Torusrand. 
	 Die exakte Sequenz
	 \gexseq{\Z}{\pi_1(K)}{G}{i}{p}{,}
	 die wir schon einmal im Beweis von Satz~\ref{thm:FICwF-Seifert} betrachtet haben, liefert uns in diesem Fall, da $K$ einen Rand hat, eine freie Untergruppe $G_1\subset G$ mit endlichem Index in $G$ (siehe Bemerkungen zu Theorem 12.2 in \cite{hempel}).
	 $p\invers(G_1)$ hat dann endlichen Index in $\pi_1(K)$.
	 Betrachte die Konjugation 
	 \begin{eqnarray*}
	 		c:p\invers(G_1)&\rightarrow &\Aut(i(\Z)) \cong \Z/2\\
							g&\mapsto& (a\mapsto gag\invers)
						\end{eqnarray*}
	 $i(\Z)$ liegt im Zentrum von $\ker(c)$ und $R\coloneqq \ker(c)/i(\Z)$ ist als Untergruppe von $G_1$ frei.
		Es gilt damit $\ker(c)\cong\Z\times R$.	$F$ ist unkomprimierbar in $K$, es ist also $\pi_1(F)$ eine Untergruppe von $\pi_1(K)$.
	Da $\Z\times R$ endlichen Index in $\pi_1(K)$ hat, liegt	$\Z\times R \cap \pi_1(F)$ mit endlichen Index in $\pi_1(F)$. Es gibt somit eine endliche "Uberlagerung $\widetilde{F}$ von $F$ mit 
	\[
		\pi_1(\widetilde{F})=\Z\times R\cap\pi_1(F).
	\]
	$\widetilde{F}$ ist selbst auch wieder eine geschlossene Fl"ache vom Geschlecht $\geq 2$. Das ist aber ein Widerspruch, wie wir gleich sehen werden, denn als Untergruppe von $\Z\times R$ kann $\Z\times R\cap\pi_1(F)$ keine Fl"achengruppe vom Geschlecht $\geq 2$ sein.
 
	Sei
	\begin{eqnarray*}
	  f: \pi_1(\widetilde{F}) \rightarrow \Z\times R, \quad g \mapsto  (f_1(g),f_2(g)) 
	\end{eqnarray*}
	ein injektiver Gruppenhomomorphismus.
	Da Untergruppen freier Gruppen wieder frei sind, 
	$\pi_1(\widetilde{F})$ aber nicht frei ist, kann $f_2:\pi_1(\widetilde{F})\rightarrow R$ nicht injektiv sein.
	Da $f$ aber injektiv ist, gibt es ein $x\in \pi_1(\widetilde{F})$, $x\neq 1$ mit $x\in \ker(f_2)$ und $x\notin \ker(f_1)$. Also ist
	$x\in f\invers(\Z\times\{1\})$ und damit aus dem Zentrum $Z(\pi_1(\widetilde{F}))$. W"ahle ein Element $y\notin Z(\pi_1(\widetilde{F}))$. Da $x$ mit allen Elementen aus $\pi_1(\widetilde{F})$ kommutiert, gilt f"ur die von $x$ und $y$ erzeugte Gruppe
	\[
	\left< x,y \right> = \Z/m \times \Z/n\subset \pi_1(\widetilde{F}). 
	\]
	$\left< x,y \right>$ ist als Untergruppe von $\pi_1(\widetilde{F})$ die Fundamentalgruppe einer "Uberlagerung $T$ von $\widetilde{F}$ und damit selbst wieder Fl"achengruppe (ausf"uhrlicher bei \cite[4.3.7]{stillwell}). Die einzige Fl"ache mit abelscher Fundamentalgruppe mit zwei Erzeugern ist der Torus, es gilt also $n=m=1$ und da der Torus kompakt ist, ist $T$ eine endliche "Uberlagerung von $\widetilde{F}$ und damit $\Z\times\Z$ von endlichem Index in $\pi_1(\widetilde{F})$. Jetzt haben wir den gesuchten Widerspruch, denn Untergruppen von Fl"achengruppen vom Geschlecht $\geq 2$ mit endlichem Index sind selbst wieder Fl"achengruppen vom Geschlecht $\geq 2$ (\cite[Theorem 1]{hoare}).

	Den Fall, dass $N$ eine nicht-triviale Graph-Mannigfaltigkeit ist, haben wir somit ausgeschlossen und f"ur alle anderen F"alle hatten wir bereits gesehen, dass $FICwF$ gilt.
\end{beweis}

Bevor wir uns an den Beweis von Satz~\ref{thm:FICwF-B-Gruppen} machen, f"uhren wir zur Fallunterscheidung folgende Definition ein.
\begin{definition}
	Sei $A=S^1\times [0,1]$ ein Annulus.
  Wir nennen eine Mannigfaltigkeit $M$ \emph{geschlaucht}\footnote{Eine etwas freie, aber anschauliche "Ubersetzung von \emph{essential embedded annulus}. Es gibt allerdings auch Beispiele, bei denen die (zum Zylinder hom"oomorphe) Vorstellung einer eingebetteten Schallplatte n"aherliegender ist als die eines Schlauches.}, wenn es eine Einbettung $f:A\rightarrow M$ gibt, so dass gilt 
  \begin{enumerate}
	       \item $f(\partial A) \subset \partial M$,
				 \item $f$ ist zu keiner Einbettung $g:A\rightarrow M$ mit $g(A)\subset \partial M$ isotop relativ zum Rand und \rem{TODO} 
		   \item $f(A)$ ist in $M$ unkomprimierbar, das hei"st die von $f$ induzierte Abbildung $\pi_1(A)\rightarrow\pi_1(M)$ ist injektiv.
  \end{enumerate}
  Einen in $M$ eingebetteten Annulus, der obige Bedingungen erf"ullt, nennen wir\\ \emph{Schlauch}.
\end{definition}

\begin{refbeweis}{Satz~\ref{thm:FICwF-B-Gruppen}}
	\begin{fallumgebung}{1}
	\fall{$M$ ist nicht geschlaucht.} 
	Sei $N\coloneqq M\cup_{\partial M} M$. Da der Rand von $M$ unkomprimierbar ist, finden wir $\pi_1(M)$ als Untergruppe in $\pi_1(N)$. Mit Lemma~\ref{lm:FICwF-subgroup} reicht es also $FICwF$ f"ur $\pi_1(N)$ zu zeigen. 
	$\partial M$ ist in $N$ eine unkomprimierbare zweiseitige Fl"ache. Damit ist $N$ also eine Haken-Mannigfaltigkeit. 
	Die Toruszerlegung von $N$ liefert uns eine endliche Menge $\mathcal{T}$ von eingebetteten unkomprimierbaren Tori, die eindeutig sind bis auf Isotopie. Ist $\T$ leer, haben wir mit Etappe~\ref{et:mit-rand-und-graph} bereits $FICwF$ f"ur $M$ gezeigt.

	Wir wollen nun zeigen, dass es eine zu $\T$ isotope Menge von Tori $\T'$ gibt mit $\T'\cap\partial M=\emptyset$. Wir nehmen 
	an $T\in\T$ sei ein Torus, der mit keiner Isotopie von $\partial M$ getrennt werden kann. 
	Wir finden auf jeden Fall einen zu $T$ isotopen Torus $T'$, der $\partial M$ transversal schneidet. $C=\partial M\cap T'$ ist eine Menge eindimensionaler geschlossener Mannigfaltigkeiten (\cite[Satz 5.12]{broecker}), also einfach geschlossene Wege auf $\partial M$, im Folgenden als Kreise bezeichnet.
	Wir w"ahlen $T'$ so, dass die Anzahl der Kreise in $C$ minimal ist. Aufgrund folgender "Uberlegung k"onnen wir davon ausgehen,
	dass kein Kreis aus $C$ eine 2-Disk auf $T'$ berandet:
	
	Sei $\gamma$ ein Kreis auf $\partial M$, der eine 2-Disk $D_1$ auf $T'\smallsetminus \partial M$ berandet. $\gamma$ ist also in $M$ nullhomotop. Da $\partial M$ unkomprimierbar in $M$ liegt, ist $\gamma$ auch in $\partial M$ nullhomotop und berandet eine 2-Disk $D_2\subset\partial M$. $D_1\cup D_2$ ist eine 2-Sph"are, die eine 3-Disk berandet, da $M$ irreduzibel ist, das hei"st die 2-Sph"are kann mit einer Isotopie "uber den Rand geschoben werden und im Widerspruch zur Wahl von $T'$ haben wir einen zu $T'$ isotopen Torus mit weniger Kreisen gefunden. 
	
	Da nach unserer Annahme $\T'\cap\partial M\neq\emptyset$ gilt und da $N\smallsetminus \partial M$ nicht zusammenh"angend ist, gibt es in $C$ mehr als einen Kreis. Sei $A$ eine Komponente von $\overline{T'\smallsetminus C}$. Dann liegt $A$ vollst"andig in einer der Kopien von $M$. $\partial A\subset C$ ist in $T'$ nicht nullhomotop und damit finden wir $\pi_1(A)$ als nicht-triviale Untergruppe in $\pi_1(T')\cong\Z\times\Z$. Fl"achen mit Rand haben freie Fundamentalgruppen, \rem{[TODO Referenz?]} also bleibt lediglich die M"oglichkeit $\pi_1(A)=\Z$. $A$ ist somit ein Annulus, der wegen $\pi_1(A)\subset\pi_1(T')\subset\pi_1(M)$ unkomprimierbar in $M$ ist und da es nach Annahme keine Isotopie gibt, die $A$ "uber den Rand von $M$ schieben kann, ist $A$ auch nicht isotop relativ zum Rand zu $\partial M$.
Im Widerspruch zur Voraussetzung haben wir also einen Schlauch in $M$ gefunden.

Wir haben damit bewiesen, dass wir die Toruszerlegung von $N$ so w"ahlen k"onnen, dass dabei nicht der Rand von $M$ zerlegt wird. Es gibt also eine Komponente in der Toruszerlegung, die eine vollst"andige Randkomponente von $M$ enth"alt. Die Voraussetzung f"ur Lemma~\ref{lm:FICwF-Toruszerlegung-die-unkomprimierbare-Flaeche-nicht-zerschneidet} ist erf"ullt und es gilt $FICwF$ f"ur $\pi_1(N)$.
\\
\rem{
\begin{satz}
Sei $M$ eine kompakte 3"=Mannigfaltigkeit. Dann gibt es ein $N_0(M)\in\N$, so dass f"ur eine Menge $F\coloneqq\{F_1,\ldots,F_n\}$ von paarweise disjunkten, unkomprimierbaren in $M$ eingebetteten Fl"achen mit $n>N_0(M)$ gilt,
\begin{enumerate}
  \item ein $F_i$ ist ein Zylinder oder eine 2-Disk und randparallel oder 
  \item es gibt $i,j\in\{1,\ldots,n\}, i\neq j$, f"ur die $F_i$ parallel zu $F_j$ ist. 
\end{enumerate}

\end{satz}
Wir interessieren uns nur f"ur den Spezialfall, dass alle Fl"achen $F_i$ Zylinder sind. Der Beweis dazu findet sich mit Satz III.24 in \cite{jaco}. Wir formulieren ihn um, zu folgender Aussage:
}

\pagebreak
	\fall{$M$ ist geschlaucht.}
	Roushon ben"otigt f"ur seinen Beweis dieses Falles, die Voraussetzung, dass man in $M$ eine maximale Menge von disjunkten Schl"auchen findet, die eindeutig bis auf Isotopie sind. Die Existenz dieser Menge ist allerdings zweifelhaft. Wir gehen auf Roushons Beweis f"ur diesen Fall sp"ater genauer ein und zeigen stattdessen erstmal einen alternativen Beweis auf, der diese Voraussetzung umgeht, aber sonst "ahnlich funktioniert: Wir verdoppeln $M$ mit Hilfe von geeigneten pseudo-Anosov-Diffeomorphismen auf den Randfl"achen, so dass aus den Schl"auchen von $M$ in der Verdopplung keine unkomprimierbaren Tori entstehen und wir wieder Lemma~\ref{lm:FICwF-Toruszerlegung-die-unkomprimierbare-Flaeche-nicht-zerschneidet} anwenden k"onnen. An dieser Stelle m"ochte ich mich bei Thilo Kuessner bedanken, der mich auf den daf"ur entscheidenden Satz von Namazi und Souto aufmerksam gemacht hat.
\\

Sei $F$ eine unkomprimierbare Randfl"ache von $M$ vom Geschlecht $\geq 2$. Da $M$ Schl"auche enth"alt, ist $M$ insbesondere eine Haken-Mannigfaltigkeit. Sei $K$ die Komponente der Toruszerlegung von $M$ die $F$ enth"alt. 
Mit einer Randfl"achen vom Geschlecht $\geq 2$ kann $K$ keine Seifert-Mannigfaltigkeit sein, muss also atorisch sein.

Wir wollen auf $K$ folgenden Satz von Namazi und Souto anwenden:
\begin{satz}[{\cite[Theorem 5.1]{souto}}]\label{thm:suoto}
	Sei $N$ eine kompakte, irreduzible, orientierbare 3"=Mannigfaltigkeit, die keine $\Z\times\Z$-Untergruppe in ihrer Fundamentalgruppe hat und die eine Randfl"ache $S$ vom Geschlecht $\geq 2$ besitzt. Dann gibt es einen pseudo-Anosov-Diffeomorphismus $f$ auf $S$, so dass $\widetilde{N} = N \cup_f N$ eine negativ gekr"ummte Metrik besitzt.
\end{satz}

$K$ erf"ullt bereits alle Voraussetzung, au"ser wenn die Toruszerlegung von $M$ nicht trivial ist. Dann findet man durch die unkomprimierbaren Torusr"ander $\Z\times\Z$ in $\pi_1(K)$. Souto, darauf angesprochen, war davon "uberzeugt, dass der Beweis von Satz~\ref{thm:suoto} in diesem  Fall trotzdem durchgeht, da die $\Z\times\Z$-Untergruppen nur vom Rand kommen, aber wir k"onnen sie auch mit Hilfe von Dehn-F"ullungen loswerden:

Als atorische Mannigfaltigkeit besitzt $K$ eine geometrisch endliche, hyperbolische Metrik (siehe beispielsweise \cite[\S6]{scott}\rem{TODO zus"atzlich Thurston? TODO}). Mit Theorem 7.3 von Bromberg \cite{bromberg} k"onnen wir Volltori in die Torusr"ander einkleben, so dass die entstehende Mannigfaltigkeit immer noch hyperbolisch ist und damit $\Z\times\Z$ nicht mehr als Untergruppe enthalten kann (\cite[Lemma 4.5]{scott}). Jetzt liefert uns Satz~\ref{thm:suoto} einen Diffeomorphismus $f$ mit dem wir $K$ zu $\widetilde{K}\coloneqq K\cup_f K$ verdoppeln k"onnen, so dass $\widetilde{K}$ eine negativ gekr"ummte Metrik besitzt und somit atorisch ist (\cite[Korollar 9.9]{cheeger}). Mit Theorem 1.2.1 von Kojima~\cite{kojima} k"onnen wir aus  $\widetilde{K}$ die eingeklebten Volltori so wieder herausnehmen, dass  $\widetilde{K}$ weiterhin atorisch bleibt. $F$ kann also in $\widetilde{K}$ von keinem Torus zerlegt werden. 

Verdoppeln wir $M$ und nutzen dabei f"ur $F$ den Diffeomorphismus $f$, erhalten wir eine geschlossene Mannigfaltigkeit, die die Voraussetzung von Lemma~\ref{lm:FICwF-Toruszerlegung-die-unkomprimierbare-Flaeche-nicht-zerschneidet} erf"ullt.
\end{fallumgebung}
\end{refbeweis}

F"ur Roushons Beweis f"ur den 2. Fall m"ussen wir etwas genauer auf pseudo-Anosov-Diffeo\-mor\-phis\-men eingehen:

	Sei $S$ eine geschlossene orientierbare Fl"ache vom Geschlecht $\geq 2$ und sei $f:S\rightarrow S$ ein orientierungserhaltender Diffeomorphismus.	Nach Thurstons Klassifizierung von Diffeomorphismen auf Fl"achen (siehe beispielsweise \cite{casson} oder explizit \cite[\S V S.175]{poenaru}) besitzt $f$ eine der drei folgenden Eigenschaften:
	\begin{enumerate}
	  	\item Es gibt eine nat"urliche Zahl $n>0$, so dass $f^n$ isotop zur Identit"at ist.
	  	\item $f$ ist isotop zu einem reduzierbaren Diffeomorphismus. 
		\item $f$ ist isotop zu einem pseudo-Anosov-Diffeomorphismus,
	\end{enumerate} wobei sich sowohl 1) und 3) als auch 2) und 3) gegenseitig ausschlie"sen, d.h. ein pseudo-Anosov-Diffeomorphismus $f$ ist weder reduzierbar noch ist eine Potenz von $f$ isotop zur Identit"at.
  Ein Diffeomorphismus $f$ hei"st reduzierbar, wenn er eine endliche Menge $\A$ von disjunkten, paarweise nicht homotopen und nicht nullhomotopen, einfach geschlossenen Wegen festh"alt, also $f(a)\in\A$ f"ur alle $a\in\A$ gilt.
	
	Wir interessieren uns f"ur die 3. Kategorie. Dabei ist f"ur uns unerheblich, was ein pseudo-Anosov-Diffeomorphismus ist. Entscheidend ist, was er nicht ist, n"amlich reduzierbar und das gilt auch f"ur eine beliebige Potenz. Ein Beweis dazu findet zum Beispiel bei Po\'{e}naru in \cite[\S VI, Korollar von Proposition 17]{poenaru}.
	
	Mit Hilfe dieser Eigenschaft k"onnen wir folgendes Lemma beweisen.
	\begin{lemma}\label{lm:pseudo-Anosov-nichtparallel}
  Sei $S$ eine geschlossene orientierbare Fl"ache vom Geschlecht $\geq 2$ und $\A$  eine endliche Menge von disjunkten, paarweise nicht homotopen und nicht nullhomotopen, einfach geschlossenen Wegen. Sei $f:S\rightarrow S$ ein pseudo-Anosov-Diffeo\-mor\-phis\-mus. Dann gibt es ein $n_0\in\N$, so dass kein Element aus $\A$ isotop zu einem Element aus $f^{n_0}(\A)$ ist.
\end{lemma}
	\begin{beweis}
  Sei $\A=\{c_1,c_2,\ldots,c_k\}$ und $l\in\{1,\ldots,k\}$. F"ur ein beliebiges $n\in\N$ ist $f^n(c_l)$ nicht isotop zu $c_l$, da sonst $f^n$ reduzierbar w"are.
  Damit ist f"ur beliebige $n,m\in\N$ auch $f^n(c_l)$ nicht isotop zu $f^m(c_l)$, da sonst $f^{|n-m|}$ reduzierbar w"are. F"ur jedes $c\in\A$ gibt es also h"ochstens ein $n\in\N$, f"ur das $f^n(c_l)$ isotop zu $c$ ist. Da $\A$ endlich ist, gibt es ein $N_l\in\N$, so dass $f^s(c_l)$ f"ur alle $s>N_l$ zu keinem $c\in\A$ isotop ist. Mit $n_0=\max\{N_1,N_2,\ldots,N_k\}$ haben wir den gesuchten Exponenten von $f$ gefunden.
\end{beweis}

Entscheidend f"ur Roushons Beweis ist die folgende Behauptung, die er in  \cite[Remark 2.1.1]{roushon-b-gruppen} aufstellt. Er fordert darin nicht explizit die Disjunktheit, ben"otigt sie aber sp"ater im Beweis.

\begin{behauptung}\label{kr:maximale-schlauchmenge}
  Sei $M$ eine kompakte Mannigfaltigkeit. Dann gibt es eine maximale Menge $S=\{S_1,\ldots,S_N\}$ von disjunkten Schl"auchen in $M$, so dass jeder Schlauch $T$ in $M$ isotop zu einem Schlauch $S_i\in S$ ist.
\end{behauptung}

Roushon verweist f"ur den Beweis dieser Behauptung auf das Splitting Theorem von Jaco-Shalen (\cite{jaco-shalen}): 
\begin{satz}[Splitting Theorem]
	Sei $M$ eine kompakte, irreduzible Haken-Mannig\-faltig\-keit. Dann existiert eine 2-dimensionale unkomprimierbare Untermannigfaltigkeit $W\subset M$, eindeutig bis auf ambiente Isotopie, mit folgenden drei Eigenschaften:
	\begin{enumerate}
		\item Die Komponenten von $W$ sind Annuli und Tori und keine Komponente ist randparallel in $M$.
		\item Die Komponenten sind entweder Seifert-Mannigfaltigkeiten oder atorisch und nicht geschlaucht, d.h. jeder unkomprimierbare Torus oder Annulus in $W$ ist schon randparallel.
		\item $W$ ist minimal bez\"uglich aller Inklusionen von 2-Mannigfaltigkeiten in $M$, die die Eigenschaften (1) und (2) haben.
	\end{enumerate}
\end{satz}

Tats"achlich gibt es folgende, deutliche schw"achere Aussage "uber maximale Mengen von Schl"auchen, die auf Hakens Theorie der Normalfl"achen zur"uckgeht und auf der letztlich auch das Splitting Theorem beruht (siehe beispielsweise Satz III.24 von Jaco und Shalen (\cite{jaco}) von allgemeinen  Fl"achen auf den Spezialfall von Annuli eingeschr"ankt):
\begin{satz}\label{thm:maximale-schlaeuche}
	Sei $M$ eine kompakte 3"=Mannigfaltigkeit. Dann gibt es ein $N_0(M)\in\N$, so dass f"ur $n>N_0(M)$ in jeder Menge $F\coloneqq\{F_1,\ldots,F_n\}$ von paarweise disjunkten, unkomprimierbaren, nicht randparallelen, in $M$ eingebetteten Annuli bereits ein Annulus $F_i\in\F$ existiert, der zu einem anderen Annulus $F_j\in F, j\neq i$ parallel ist.
\end{satz}

Die Menge von Schl"auchen, die uns das Splitting Theorem liefert, k"onnen wir zwar mit Satz~\ref{thm:maximale-schlaeuche} zu einer maximalen Menge erweitern, aber im Allgemeinen nicht eindeutig bis auf Isotopie. Die Zerlegung entlang von Tori und Annuli ist gerade nur dann eindeutig, wenn man unter Ber"ucksichtigung der charakteristischen Untermannigfaltigkeit rechtzeitig mit dem zerlegen aufh"ort. In der charakteristischen Untermannigfaltigkeit findet man aber im Allgemeinen durchaus unterschiedliche, nicht disjunkte und nicht isotope Schl"auche.  

Gehen wir dennoch von der Existenz dieser maximalen Menge aus, k"onnen wir folgenden Beweis f"uhren:
\begin{proof}[Roushons Beweis des 2.Falls von Satz~\ref{thm:FICwF-B-Gruppen}]
Seien $\partial_1,\ldots,\partial_k$ die Komponenten von $\partial M$. Dann ist jedes der $\partial_i$ eine geschlossene, orientierbare Fl"ache vom Geschlecht $\geq 2$. Nach Voraussetzung besitzt $M$ mindestens einen Schlauch. 
	Wir w"ahlen nach Behauptung~\ref{kr:maximale-schlauchmenge} die maximale Menge $\mathcal{S}$ von Schl"auchen in $M$, die uns eine Menge $\A\coloneqq S\cap \partial M=\partial S$ von paarweise disjunkten und paarweise nicht homotopen, nicht nullhomotopen, einfach geschlossenen Wegen auf $\partial M$ liefert. 
	Sei $\A_i=\A\cap \partial_i$. Mit \rem{TODO gibt es auf jeder Fl"ache einen pseudo-Anosov-Diffeomorphismus TODO}Lemma~\ref{lm:pseudo-Anosov-nichtparallel} finden wir f"ur jedes $i=1,\ldots,k$ einen pseudo-Anosov-Diffeomorphismus  $f_i:\partial_i\rightarrow\partial_i$, so dass 
		$a\nsimeq b \mbox{ f"ur alle  } a\in\A \mbox{ und alle } b\in f(\A _i)$
 	gilt.	

	Seien $M^1$ und $M^2$ Kopien von $M$ und $N\coloneqq M^1\cup_{f_i}M^2$ die Verdopplung von $M$, die mit Hilfe der Diffeomorphismen $f_i$ verklebt wurde. $\partial M$  ist in $N$ eine zweiseitige, unkomprimierbare Fl"ache. $N$ ist also eine Haken-Mannigfaltigkeit und besitzt eine Toruszerlegung entlang von endlich vielen Tori $\T$. "Ahnlich wie in Fall A wollen wir wieder eine Isotopie von $\T$ finden, die $\partial M$ nicht schneidet.
	
	Wir nehmen an, es existiere ein Torus $T\in\T$, der sich durch keine Isotopie von $\partial M$ trennen l"asst. 
	Analog zu Fall A, finden wir eine Isotopie von $T'$ von $T$, so dass $T'\cap M^j$, $j=1,2$ nur aus Schl"auchen besteht.
	Da $\mathcal{S}$ bis auf Isotopie \rem{TODO vorher stand hier:Parallelit"at} maximal gew"ahlt war, ist $T\cap M^j$ isotop zu einer Teilmenge von $\mathcal{S}$ und damit liefert $T\cap \partial M^j$ f"ur $j=1,2$ Mengen $\A^j$ von Wegen die homotop zu Teilmengen von $\A$ sind.
	$T$ geht aus den Schl"auchen hervor, indem man die Randkomponenten $\A^1$ der Schl"auche $T\cap M^1$ mit den Randkomponenten $A^2$ der Schl"auche $T\cap M^2$ mit den Diffeomorphismen $f_i$ verklebt.
	F"ur $\gamma\in \A^1$ gilt also $f(\gamma)\in \A^2$, das hei"st eine Teilmenge von $\A$ ist isotop zu einer Teilmenge von $f_i(\A)$ ist, was im Widerspruch zur Wahl der $f_i$ steht. 
	
	Es gibt somit eine Isotopie $\T'$ von $\T$ mit $\T'\cap \partial=\emptyset$ und wir k"onnen Lemma~\ref{lm:FICwF-Toruszerlegung-die-unkomprimierbare-Flaeche-nicht-zerschneidet} anwenden. $\pi_1(N)$ erf"ullt damit $FICwF$. Da der Rand, entlang dem wir die Kopien von $M$ zu $N$ verklebt haben, unkomprimierbar ist, liegt $\pi_1(M)$ als Untergruppe in $\pi_1(N)$ und erf"ullt damit nach Lemma~\ref{lm:FICwF-subgroup} ebenfalls $FICwF$.
\end{proof}

Mit der Arbeit von Neumann und Swarup (\cite{neumann-swarup}) zur Toruszerlegung wird noch deutlicher, was das Problem dieses Beweises ist. Sie konstruieren eine Menge von unkomprimierbaren, disjunkten Annuli und Tori, die sie W-System nennen. Das W-System ist eine Obermenge der Tori und Annuli, die f"ur die Zerlegung nach Jaco-Shalen und Johannson notwendig sind. Es besitzt dar"uber hinaus die Eigenschaft, dass jeder weitere unkomprimierbare Torus oder Annulus sich mit einer Isotopie disjunkt von diesem W-System trennen l"asst (\cite[Lemma 2.2]{neumann-swarup}). Diese zus"atzlichen Annuli lassen sich  untereinander aber im Allgemeinen nicht disjunkt trennen. Wir haben in dem Beweis Tori gew"ahlt, die aus dem W-System von $N$ stammen, aber die Annuli, in die diese Tori in $M$ zerfallen, m"ussen nicht zum W-System von $M$ geh"oren und damit auch nicht disjunkt trennbar sein.\\ 

Mit den Ergebnissen aus diesem Kapitel nimmt das Baumdiagramm nun folgende Gestalt an:
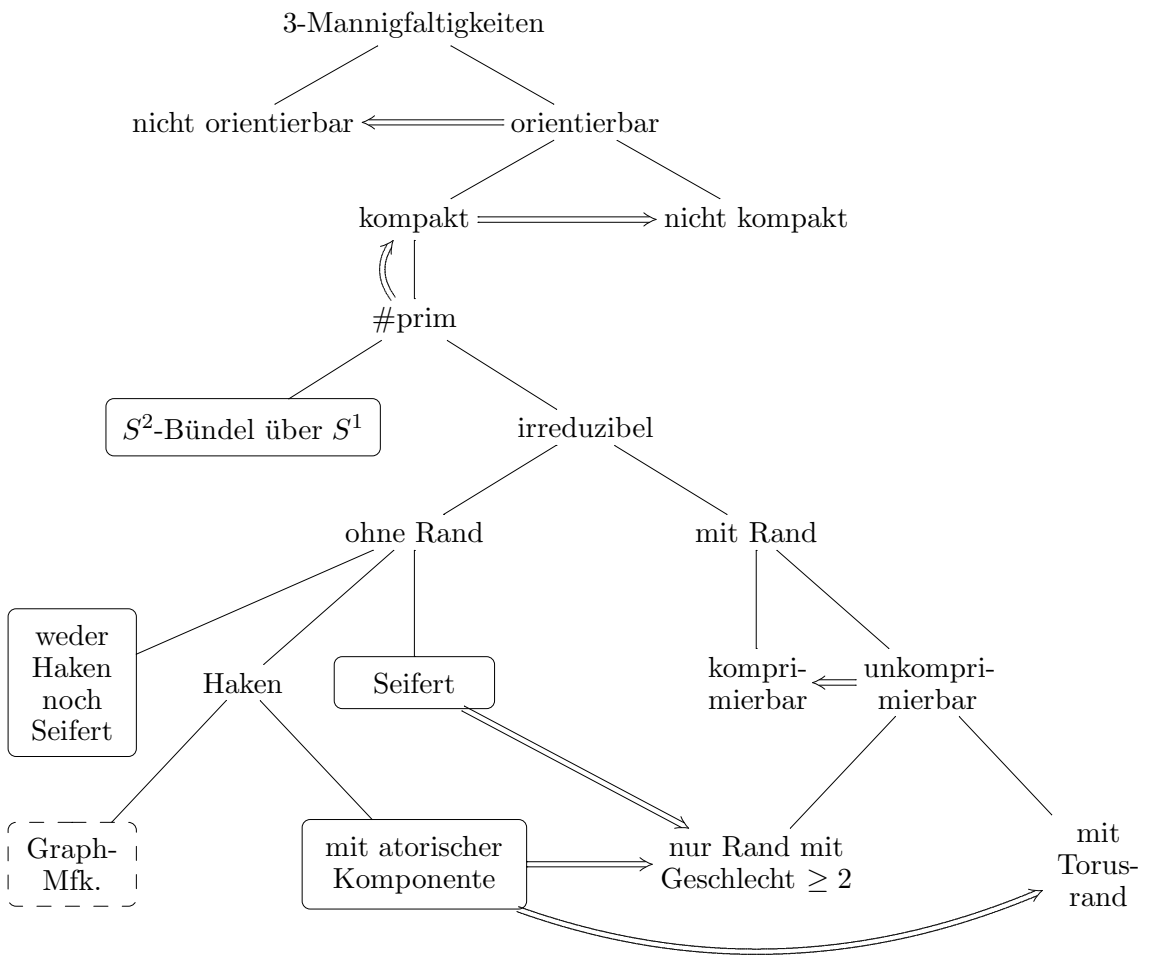
\begin{align*}
	\xymatrix@!C=40pt@R=22pt{
	& &  
		{\txt{3"=Mannigfaltigkeiten}}
		\ar@{-}[dl]
		\ar@{-}[dr]  
	& & & & 
	\\
    & 
		{\text{nicht orientierbar}} 
	& &  
		{\text{orientierbar}}
		\ar@{-}[dl] 
		\ar@{-}[dr]\ar@{=>}[ll] 
	& & &  
	\\
	& &  
		{\text{kompakt}} 
		\ar@{-}[d] 
		\ar@{=>}[rr]
		& & 
		{\text{nicht kompakt}} 
	& & 
	\\
	& &  
		{\text{\#prim}}  
		\ar@{=>}@/^1pc/[u]
		\ar@{-}[dl]
		\ar@{-}[dr]  
	& & & & 
	\\
	& 
	*++[F-:<3pt>]{\txt{$S^2$-B"undel "uber $S^1$}} 
	& & 
		{\text{irreduzibel}} 
		\ar@{-}[dl]
		\ar@{-}[dr] 
	& & &  
	\\ 
	& & 
		{\text{ohne Rand}} 
		\ar@{-}[dll]
		\ar@{-}[dl]
		\ar@{-}[d] 
	& & 
		{\text{mit Rand}}
		\ar@{-}[d] 
		\ar@{-}[dr]
	& & 
	\\
		*++[F-:<3pt>]{\txt<3pc>{weder Haken noch Seifert}} 
	& 
		{\txt{Haken}}
		\ar@{-}[dl]
		\ar@{-}[dr] 
	& 
		*++[F-:<3pt>]{\txt<4pc>{Seifert}} 
		\ar@{=>}[drr]
	& &  
		{\txt<3pc>{kompri\-mierbar}} 
	&
		{\txt<4pc>{unkompri\-mierbar}}
		\ar@{-}[dr]
		\ar@{-}[dl] 
		\ar@{=>}[l]
	\\
 		*++[F--:<3pt>]{\txt<3pc>{Graph-Mfk.}}
	& &
		*++[F-:<3pt>]{\txt<6pc>{mit atorischer Komponente}} 
		\ar@{=>}[rr]
		\ar@{=>} @/_7ex/[rrrr]
	& & 
		{\txt<6pc>{nur Rand mit\\Geschlecht $\geq 2$}} 
	& & 
		{\txt<3pc>{mit Torusrand}}  
	\\
	& & & & & &
	\\
	}
\end{align*}
\vfill\pagebreak

Mit den bereits erreichten Etappen~\ref{et:mit-rand-und-graph} und \ref{et:rand-etappe} k"onnen wir somit festhalten:
\begin{etappe}\label{et:FICwF-nichtkompakt-oder-mit-Rand}
  $FICwF$ ist wahr f"ur geschlossene Haken-Mannigfaltigkeiten, f"ur geschlossene Seifert-Mannigfaltigkeiten und f"ur kompakte, irreduzible 3"=Mannigfaltigkeiten mit Rand.
  Gilt $FICwF$ f"ur geschlossene Graph-Mannigfaltigkeiten, dann gilt $FICwF$ f"ur beliebige 3"=Mannigfaltigkeiten.
\end{etappe}

\section{Graph-Mannigfaltigkeiten-Ast}
Jetzt fehlt uns also nur noch $FICwF$ f"ur geschlossene nicht-triviale Graph-Mannig\-faltig\-keitsgruppen.
Der Schl"ussel dazu, ist der Hauptsatz aus \cite{roushon}. Mit Korollar~\ref{kr:FICwF-Bettizahl} erlaubt er uns
unseren Beweisbaum auf folgende Weise zu vervollst"andigen:
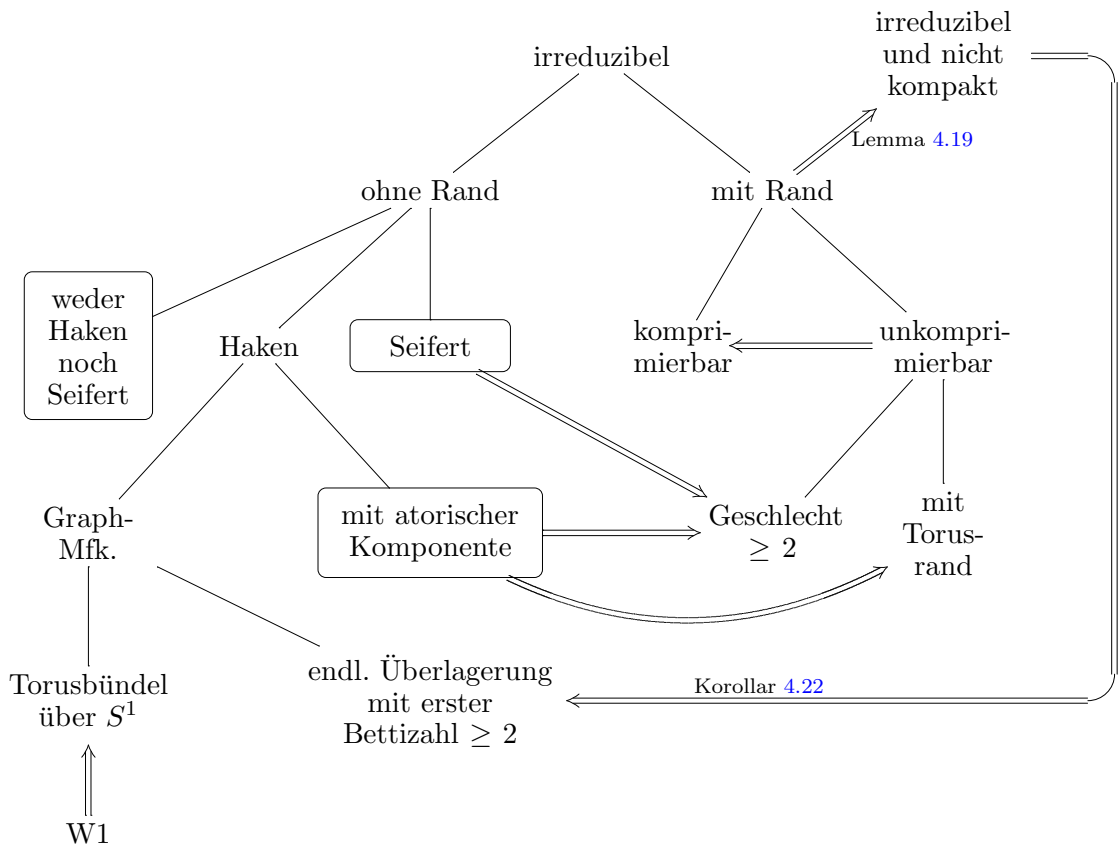
\begin{align*}
	\xymatrix@!C=40pt{
	\\
	& 
	& & 
		{\text{irreduzibel}} 
		\ar@{-}[dl]
		\ar@{-}[dr] 
	& & 
	{\txt<5pc>{irreduzibel und nicht kompakt}}
		\ar@{=>} `r[r] `[dddd] [ddddlll]_{\txt{\scriptsize Korollar~\ref{kr:FICwF-Bettizahl}}}
	&  
	\\ 
	& & 
		{\text{ohne Rand}} 
		\ar@{-}[dll]
		\ar@{-}[dl]
		\ar@{-}[d] 
	& & 
		{\text{mit Rand}}
		\ar@{-}[dr]
		\ar@{=>}[ru]_{\txt{\scriptsize Lemma~\ref{lm:FICwF-irreduzibel-nicht-kompakt}}}
	& & 
	\\
		*++[F-:<3pt>]{\txt<3pc>{weder Haken noch Seifert}}
	& 
		{\txt{Haken}}
		\ar@{-}[dl]
		\ar@{-}[dr] 
	& 
		*++[F-:<3pt>]{\txt<4pc>{Seifert}} 
		\ar@{=>}[drr]
	& 
	\save+<30pt,0pt>*{\txt<3pc>{kompri\-mierbar}}\ar@{<=}[rr]\ar@{-}[ur] \restore
	& & 
		{\txt<4pc>{unkompri\-mierbar}}
		\ar@{-}[d]
		\ar@{-}[dl] 
	\\
 		{\txt<4pc>{Graph-Mfk.}}
		\ar@{-}[d]
		\ar@{-}[drr]
	& &
		*++[F-:<3pt>]{\txt<6pc>{mit atorischer Komponente}} 
		\ar@{=>}[rr]
		\ar@{=>} @/_7ex/[rrr]
	& & 
		{\txt<4pc>{Geschlecht $\geq 2$}} 
	& 
		{\txt<3pc>{mit Torusrand}}
	&
	\\
			\txt<5pc>{Torusb"undel "uber $S^1$}
	\save+<0pt,-50pt>*+{
	\txt{W1}}
	\ar@{ =>}[u]+<0pt,-80pt>
	\restore
	& &  
		{\txt<8pc>{endl. "Uberlagerung mit erster Bettizahl $\geq 2$}}
	& & & & \\
	}
\end{align*}

\begin{lemma}\label{lm:FICwF-irreduzibel-nicht-kompakt}
	Sei $M$ eine irreduzible, nicht-kompakte 3"=Mannigfaltigkeit. Dann gilt $FICwF$ f"ur $M$.
\end{lemma}
\begin{beweis}
	Analog zu dem Beweis von Lemma~\ref{lm:mit-Rand=>nicht-kompakt} finden wir eine Folge von kompakten 3"=Mannigfaltigkeiten $M_i$ mit Rand, so dass
	$\pi_1(M) = \colim_{i\in I} M_i$ gilt.
	Mit Lemma~\ref{lm:Zerlegung-Mfkgruppe-in-freies-Produkt} k"onnen wir die $M_i$ in zusammenh"angende Summen zerlegen, so dass gilt
	\[\pi_1(M_i) \cong \pi_1(M^1_i)*\ldots*\pi_1(M^{k_i}_i)*F_i,\]

	wobei alle $M^j_i, 1\leq j \leq k_i$ irreduzibel und kompakt sind. 
	Wir nutzen ein Argument, dass wir schon im Beweis von Lemma~\ref{lm:Zerlegung-Mfkgruppe-unkomprimierbar} angewendet haben. In diesem Fall brauchen wir lediglich noch einen Index mehr:

	Seien $\sd{1}{D^j_i},\ldots\sd{n_{ji}}{D^j_i}$ die 3-Disks, an denen $M^j_i$ in der Primzerlegung von $M_i$ verklebt wird.
	Dann ist $M^j_i\smallsetminus(\sd{1}{\mathring D^j_i}\cup\ldots\cup\sd{n_{ji}}{\mathring D^j_i})$ eine Untermannigfaltigkeit von $M$ und besitzt nach Lemma~\ref{lm:Rand-von-Untermfk} mindestens eine Randkomponente vom Geschlecht $\geq 1$. Damit haben alle $M^j_i$ einen Rand.
	Mit Etappe~\ref{et:FICwF-nichtkompakt-oder-mit-Rand} wissen wir bereits, dass $FICwF$ f"ur irreduzible Mannigfaltigkeiten mir Rand gilt, also auch f"ur alle $M^j_i$ und damit auch f"ur alle $M_i$ und schlie"slich mit *W4 auch f"ur $M$.
\end{beweis}

\begin{satz}[{\cite[Main Theorem]{roushon}}]\label{thm:maintheorem}
	Sei $M$ eine geschlossene, irreduzible 3"=Mannigfaltigkeit. Gibt es eine Gruppe $H$ mit einem surjektiven Gruppenhomomorphismus 
	$p:\pi_1(M)\rightarrow H$ und folgenden Eigenschaften
	\begin{enumerate}
		\item $H$ besitzt eine nicht-triviale, torsionsfreie Untergruppe mit endlichem Index.
		\item $H$ erf"ullt $FICwF$.
		\item Jede unendliche, zyklische Untergruppe von $H$ hat unendlichen Index in $H$.
	\end{enumerate}
	Dann wird $FICwF$ von $\pi_1(M)$ erf"ullt.
\end{satz}
\begin{bemerkung}
	Da $H$ wegen der ersten Eigenschaft nicht endlich sein kann, ist die dritte Eigenschaft "aquivalent zu der Forderung, dass $H$ nicht virtuell zyklisch ist.
\end{bemerkung}
\begin{refbeweis}{Satz~\ref{thm:maintheorem}}
	F"ur eine Untergruppe $L$ von $\pi_1(M)$ bezeichnen wir im Folgenden mit $M_L$ die "Uberlagerung von $M$ mit Fundamentalgruppe $L$.   
	Hat $L$ unendlichen Index, dann ist $M_L$ eine nicht-kompakte 3"=Mannigfaltigkeit.\newline {}\newline
	\begin{fallumgebung}{1}
		\fall{H ist torsionsfrei.}
			Da $FICwF$ nach Voraussetzung von $H$ erf"ullt wird, reicht es aufgrund *W5 zu zeigen, dass f"ur jede virtuell zyklische Untergruppe $V$ von $H$ $FICwF$ f"ur $p\invers(V)$ wahr ist. Da $H$ torsionsfrei ist, ist auch $V$ torsionsfrei und kann deshalb als virtuell zyklische Gruppe nur trivial oder isomorph zu $\Z$ sein. 
    Da $H$ nicht virtuell zyklisch ist, muss $|H/V|$ unendlich sein. 
		Sei $pr:H\rightarrow H/V$ die Projektion. Betrachte folgende Mengenabbildung:
		\[
				\pi_1(M)\xrightarrow{p}H\xrightarrow{pr}H/V.
		\]
		$pr\circ p$ ist surjektiv und es gilt $\ker(pr\circ p)=p\invers(V)$.	$V$ ist zwar nicht normal in $H$, so dass uns der Homomorphiesatz zwar keine Isomorphie- daf"ur aber noch die M"achtigkeitsaussage $|\pi_1(M)/p\invers(V)|= |H/V|$ liefert.  
$M_{p\invers(V)}$ ist damit nicht-kompakt. Als "Uberlagerung einer irreduziblen Mannigfaltigkeit, ist $M_{p\invers(V)}$ au"serdem wieder irreduzibel (\cite[Theorem 7]{meeks}) und erf"ullt somit $FICwF$ aufgrund von Lemma~\ref{lm:FICwF-irreduzibel-nicht-kompakt}. F"ur jede beliebige virtuell zyklische Untergruppe $V$ von $H$ gilt also $FICwF$ f"ur $\pi_1(M_{p\invers(V)})=p\invers(V)$.\newline {}\newline 
		
		\fall{H ist nicht torsionsfrei.}
			Nach Voraussetzung existiert eine torsionsfreie Untergruppe $J$ von $H$ mit endlichem Index. Analog zum vorherigen Fall haben wir wieder eine surjektive Mengenabbildung 
			\[
				\pi_1(M)\xrightarrow{p}H\xrightarrow{pr}H/J,
			\]
			zu der uns der Homomorphiesatz
				\[
				|\pi_1(M)/p\invers (J)| = |H/J|
			\]
			liefert.
			Die "Uberlagerung $M_{p\invers(J)}$ von $M$ ist also endlich. 
			Aufgrund von Korollar~\ref{kr:FICwF-Ueberlagerung} brauchen wir $FICwF$ nur noch f"ur $\pi_1(M_{p\invers(J)})=p\invers(J)$ zu zeigen. Daf"ur ziehen wir uns einfach auf den vorherigen Fall zur"uck:
			
			Wir haben mit $p$ einen surjektiven Gruppenhomomorphismus von $p\invers(J)$ auf $J$. 
			Als torsionsfreie Gruppe ist $J$ selbst eine torsionsfreie Untergruppe von $J$ mit endlichem Index, als Untergruppe von $H$ erf"ullt $J$ mit Lemma~\ref{lm:FICwF-subgroup} $FICwF$ und weil $J$ endlichen Index in $H$ hat, "ubertr"agt sich die Eigenschaft von $H$, dass jede unendliche zyklische  Untergruppe unendlichen Index hat, auf $J$. Die Voraussetzungen des Satzes werden also von $J$ erf"ullt. Da $J$ torsionsfrei ist, k"onnen wir den ersten Fall dieses Beweises anwenden und wissen damit, dass $FICwF$ von $p\invers(J)$ erf"ullt wird. Der Hauptsatz ist damit bewiesen.
	\end{fallumgebung}
\end{refbeweis}

\begin{korollar}\label{kr:FICwF-Bettizahl}
	Sei $M$ eine geschlossene, irreduzible 3-Mannig\-faltig\-keit mit endlicher "Uberlagerung mit erster Bettizahl $\geq 2$. Dann gilt $FICwF$ f"ur $\pi_1(M)$.
\end{korollar}
\begin{beweis}
	Sei $N$ eine endliche "Uberlagerung von $M$ mit erster Bettizahl $\geq 2$.  
	Mit Korollar~\ref{kr:FICwF-Ueberlagerung} reicht es $FICwF$ f"ur $\pi_1(N)$ zu zeigen.

	Sei $H\coloneqq H_1(N;\Z)=\pi_1(N)/[\pi_1(N),\pi_1(N)]$. 
	Wir pr"ufen nach, dass $H$ alle Voraussetzungen des Hauptsatzes erf"ullt.
	$H$ ist eine endlich erzeugte abelsche Gruppe, hat also die Gestalt $H=\Z^r\oplus\Z/k_1\Z\oplus\ldots\oplus\Z/k_n\Z$.
	\begin{enumerate}
		\item	$\Z^r$ ist eine torsionsfreie Untergruppe von $H$ und $H/\Z^r\cong \Z/k_1\Z\oplus\ldots\oplus\Z/k_n\Z$ ist endlich.
		\item   $H$ als abelsche Gruppe erf"ullt nach Lemma~\ref{FICwF-virtuell-abelsch} $FICwF$  
		\item   Nach Voraussetzung ist $r\geq 2$ und damit hat jede unendliche zyklische Untergruppe von $H$ unendlichen Index in $H$.\hfill\qedsymbol
	\end{enumerate}
			\renewcommand{\qedsymbol}{}
\end{beweis}
\begin{satz}
	Sei $M$ eine geschlossene, irreduzible Graph-Mannig\-faltig\-keit. Dann ist $FICwF$ erf"ullt f"ur $\pi_1(M)$.
\end{satz}
\begin{beweis}
	Eine triviale Graph-Mannigfaltigkeit ist eine Seifert-Mannigfaltigkeit, f"ur die wir $FICwF$ bereits bewiesen haben (siehe Etappe~\ref{et:mit-rand-und-graph}). Sei $M$ also eine nicht-triviale geschlossene Graph-Mannigfaltigkeit. Dann enth"alt $M$ einen unkomprimierbaren Torus. Satz 1.1 aus \cite{luecke} sagt uns, dass wir f"ur $M$ eine endliche "Uberlagerung $N$ finden, die ein Torusb"undel "uber $S^1$ ist oder die eine erste Bettizahl $\geq 2$ hat. Gilt $FICwF$ f"ur $\pi_1(N)$, dann "ubertr"agt sich $FICwF$ mit Korollar~\ref{kr:FICwF-Ueberlagerung} auf $\pi_1(M)$.  Im letzten Fall folgt $FICwF$ f"ur $\pi_1(N)$ direkt aus Korollar~\ref{kr:FICwF-Bettizahl}.

  Bleibt der Fall, dass $N$ ein Torusb"undel "uber $S^1$ ist. F"ur das triviale Torusb"undel $(S^1\times S^1)\times S^1$ ist $\pi_1(N)\cong \Z\times\Z\times\Z$ abelsch und erf"ullt somit $FICwF$. Ein nicht triviales Torusb"undel "uber $S^1$ hat die Gestalt $(T^1\times [0,1])/\!\!\sim_f$, wobei $f:T^1\rightarrow T^1$ ein Hom"oomorphismus auf dem Torus ist und die "Aquivalenzrelation  $\sim_f$ identifiziert $(x,0)$ mit $(f(x),1)$ f"ur alle $x\in T^1$. F"ur die Fundamentalgruppe von $N$ bekommen wir dann die Pr"asentation 
	\[
	\pi_1(N)=\left< a,b,t; \;aba\invers b\invers, tat\invers = \pi_1(f)(a),tbt\invers = \pi_1(f)(b) \right>.
	\]
  Teilen wir $\left< a,b; \; aba\invers b\invers \right>$ aus $\pi_1(N)$ heraus, erhalten wir folgende exakte Sequenz
  \begin{align*}
  \xymatrix@!R=10pt{
  		1 
		\ar[r] 
	& 
		\left< a,b | aba\invers b\invers \right> 
		\ar[r] 
		\ar@{-}^{\cong}[d] 
	& 
		\pi_1(T\times [0,1])/\!\!\sim_f 
		\ar[r]
		\ar@{-}^{\cong}[d] 
	& 	
		\left< t\right> 
		\ar[r]
		\ar@{-}^{\cong}[d] 
	& 
		1 
	\\
  		1
		\ar[r] 
	& 
		\Z^2
		\ar[r] 
	& 
		\pi_1(N) 
		\ar[r] 
	& 
		\Z 
		\ar[r] 
	& 
		1 
	\\
  }
\end{align*}
Es gilt also $\pi_1(N)\cong \Z^2\sdp Z$, wof"ur wir mit W1 $FICwF$ gefordert hatten.
\end{beweis}

Damit hat unser FICwF-3"=Mannigfaltigkeiten-Baum folgende Gestalt angenommen

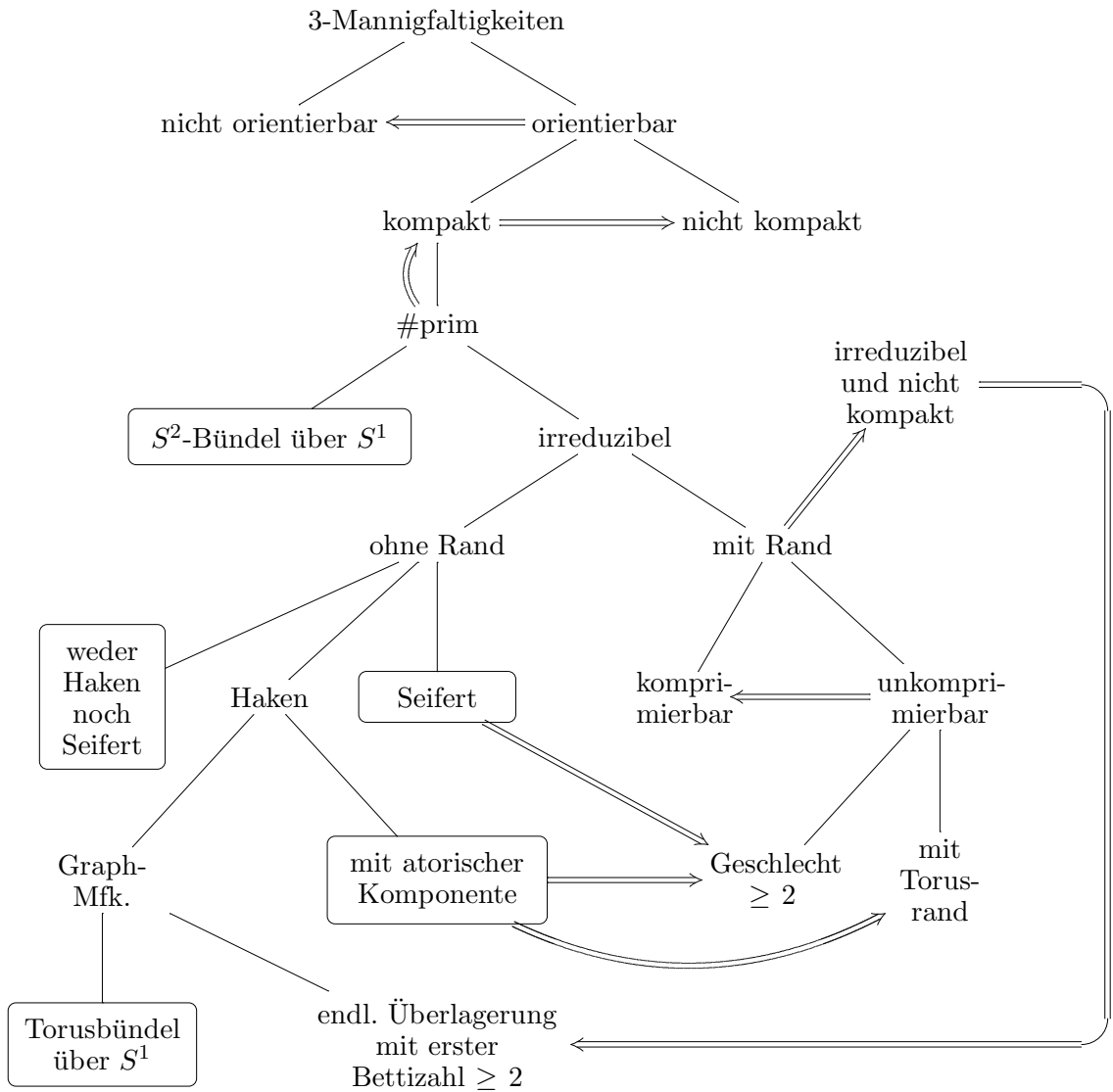
\begin{align*}
	\xymatrix@!C=40pt{
	& &  
		{\txt{3"=Mannigfaltigkeiten}}
		\ar@{-}[dl]
		\ar@{-}[dr]  
	& & & & 
	\\
    & 
		{\text{nicht orientierbar}} 
	& &  
		{\text{orientierbar}}
		\ar@{-}[dl] 
		\ar@{-}[dr]\ar@{=>}[ll] 
	& & &  
	\\
	& &  
		{\text{kompakt}} 
		\ar@{-}[d] 
		\ar@{=>}[rr]
		& & 
		{\text{nicht kompakt}}
	& & 
	\\
	& &  
		{\text{\#prim}}  
		\ar@{=>}@/^1pc/[u]
		\ar@{-}[dl]
		\ar@{-}[dr]  
	& & & & 
	\\
	& 
		*++[F-:<3pt>]{\txt<8pc>{$S^2$-B"undel "uber $S^1$}} 
	& & 
		{\text{irreduzibel}} 
		\ar@{-}[dl]
		\ar@{-}[dr] 
	& &
	\save+<-15pt,20pt>*{\txt<5pc>{irreduzibel und nicht kompakt}}
		\ar@{=>} `r[r] `[dddd] [ddddlll]
		\ar@{<=}[dl]
		\restore
&  
	\\ 
	& & 
		{\text{ohne Rand}} 
		\ar@{-}[dll]
		\ar@{-}[dl]
		\ar@{-}[d] 
	& & 
		{\text{mit Rand}}
		\ar@{-}[dr]
	& & 
	\\
		*++[F-:<3pt>]{\txt<3pc>{weder Haken noch Seifert}}
	& 
		{\txt{Haken}}
		\ar@{-}[dl]
		\ar@{-}[dr] 
	& 
		*++[F-:<3pt>]{\txt<4pc>{Seifert}} 
		\ar@{=>}[drr]
	& 
	\save+<30pt,0pt>*{\txt<3pc>{kompri\-mierbar}}\ar@{<=}[rr]\ar@{-}[ur] \restore
	& & 
		{\txt<4pc>{unkompri\-mierbar}}
		\ar@{-}[d]
		\ar@{-}[dl] 
	\\
 		{\txt<4pc>{Graph-Mfk.}}
		\ar@{-}[d]
		\ar@{-}[drr]
	& &
		*++[F-:<3pt>]{\txt<6pc>{mit atorischer Komponente}} 
		\ar@{=>}[rr]
		\ar@{=>} @/_7ex/[rrr]
	& & 
		{\txt<4pc>{Geschlecht $\geq 2$}} 
	& 
		{\txt<3pc>{mit Torusrand}}
	&
	\\
		*++[F-:<3pt>]{\txt<5pc>{Torusb"undel "uber $S^1$}} 
	& &  
		{\txt<8pc>{endl. "Uberlagerung mit erster Bettizahl $\geq 2$}}
	& & & & 
	}
\end{align*}
\\
\\
\\
\\
\\
\\
und wir formulieren unseren 
\begin{schlusssatz}
	Sei $FICwF$ eine endlich erweiterbare, gefaserte Isomorphismusvermutung mit Werkzeugkasten.
  Dann gilt $FICwF$ f"ur Fundamentalgruppen von 3-Mannig\-faltig\-keiten.

	Gilt $K$-$FICwF$ f"ur Gruppen der Gestalt $\Z^2\sdp\Z$ und $K$-$FIC$ f"ur $n\geq1$ f"ur Gruppen, die proper, kokompakt und isometrisch auf einem $\CAT(0)$-Raum operieren, dann gilt $K$-$FICwF$ f"ur Fundamentalgruppen von 3-Mannigfaltigkeiten.

	Gilt $L$-$FICwF$ f"ur Gruppen der Gestalt $\Z^2\sdp\Z$ und $L$-$FIC$ f"ur virtuell abelsche Gruppen, dann gilt $L$-$FICwF$ f"ur Fundamentalgruppen von 3-Mannigfaltigkeiten.
\end{schlusssatz}

\begin{appendix}
	\chapter{Werkzeug"ubersicht}
	Da wir sie h"aufig brauchen, geben wir zur "Ubersicht im Folgenden nochmal alle verwendeten Werkzeuge an:
	\begin{enumerate}
	\item[W1:] $FICwF$ gilt f"ur alle Gruppen der Gestalt $\Z^2\sdp\Z$.
	\item[W2:] $FICwF$ gilt f"ur Fundamentalgruppen von nicht-positiv ge\-kr"umm\-ten, geschlossenen Mannigfaltigkeiten. 
	\item[W3:] $FIC$ gilt f"ur virtuell abelsche Gruppen.
	\item[W4:] Sei $\{G_i\; | \; i\in I\}$ ein gerichtetes System von Gruppen und f"ur jedes $G_i$ gelte $FIC$. Dann gilt $FIC$ auch f"ur $\colim_{i\in I}G_i$. 
	\item[W5:] Sei $p:L\rightarrow G$ ein Gruppenhomomorphismus und $FIC$ gelte f"ur $G$. Ist $FIC$ wahr f"ur $p\invers(V)$ f"ur jede virtuell zyklische Untergruppe $V$ von $G$, dann gilt $FIC$ auch f"ur $L$.
\end{enumerate}

\begin{enumerate}
	\item[*W3:] $FICwF$ gilt f"ur virtuell abelsche Gruppen.
	\item[*W4:] Sei $\{G_i\; | \; i\in I\}$ ein gerichtetes System von Gruppen und f"ur jedes $G_i$ gelte $FICwF$. Dann gilt $FICwF$ auch f"ur $\colim_{i\in I}G_i$. 
	\item[*W5:] Sei $p:L\rightarrow G$ ein Gruppenhomomorphismus und $FICwF$ gelte f"ur $G$. Ist $FICwF$ wahr f"ur $p\invers(V)$ f"ur jede virtuell zyklische Untergruppe $V$ von $G$, dann gilt $FICwF$ auch f"ur $L$.
\end{enumerate}

\begin{enumerate}
	\item[*W6:] $FICwF$ gilt f"ur endliche Gruppen.
  \item[*W7:] $FICwF$ gilt f"ur abz"ahlbare, freie Gruppen.	
	\item[*W8:] Seien $G_1$ und $G_2$ Gruppen die $FICwF$ erf"ullen. Dann gilt $FICwF$ auch f"ur $G_1\times G_2$.
	\item[*W9:] Seien $G_1$ und $G_2$ abz"ahlbare Gruppen die $FICwF$ erf"ullen. Dann gilt $FICwF$ auch f"ur $G_1*G_2$.
\end{enumerate}

\end{appendix}

\end{document}